\documentclass[review,onefignum,onetabnum]{siamart251216}

\usepackage{lipsum}
\usepackage{amsfonts}
\usepackage{graphicx}
\usepackage{epstopdf}
\ifpdf
  \DeclareGraphicsExtensions{.eps,.pdf,.png,.jpg}
\else
  \DeclareGraphicsExtensions{.eps}
\fi

\newsiamremark{remark}{Remark}

\headers{Neural Variational Method for the Anisotropic Eikonal Equation}{Y. Park and J. Hahn}

\title{A Unit-Gradient Neural Variational Method for Anisotropic Eikonal Equations\thanks{Submitted to the editors DATE.
\funding{This document is the supported by the research project funded by NRF Grant RS-2024-00343226.}}}

\author{Yesom Park\thanks{Department of Mathematics, University of California, Los Angeles, Los Angeles, CA 90095, USA 
  (\email{yeisom@math.ucla.edu}).}
\and Jooyoung Hahn\thanks{Department of Software Engineering, Faculty of Nuclear Sciences and Physical Engineering, Czech Technical University in Prague, Trojanova 339/13, Prague 120 00, Czech Republic
  (\email{jooyoung.hahn@fjfi.cvut.cz}) Corresponding author.}}

\usepackage{amsopn}

\ifpdf
\hypersetup{
  pdftitle={Variational Neural Method for Anisotropic Eikonal Equations},
  pdfauthor={Y. Park and J. Hahn}
}
\fi
\usepackage{wrapfig}
\usepackage{mathtools}
\usepackage{amssymb}
\usepackage{booktabs}
\usepackage{subfigure}
\usepackage{multirow}
\usepackage{kotex}
\usepackage{algorithm}
\usepackage{algpseudocode}
\usepackage{tikz}
\usetikzlibrary{shapes.geometric, arrows}
\usetikzlibrary{positioning}  
\usetikzlibrary{calc}
\usetikzlibrary{fit}
\usetikzlibrary{decorations.pathreplacing}
\usetikzlibrary{arrows,shapes,backgrounds}
\usetikzlibrary{arrows.meta}

\newtheorem{example}{Example}[section]

\makeatletter
\newcommand*{\addFileDependency}[1]{\typeout{(#1)}\@addtofilelist{#1}\IfFileExists{#1}{}{\typeout{No file #1.}}}
\makeatother
\addFileDependency{SISC2026_park_hahn_supplement.tex}
\addFileDependency{SISC2026_park_hahn_supplement.aux}

\newcommand{\cH}{\mathcal{H}}

\newcommand{\cL}{\mathcal{L}}

\newcommand{\cU}{\mathcal{U}}

\newcommand{\bR}{\mathbb{R}}

\newcommand{\dV}[1][u]{\cU\left(\Omega;\bR^{d_v}\right)}
\newcommand{\dVt}[1][u]{\cU\left(\left[0,\infty\right)\times\Omega;\bR^{d_v}\right)}

\newcommand{\bx}{\mathbf{x}}
\newcommand{\by}{\mathbf{y}}

\newcommand{\bp}{\mathbf{p}}

\newcommand{\pnet}{\mathbf{p}_{\theta'}}

\DeclareMathOperator*{\argmin}{\arg\min}

\newcommand*\diff{\mathop{}\!\mathrm{d}}

\definecolor{forestgreen}{RGB}{34,159,44}
\definecolor{deepred}{RGB}{235,0,20} 
\definecolor{myred}{HTML}{DC133B}
\definecolor{myblue}{HTML}{1F77B4}
\definecolor{myblack}{HTML}{000000}
\definecolor{mygray}{RGB}{125,125,125}

\begin{document}

\maketitle

\begin{abstract}
We develop a variational neural method for computing viscosity solutions of the anisotropic eikonal equation. Building on the vanishing-viscosity principle, we minimize an energy functional whose Euler--Lagrange equation is the viscous anisotropic eikonal equation and introduce an auxiliary variable for the gradient, enforcing the constraint through a penalty.
A key challenge in this formulation is to handle nonlinearity of the functional while controlling the interaction between the viscosity and penalty parameters. To address this challenge, we impose a unit-norm constraint on the auxiliary variable exactly, as prescribed by the eikonal equation.
This removes the gradient from the exponentially weighted part of the functional and allows the viscosity and penalty limits to be studied separately. We establish two convergence results. First, for each fixed viscosity parameter $\epsilon>0$, we prove convergence to the viscosity solution as the penalty parameter $\lambda\rightarrow\infty$. Second, for the distance function to the boundary of a convex or smooth domain, we prove convergence as $\epsilon\rightarrow0$ for any fixed $\lambda>0$. 
The latter result motivates our stable numerical strategy, which keeps $\lambda$ fixed while decreasing $\epsilon$ whenever the training loss stagnates. The formulation is implemented using a neural network that satisfies the unit-norm constraint by construction and a mesh-free approach that evaluates the singularly scaled functional at randomly sampled points without requiring mesh refinement as $\epsilon\to0$. 
Numerical experiments demonstrate the effectiveness and robustness of the method on diverse isotropic and anisotropic problems, including three-dimensional open-surface reconstruction from point clouds.
\end{abstract}

\begin{keywords}
Anisotropic eikonal equation; Viscosity solution; Vanishing viscosity; Variational method; Penalty method; Scientific machine learning
\end{keywords}

\begin{MSCcodes}
35F21, 35D40, 35B25, 65K10, 68T07
\end{MSCcodes}
\section{Introduction}
We consider the anisotropic eikonal equation
\begin{equation}\label{eq:eikonal}
\begin{cases}
  \left\Vert \nabla u \right\Vert_M = 1
  & \text{in } \Omega\setminus\Gamma,\\
  u = g
  & \text{on } \Gamma.
\end{cases}
\end{equation}
where $\Omega\subset\bR^d$ is a bounded domain, $\Gamma\subset\overline{\Omega}$ is compact, $g:\Gamma\to\bR$ is nonnegative, and $\|\by\|_M=\sqrt{\by^{\mathrm T}M\by}$ with a spatially varying Riemannian metric $M^{-1}\in\mathcal C^0(\overline{\Omega},\mathcal S_d^+)$, $\mathcal S_d^+$ being the symmetric positive definite $d\times d$ matrices. We assume that $g$ is compatible with the metric, $|g(\bx)-g(\by)|\le d_M(\bx,\by)$ for $\bx,\by\in\Gamma$, where $d_M$ is the geodesic distance in $\overline\Omega$ of the metric $M^{-1}$. For $\Gamma=\partial\Omega$, this compatibility condition is necessary and sufficient for the existence of a viscosity solution satisfying the boundary condition~\cite{lions1982generalized}. Equation~\eqref{eq:eikonal} models first-arrival travel times in anisotropic media~\cite{alkhalifah2002traveltime} and geodesic distances in imaging~\cite{peyre2010geodesic}. 
The equation admits multiple weak solutions, namely, Lipschitz functions that satisfy~\eqref{eq:eikonal} almost everywhere and the boundary condition $u=g$ on $\Gamma$~\cite{tran2021hamilton}.
Among them, the physically meaningful one is the viscosity solution $u^*$, the limit as $\epsilon\to0$ of the solutions of the viscous equation obtained by adding the diffusion term $-\epsilon\nabla\cdot(M\nabla u)$~\cite{varadhan1967behavior,lions1982generalized,crandall1983viscosity,bonnans2022linear}. We compute $u^*$ through this viscous approximation.

Classical numerical methods for the eikonal equation include fast marching and ordered upwind methods~\cite{sethian1996fast,sethian2001ordered,mirebeau2014anisotropic}, fast sweeping methods~\cite{zhao2005fast,qian2007fast}, semi-Lagrangian schemes~\cite{cristiani2007fast}, and vanishing-viscosity schemes~\cite{churbanov2019numerical,yang2022continuous,hahn2024laplacian}. These methods achieve high accuracy and efficiency on appropriately constructed meshes. However, mesh generation remains a major bottleneck in industrial simulations, particularly for complex geometries, and fully automated meshing and mesh adaptation remain challenging research objectives~\cite{slotnick2014cfd}. Anisotropic mesh adaptation, in particular, continues to receive considerable attention~\cite{alauzet2016decade}. Moreover, on polyhedral meshes, which are commonly used for complex geometries~\cite{garimella2014polyhedral}, the upwind updates underlying classical eikonal solvers, typically derived for simplicial meshes, cannot be directly applied~\cite{hahn2024laplacian}.  
These challenges motivate mesh-free approaches that avoid explicit mesh generation.
For surface reconstruction from the unorganized point clouds, the complexity lies in the data $\Gamma$ rather than $\Omega$, and Cartesian-grid level-set methods also handle such data without unstructured mesh generation~\cite{zhao2000implicit}; extending our approach to domains that require unstructured meshes remains a direction for future research.

Physics-informed neural networks (PINNs)~\cite{raissi2019physics} minimize the PDE residual and have been applied to isotropic~\cite{gropp2020implicit,bin2021pinneik} and anisotropic~\cite{waheed2020anisotropic} eikonal equations, but a small residual of~\eqref{eq:eikonal} does not single out the viscosity solution~\cite{park2025nsp,han5352686stochastic,kim2027115466neural}. The Deep Ritz method~\cite{e2018deep}, primarily developed for elliptic PDEs, minimizes variational energies over neural networks with integrals estimated at random points, as we do for a functional derived from the viscous eikonal equation. Least-squares residuals of monotone finite-difference schemes select $u^*$ through the numerical diffusion of the scheme, with an accuracy set by its step size~\cite{esteve2025finite,bokanowski2026solving}. Another approach employs a new implicit formulation to obtain viscosity solutions of time-dependent Hamilton--Jacobi equations, with a primary focus on high-dimensional problems~\cite{park2025neural}. Other neural solvers target $u^*$ through a stochastic reading of the vanishing-viscosity limit~\cite{han5352686stochastic} or, for distance functions ($M=I$, $g=0$), through viscosity-type and $p$-Poisson formulations~\cite{lipman2021phase,park2024p,wang2024hotspot} and an augmented Lagrangian method (ALM) for a maximal-subsolution characterization~\cite{kim2027115466neural}. Closest to our work, HotSpot~\cite{wang2024hotspot} derives a variational problem from the viscous eikonal equation by the Hopf--Cole transformation, only for $M=I$ and $g=0$, and adds the additional eikonal residual as a soft penalty.

\begin{figure}[t]
\centering
\begin{tikzpicture}[
    node distance=0.3cm and 0.6cm,
    block/.style={draw=black!60, rounded corners=2pt, minimum width=3.0cm,
        minimum height=0.7cm, align=center, font=\scriptsize, fill=black!3},
    keyblock/.style={block, draw=teal!70!black, fill=teal!5, thick},
    arrow/.style={->, >=stealth, thick, black!60},
    outcome/.style={font=\tiny, text=black!70, align=center, text width=3.0cm}
]
\node[block] (a) {\textbf{Anisotropic eikonal}\\ $\|\nabla u\|_M=1$ \eqref{eq:eikonal}};
\node[block, right=of a] (b) {\textbf{Vanishing viscosity}\\ $-\epsilon\nabla\cdot(M\nabla u)$ \eqref{eq:viscosity_eikonal}};
\node[block, right=of b] (c) {\textbf{Variational formulation}\\ $\mathcal J_\epsilon(u)$ \eqref{eq:direct_ftnal}};
\node[block, below=0.7cm of c] (d) {\textbf{Variable splitting}\\ $\mathbf p=\nabla u$ \eqref{eq:constraint_prob}};
\node[keyblock, left=of d] (e) {\textbf{Unit-norm constraint}\\ $\|\mathbf p\|_M=1$ \eqref{eq:proposed_formulation}};
\node[block, left=of e] (f) {\textbf{Deep learning}\\ $\mathbf p=\hat{\mathbf p}/\|\hat{\mathbf p}\|_M$ \eqref{eq:eik_auxiliary_variable}};
\draw[arrow] (a) -- (b);
\draw[arrow] (b) -- (c);
\draw[arrow] (c.east) -- ++(0.25,0) |- (d.east);
\draw[arrow] (d) -- (e);
\draw[arrow] (e) -- (f);
\node[outcome, below=0.08cm of a] {Many weak solutions};
\node[outcome, below=0.08cm of b] {Viscosity solution as $\epsilon\to0$};
\node[outcome, below=0.08cm of c] {No second derivatives};
\node[outcome, below=0.08cm of d] {Gradient convergence};
\node[outcome, below=0.08cm of e] {Decoupling of limits};
\node[outcome, below=0.08cm of f] {Exact constraint, mesh-free};
\end{tikzpicture}
\caption{Ingredients of the method (boxes, in the order of the arrows) and, below each box, the difficulty it poses or what it provides.}
\label{fig:method_overview}
\end{figure}
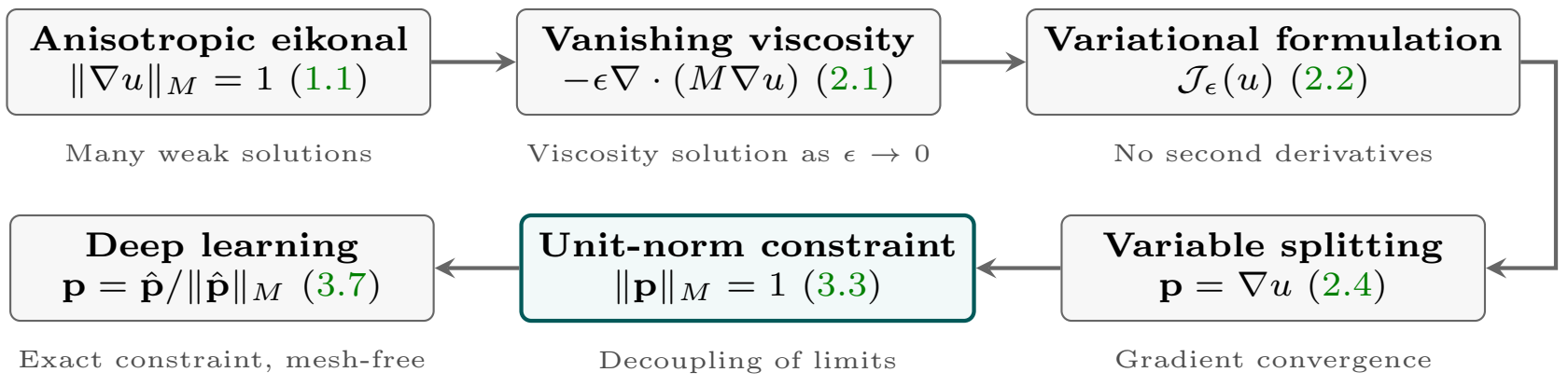

As outlined in Figure~\ref{fig:method_overview}, we minimize an energy functional whose Euler--Lagrange equation is the viscous anisotropic eikonal equation~\cite{varadhan1967behavior,bonnans2022linear}. Unlike the PDE residual, the energy involves no second derivatives, but features an exponentially weighted gradient term, $\mathrm e^{-2u/\epsilon}\|\nabla u\|_M^2$. To act directly on the gradient, we introduce an auxiliary variable $\mathbf p$ and enforce the constraint $\mathbf p=\nabla u$ through a penalty with parameter~$\lambda$. The resulting functional couples $u$ and $\mathbf p$ through the exponentially weighted term $\mathrm e^{-2u/\epsilon}\|\mathbf p\|_M^2$. Recovering the viscosity solution requires the joint limits $\epsilon\to0$ and $\lambda\to\infty$, whose interaction poses numerical challenges: large $\lambda$ leads to ill-conditioning, while small $\epsilon$ concentrates the exponential weight near~$\Gamma$. In neural network training, this parameter interaction also causes training oscillations. ALM and related multiplier approaches~\cite{belyaev2020admm,ennaji2021augmented,kim2027115466neural} avoid large penalty parameters but introduce an additional multiplier function that must be updated at sampled points through an outer iteration. Instead of strengthening the enforcement of $\mathbf p=\nabla u$, we reformulate the constraint to reduce the burden on the penalty.

The eikonal equation~\eqref{eq:eikonal} prescribes the magnitude of the gradient and leaves only its direction unknown. We therefore impose~$\|\mathbf p\|_M=1$ exactly, for every~$\epsilon$ and~$\lambda$. The exponential weight then no longer multiplies the gradient: the viscosity term depends on $u$ alone and favors larger $u$, while the penalty aligns $\nabla u$ with the unit field $\mathbf p$ and so penalizes slopes larger than one. Under the compatibility condition, the largest function with $\|\nabla u\|_M\le1$ and $u=g$ on $\Gamma$ is the viscosity solution~\cite{belyaev2020admm,ennaji2021augmented}; the neural method of~\cite{kim2027115466neural} enforces this characterization with Lagrange multipliers for $M=I$, $g=0$, whereas here the exponential weight and the unit-norm constraint on $\mathbf p$ realize it without multipliers. Hence the viscosity term selects the correct solution, while the eikonal constraint is carried by $\mathbf p$ instead of emerging only as $\epsilon\to0$.

As a result, the two limits decouple: the unit-norm constraint enforces the eikonal condition on $\mathbf p$, while $\lambda\to\infty$ enforces $\nabla u=\mathbf p$ at any fixed $\epsilon>0$. Proposition~\ref{prop:penalty_convergence} proves this for a general metric and $\Gamma=\partial\Omega$: global minimizers satisfy $u_{\epsilon,\lambda}\to u^*$ and $\mathbf p_{\epsilon,\lambda},\nabla u_{\epsilon,\lambda}\to\nabla u^*$ strongly in $L^2$ as $\lambda\to\infty$. This route still requires a large $\lambda$. With $\lambda$ fixed instead, the exponential weight confines the violation of $\nabla u=\mathbf p$ to a layer of width $O(\epsilon)$ near $\Gamma$, where the penalty limits it to an excess slope of $O(1/\lambda)$. For $M=I$ and $g=0$ on $\Gamma=\partial\Omega$ with $\Delta u^*\le K$ (e.g., convex or smooth domains), Proposition~\ref{prop:bound} accordingly bounds $\|u_{\epsilon,\lambda}-u^*\|_{L^\infty}$ and $\|\nabla u_{\epsilon,\lambda}-\nabla u^*\|_{L^2}^2$ by $O(\epsilon/\lambda)$. A moderate $\lambda$ thus suffices as $\epsilon\to0$.

Deep learning is well suited to this formulation for two reasons. First, the unit-norm constraint can be enforced exactly through the network architecture. Second, its mesh-free nature eliminates the need for mesh refinement as the functional becomes increasingly concentrated near~$\Gamma$ with decreasing $\epsilon$. Instead, we resample random points at every training epoch, ensuring that points near~$\Gamma$ retain a nonzero probability of being sampled, even as this probability decreases with~$\epsilon$. We keep $\lambda$ fixed and decrease $\epsilon$ whenever the training loss stagnates. Numerical experiments demonstrate that this strategy recovers the viscosity solution across isotropic, discontinuous, and anisotropic metrics, with both point and extended sources, extending beyond the setting covered by our convergence analysis. Ablation studies comparing against a HotSpot-type loss (Section~\ref{ex:random_pnts}) and alternative $\lambda$ strategies, together with surface reconstruction from the unorganized point clouds, further support the effectiveness of the proposed approach.


The remainder of the paper is organized as follows. Section~\ref{sec:prelim} reviews the variational formulation of the vanishing-viscosity approximation and its penalized variable splitting. Section~\ref{sec:method} introduces the unit-norm gradient variable, proves the convergence results, and describes the neural network implementation. Section~\ref{sec:experiments} presents numerical results and ablation studies, and Section~\ref{sec:conclusion} concludes. The proofs are given in Appendices~\ref{appen:proof_fix_epsilon} and~\ref{appen:proof_fix_lambda}.

\section{Preliminary}\label{sec:prelim}
\subsection{Variational Formulation for Vanishing Viscosity}\label{sec:VarAppOpSp}
The \textit{method of vanishing viscosity} approximates~\eqref{eq:eikonal} by adding a diffusive term with a viscosity parameter~$\epsilon>0$:
\begin{equation}\label{eq:viscosity_eikonal}
	\begin{cases}
		-\epsilon\nabla\cdot\left(M\nabla u_\epsilon\right) +\left\Vert\nabla u_\epsilon\right\Vert_M^2 = 1 & \ \text{in }\Omega \setminus \Gamma,\\
		u_\epsilon=g & \ \text{on }\Gamma,
	\end{cases}
\end{equation}
and $u_\epsilon$ converges uniformly to the viscosity solution $u^*$ of~\eqref{eq:eikonal} as $\epsilon\to0$~\cite{varadhan1967behavior,crandall1983viscosity}.
Multiplying~\eqref{eq:viscosity_eikonal} by the positive function $\frac{1}{\epsilon}\mathrm{e}^{-2u/\epsilon}$ gives
\begin{equation*}
- \frac{1}{\epsilon}\mathrm{e}^{ -\frac{2u}{\epsilon} } \left(1 + \left\Vert\nabla u \right\Vert_M^2\right) -  \nabla \cdot \left( \mathrm{e}^{ -\frac{2u}{\epsilon} } M \nabla u\right) = 0,
\end{equation*}
which is the Euler--Lagrange equation of the functional
\begin{equation}\label{eq:direct_ftnal}
	\mathcal{J}_\epsilon(u)  = \frac{1}{2}\int_{\Omega} s^{u}_\epsilon \left(1+\left\Vert\nabla u\right\Vert_M^2 \right)\diff\bx, \qquad s^{u}_\epsilon(\bx) = \mathrm{e}^{ -\frac{2 u(\bx) }{\epsilon} }.
\end{equation}
Under the Hopf--Cole transformation $w=\mathrm{e}^{-u/\epsilon}$, $\mathcal{J}_\epsilon(u)=\frac12\int_\Omega\bigl(w^2+\epsilon^2\|\nabla w\|_M^2\bigr)\diff\bx$ is the energy of the linear equation $-\epsilon^2\nabla\cdot(M\nabla w)+w=0$~\cite{varadhan1967behavior,bonnans2022linear} and is strictly convex in $w$. The solution of~\eqref{eq:viscosity_eikonal} is therefore characterized as the minimizer
\begin{equation}\label{eq:var_approach}
	u_\epsilon = \argmin_{u\in H^1\left(\Omega, \mathbb{R}\right)} \Bigl\{ \mathcal{J}_\epsilon(u) : u=g\text{ on }\Gamma \Bigr\},
\end{equation}
which extends the formulation used in HotSpot~\cite{wang2024hotspot} from $M=I$ and $g=0$ to general $M$ and $g$. The functional contains only first derivatives of $u$, but its weight $s_\epsilon^u$ concentrates exponentially near $\Gamma$ as $\epsilon\to0$. Consequently, as $\epsilon\to0$, direct minimization of $\mathcal{J}_\epsilon$ becomes increasingly ineffective away from the sources $\Gamma$, since the energy is exponentially localized near $\Gamma$ (see Section~\ref{ex:random_pnts}).

The trace condition $u=g$ on $\Gamma$ requires $\Gamma$ to have positive capacity, e.g., $\Gamma=\partial\Omega$ or, for $d=2$, a union of curves. Since isolated points have zero capacity for $d\geq2$, the source condition cannot be enforced through an $H^1$ trace constraint. 
Instead, we impose it directly through the neural representation~\eqref{eq:eik_primary_variable} (Section~\ref{sec:NN_method}), which restricts the admissible functions to a regularized class and limits arbitrarily sharp variations near the source that can cause degeneracy in the $H^1$ variational formulation.

\subsection{Penalized Variable-Splitting Formulation}
To facilitate the optimization of the variational formulation, we introduce an \textit{auxiliary variable} $\mathbf{p}:\Omega\rightarrow\mathbb{R}^d$ for $\nabla u$ and rewrite~\eqref{eq:var_approach} as
\begin{align}
	&\min_{u \in H^1(\Omega, \mathbb{R}),\ \mathbf{p} \in L^2(\Omega, \mathbb{R}^d)}
	 \Bigl\{ \mathcal{V}_\epsilon(u, \mathbf{p}) : u=g\text{ on }\Gamma, \: \nabla u = \mathbf{p}\text{ in }\Omega \Bigr\}, \label{eq:constraint_prob}\\
	&\mathcal{V}_\epsilon(u, \mathbf{p})  = \frac{1}{2}\int_{\Omega} s^{u}_\epsilon \left(1+\left\Vert \mathbf{p} \right\Vert_M^2 \right)\diff\bx. \label{eq:loss_vis}
\end{align}
By itself, the splitting does not remove the nonlinearity: $u$ and $\mathbf p$ remain coupled through $s_\epsilon^u\|\mathbf p\|_M^2$. A standard penalty method with parameters $\lambda>0$ and $\lambda_\Gamma>0$ gives
\begin{equation}\label{eq:unconstrained_prob}
(u_{\epsilon,\lambda,\lambda_\Gamma},\mathbf{p}_{\epsilon,\lambda,\lambda_\Gamma}) = \argmin_{u \in H^1(\Omega, \mathbb{R}),\ \mathbf{p} \in L^2(\Omega, \mathbb{R}^d)}
 \cH_\epsilon(u, \mathbf{p}; \lambda, \lambda_{\Gamma}),
\end{equation}
where
\begin{subequations}\label{eq:pen_obj_ftnal}
	\begin{align}
	\cH_\epsilon\left(u, \mathbf{p}; \lambda, \lambda_{\Gamma}\right) &= \mathcal{V}_\epsilon(u, \mathbf{p}) +  \frac{\lambda}{2}\mathcal{L}_{\text{GM}}(u, \mathbf{p}) + \frac{\lambda_\Gamma}{2}\mathcal{L}_{\Gamma}(u), \\ \mathcal{L}_{\text{GM}}(u, \mathbf{p}) &= \int_\Omega \left\Vert \nabla u - \mathbf{p} \right\Vert^2 \diff\bx, \label{eq:loss_gm} \\
	\mathcal{L}_{\Gamma}(u) &= \int_\Gamma\left(u-g\right)^2\diff S(\bx),
    \label{eq:loss_boundary}
	\end{align}
\end{subequations}
with the Euclidean norm $\|\cdot\|$ in~\eqref{eq:loss_gm}. We call $\mathcal{L}_{\text{GM}}$ and $\mathcal{L}_{\Gamma}$ the gradient-matching (as in~\cite{park2024ReSDF}) and boundary-matching functionals. The penalty formulation recovers~\eqref{eq:constraint_prob} only as $\lambda\to\infty$ and the viscosity solution to~\eqref{eq:eikonal} only as $\epsilon\to0$, and the two limits interact unfavorably. In particular, increasing $\lambda\rightarrow\infty$ may cause numerical instability.


\section{Method}\label{sec:method}

\subsection{Variational Formulation with a Unit-Norm Gradient Variable}\label{sec:unitnorm_variation}
To prevent the gradient-matching term from inheriting the exponential localization of $s_\epsilon^u$, we impose the unit norm prescribed by~\eqref{eq:eikonal} on the auxiliary variable exactly,
\begin{equation}\label{eq:p_unit_constraint}
	\|\mathbf{p}\|_M = 1  \qquad\text{in }\Omega,
\end{equation}
for every $\epsilon>0$ and $\lambda>0$. Then $1+\|\mathbf p\|_M^2=2$, and $\mathcal{V}_\epsilon$ in~\eqref{eq:loss_vis} reduces exactly to
\begin{equation}\label{eq:reduced_vis_functional}
	\mathcal{F}_\epsilon(u)\coloneqq\int_\Omega s_\epsilon^u\,\diff\bx.
\end{equation}
The \textbf{variational formulation with a unit-norm gradient variable} reads
\begin{align}\label{eq:proposed_formulation}
	\min_{\substack{u\in H^1(\Omega,\mathbb{R}),\; \mathbf{p}\in L^2(\Omega,\mathbb{R}^d)}} \quad &
	\mathcal E_{\epsilon,\lambda}(u,\mathbf p)\coloneqq \mathcal{F}_\epsilon(u) + \frac{\lambda}{2}\mathcal{L}_{\mathrm{GM}}(u,\mathbf{p})\\
    \nonumber
	\text{s.t.} \qquad\qquad\ \  & u=g\quad\text{on }\Gamma, 	\qquad 	\left\Vert\mathbf{p}\right\Vert_M=1 \quad\text{in }\Omega.
\end{align}
The constraint removes the gradient from the exponentially weighted term: $\mathcal F_\epsilon$ depends on $u$ alone and decreases as $u$ increases, while $\mathcal L_{\mathrm{GM}}$ aligns $\nabla u$ with the unit field $\mathbf p$.
\begin{remark}
The heat method~\cite{crane2017The} also fits a potential to a unit field by least squares, i.e., minimizes $\mathcal L_{\mathrm{GM}}(\cdot,\mathbf p)$, but takes $\mathbf p=-\nabla w/\|\nabla w\|$ fixed from a heat step $(\mathrm{id}-t\Delta)w=\delta_\Gamma$, related to $(\mathrm{id}-t\Delta)w=0$, $w=1$ on $\Gamma$, the Hopf--Cole equation of Section~\ref{sec:VarAppOpSp} with $M=I$, $g=0$, $t=\epsilon^2$; here $\mathbf p$ is minimized jointly with $u$, and $\mathcal F_\epsilon$ selects the viscosity solution within the same functional. 
\end{remark}

\subsection{Analysis on Decoupling of Asymptotic Limits}\label{sec:theory}
With~\eqref{eq:p_unit_constraint} in place, the viscosity solution is reached along two separate routes, both for $\Gamma=\partial\Omega$: (i) $\lambda\to\infty$ with $\epsilon>0$ fixed, for a general metric $M$, and (ii) $\epsilon\to0$ with $\lambda>0$ fixed, for $M=I$ and $g=0$.

\begin{proposition}[Convergence for fixed viscosity]
\label{prop:penalty_convergence}
Let $\Omega\subset\mathbb{R}^d$ be a bounded Lipschitz domain with
$\Gamma=\partial\Omega$. Assume that
$M\in C(\overline{\Omega};\mathbb{R}^{d\times d})$ is symmetric and uniformly positive definite and that~\eqref{eq:eikonal} admits a unique Lipschitz viscosity solution $u^*\in W^{1,\infty}(\Omega)$ with $u^*=g$ on $\Gamma$.

Fix $\epsilon>0$ and assume a global minimizer of $\mathcal{E}_{\epsilon,\lambda}$ in~\eqref{eq:proposed_formulation} exists and let
\[
(u_{\epsilon,\lambda},\boldsymbol{p}_{\epsilon,\lambda}) = \operatorname*{arg\,min}_{(u,\boldsymbol{p})\in\mathcal{A}} \mathcal{E}_{\epsilon,\lambda}(u,\boldsymbol{p}),\quad \forall\lambda>0,
\]
where $
\mathcal{A} := \left\{ (u,\boldsymbol{p})\in W^{1,2}(\Omega) \times L^2(\Omega;\mathbb{R}^d) \;\middle|\; 
u=g\ \text{on }\Gamma,\ 
\left\|\boldsymbol{p}\right\|_M=1 \ \text{a.e. in }\Omega
\right\}$.

Then, as $\lambda\to\infty$,
\[
\begin{aligned}
u_{\epsilon,\lambda}&\to u^*
&&\text{strongly in }L^2(\Omega),\\
\boldsymbol{p}_{\epsilon,\lambda}&\to\nabla u^*
&&\text{strongly in }L^2(\Omega;\mathbb{R}^d),\\
\nabla u_{\epsilon,\lambda}&\to\nabla u^*
&&\text{strongly in }L^2(\Omega;\mathbb{R}^d).
\end{aligned}
\]
\end{proposition}
The proof is given in Appendix~\ref{appen:proof_fix_epsilon}. As $\lambda\to\infty$, the penalty enforces $\mathbf p=\nabla u$, so the limit is a subsolution, $\|\nabla u\|_M\le1$ with $u=g$ on $\Gamma$; since $\mathcal F_\epsilon$ favors larger $u$, the maximal-subsolution property of the viscosity solution identifies the limit with $u^*$. Numerically, however, a large $\lambda$ makes the optimization ill-conditioned.

The error bounds of route~(ii) assume $\Delta u^*\le K$, which holds with $K=0$ on convex domains and with a finite $K$ on domains satisfying a uniform exterior sphere condition, such as $C^2$ domains (Appendix~\ref{appen:proof_fix_lambda}).

\begin{proposition}[Convergence for fixed penalty]
\label{prop:bound}
Let $\Omega\subset\mathbb{R}^d$ be a bounded Lipschitz domain with $\Gamma=\partial\Omega$, let $M=I$ and $g=0$, and let $u^*\in W^{1,\infty}(\Omega)\cap C(\overline{\Omega})$ be the viscosity solution of~\eqref{eq:eikonal}.

\textup{(i)} For all $\epsilon,\lambda>0$, $\mathcal E_{\epsilon,\lambda}$ has a unique global minimizer $(u_{\epsilon,\lambda},\mathbf p_{\epsilon,\lambda})$ on $\mathcal A$; it satisfies $u_{\epsilon,\lambda}\ge u^*$ a.e.\ in $\Omega$ and $\mathcal L_{\mathrm{GM}}(u_{\epsilon,\lambda},\mathbf p_{\epsilon,\lambda})\le\frac{2C_\Omega}{\mathrm e}\frac{\epsilon}{\lambda}$, where $C_\Omega=\sup_{t>0}t^{-1}\bigl|\{\bx\in\Omega:u^*(\bx)<t\}\bigr|<\infty$ depends only on $\Omega$.

\textup{(ii)} If, in addition, $\Delta u^*\le K$ in $\mathcal D'(\Omega)$ for some $K\ge0$, then, for $0<\epsilon\le1/K$ (any $\epsilon>0$ if $K=0$),
\begin{equation}\label{eq:err_bound_fix_lambda}
	\left\|u_{\epsilon,\lambda}-u^*\right\|_{L^\infty(\Omega)}\le\frac{\epsilon}{2\lambda}(1+K\epsilon),
\end{equation}
\begin{equation}\label{eq:grad_bound_fix_lambda}
	\left\|\mathbf p_{\epsilon,\lambda}-\nabla u^*\right\|_{L^2(\Omega)}^2\le\left\|\nabla u_{\epsilon,\lambda}-\nabla u^*\right\|_{L^2(\Omega)}^2\le\Bigl(\frac{4C_\Omega}{\mathrm e}+2K|\Omega|\Bigr)\frac{\epsilon}{\lambda}.
\end{equation}
\end{proposition}
The proof is provided in Appendix~\ref{appen:proof_fix_lambda}. Since $s_\epsilon^u$ is exponentially small away from $\Gamma$, the violation of $\nabla u=\mathbf p$, an excess slope of order at most $1/\lambda$, is confined to a layer of width $O(\epsilon)$ near $\Gamma$. The resulting error $O(\epsilon/\lambda)$ vanishes as $\epsilon\to0$ for every fixed $\lambda>0$. This result motivates the method of the next section, which keeps a moderate $\lambda$ fixed while $\epsilon$ decreases and thereby avoids the ill-conditioning associated with increasing $\lambda$.

\subsection{Deep Learning Framework with a Unit-Norm Gradient Variable} \label{sec:NN_method}

\begin{figure}[t]
	\centering
		\begin{tikzpicture}[font=\small]
			\node (input) [rounded corners=4pt, text centered, align=center, inner sep=0pt, draw, thick, minimum size=0.6cm, fill=green!30] {$\mathbf{x}$};
			\node[shift={(0.0,-0.3)}] at (input.south) {\footnotesize \texttt{input}};
			\node (layer1) [right=0.6cm of input, draw, thick, minimum width=0.35cm, minimum height=1.5cm, fill=cyan!30] {};
			\node (layer2) [right=0.4cm of layer1, draw, thick, minimum width=0.35cm, minimum height=1.5cm, fill=cyan!30] {};
			\node (dots) [right=0.4cm of layer2] {...};
			\node (layerN) [right=0.4cm of dots, draw, thick, minimum width=0.35cm, minimum height=1.5cm, fill=cyan!30] {};
			\coordinate (common_output) at ($(layerN.east) + (0.7,0)$);
			\coordinate (common_output_N) at ($(common_output) + (0,0.6)$);
			\coordinate (common_output_S) at ($(common_output) + (0,-0.6)$);
			\begin{pgfonlayer}{background}
				\node (back) [rounded corners=6pt, fit=(layer1) (layer2) (layerN), fill=gray!20, inner sep=0.3cm] {};
				\node[shift={(0.0,-0.25)}] at (back.south) {\footnotesize \texttt{Network}};
			\end{pgfonlayer}
			\node (hat_p) [rounded corners=4pt, text centered, align=center, inner sep=0pt, draw, thick, minimum width=1.2cm, minimum height=0.6cm, fill=gray!30] at ($(common_output_S) + (1.3,0)$) {$\hat{\mathbf{p}}_{\theta'}(\mathbf{x})$};
			\node (p) [rounded corners=4pt, text centered, align=center, inner sep=0pt, draw, thick, minimum width=1.2cm, minimum height=0.6cm, fill=pink!30] at ($(hat_p) + (3.0,0)$) {$\pnet(\mathbf{x})$};
			\node[anchor=west] at (p.east) {\footnotesize \texttt{Auxiliary Variable}};
			\node (u) [rounded corners=4pt, text centered, align=center, inner sep=0pt, draw, thick, minimum width=1.2cm, minimum height=0.6cm, fill=yellow!30] at ($(common_output_N) + (1.3,0)$) {$\hat{u}_\theta(\mathbf{x})$};
			\node (usol) [rounded corners=4pt, text centered, align=center, inner sep=0pt, draw, thick, minimum width=1.2cm, minimum height=0.6cm, fill=pink!30] at ($(u) + (3.0,0)$) {$u_\theta(\mathbf{x})$};
			\node[anchor=west] at (usol.east) {\footnotesize \texttt{Primary Variable}};
			\draw[-Latex, thick] (input) -- (layer1);
			\draw[-Latex, thick] (layer1) -- (layer2);
			\draw[-Latex, thick] (layer2) -- (dots);
			\draw[-Latex, thick] (dots) -- (layerN);
			\draw[-, thick] (layerN) -- (common_output);
			\draw[-, thick] (common_output) -- (common_output_N);
			\draw[-Latex, thick] (common_output_N) -- (u);
			\draw[-, thick] (common_output) -- (common_output_S);
			\draw[-Latex, thick] (common_output_S) -- (hat_p);
			\draw[-Latex, thick] ($(hat_p.east)$) -- ($(p.west)$) node[midway, below] {\footnotesize $\hat{\mathbf{p}}_{\theta'}/\| \hat{\mathbf{p}}_{\theta'} \|_{M}$};
			\draw[-Latex, thick] ($(u.east)$) -- ($(usol.west)$) node[midway, above] {\footnotesize $\tilde g+\phi\hat{u}_\theta$};
		\end{tikzpicture}
		\caption{Network architecture: one network with $d+1$ outputs and the hard constraints~\eqref{eq:eik_primary_variable} and~\eqref{eq:eik_auxiliary_variable}.}\label{fig:network}
\end{figure}
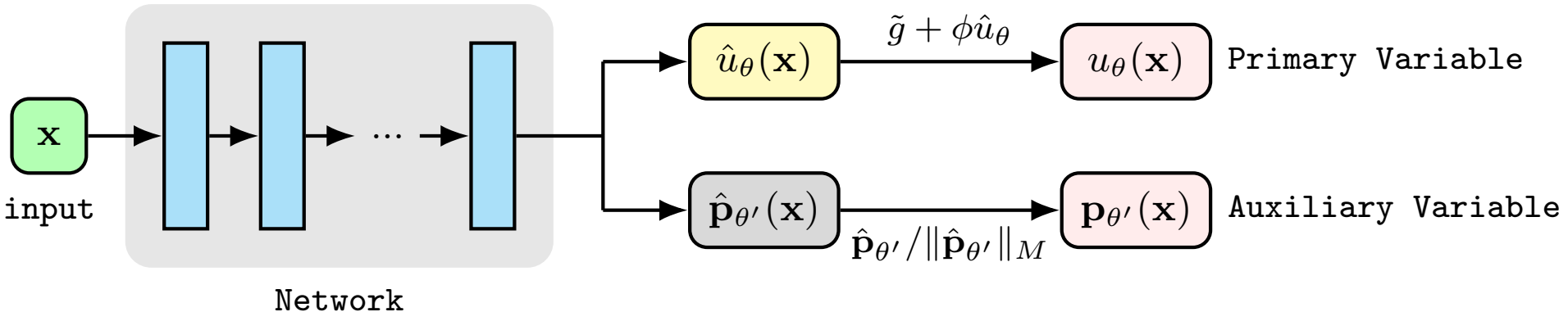

\paragraph{Hard-Constrained Network Representation}
As in ReSDF~\cite{park2024ReSDF}, a single fully connected network with $d+1$ outputs $(\hat u_\theta,\hat{\mathbf p}_{\theta'})$ represents both variables (Figure~\ref{fig:network}); $\theta$ and $\theta'$ distinguish the two output branches, which share all hidden layers. For a function $\tilde g\in\mathcal C^1(\overline\Omega)$ that extends $g$ and a level set function $\phi:\Omega\rightarrow\mathbb{R}$ whose zero level set represents $\Gamma$, we set
\begin{align}
	u_\theta(\bx) &= \tilde{g}(\bx)+\phi(\bx)\hat{u}_\theta(\bx), \label{eq:eik_primary_variable}\\
	\pnet(\bx) &= \frac{\hat{\mathbf{p}}_{\theta'}(\bx)}{\|\hat{\mathbf{p}}_{\theta'}(\bx)\|_M}, \label{eq:eik_auxiliary_variable}
\end{align}
where $\hat{u}_\theta$ is constrained to be nonnegative, e.g., by applying a softplus activation function at the output layer of the network. To avoid division by zero in~\eqref{eq:eik_auxiliary_variable}, we set the normalization parameter $\delta$ to $10^{-12}$, which is the default value used in PyTorch's \texttt{torch.nn.functional.normalize}. Then $u_\theta=g$ on $\Gamma$ exactly, as in the distance-function construction of~\cite{sukumar2022exact}, and $\|\pnet(\bx)\|_M=1$ at every $\bx\in\Omega$ up to $\delta$, independently of $\epsilon$ and $\lambda$. The training objective is
\begin{equation}\label{eq:pen_obj_ftnal_final}
	\mathcal{E}_{\epsilon,\lambda}(\theta,\theta') = \mathcal{F}_\epsilon(u_\theta) + \frac{\lambda}{2}\mathcal{L}_{\mathrm{GM}}(u_\theta,\pnet),
\end{equation}
which needs no penalty for the unit norm or the boundary condition. When no such $\phi$ is available (Examples~\ref{ex:shape_shading} and~\ref{ex:surface_recon}), we add the boundary-matching term $\frac{\lambda_\Gamma}{2}\mathcal{L}_\Gamma(u_\theta)$ of~\eqref{eq:loss_boundary} to~\eqref{eq:pen_obj_ftnal_final}, with $\mathcal L_\Gamma$ approximated by the mean of $(u_\theta-g)^2$ over a set $m_\Gamma\subset\Gamma$ of boundary collocation points. We set $\lambda_\Gamma=1$ in these experiments.

\paragraph{Mesh-Free Resolution of the Localized Functional}
We approximate the integrals in~\eqref{eq:pen_obj_ftnal_final} by Monte Carlo sampling,
\begin{equation}\label{eq:discrete_loss}
\widehat{\mathcal{E}}_{\epsilon,\lambda}(\theta,\theta')
=\frac{1}{|m_\Omega|}\sum_{\bx\in m_\Omega}\Bigl(s_\epsilon^{u_\theta}(\bx)
+\frac{\lambda}{2}\|\nabla u_\theta(\bx)-\pnet(\bx)\|^2\Bigr),
\end{equation}
where $m_\Omega=\{\bx_i\}_{i=1}^{N_c}$ consists of i.i.d.\ samples from the uniform distribution $U(\Omega)$, redrawn at every epoch.

This mesh-free sampling strategy is particularly useful for resolving the exponential concentration of $\mathcal{F}_\epsilon$ near~$\Gamma$ as $\epsilon\to0$. For $M=I$ and $g=0$, the viscosity solution is the distance function $d_\Gamma(\bx):=\operatorname{dist}(\bx,\Gamma)$, and $s_\epsilon^{d_\Gamma}=\mathrm{e}^{-2d_\Gamma/\epsilon}$. For a numerical threshold $\eta\in(0,1)$, the region where $s_\epsilon^{d_\Gamma}>\eta$ is
\begin{equation}\label{eq:safe_region}
\Gamma_\epsilon
=\left\{\bx\in\Omega:
d_\Gamma(\bx)<\frac{\epsilon\log(\eta^{-1})}{2}\right\},
\end{equation}
which forms a layer of thickness $O(\epsilon)$ around $\Gamma$. Mesh-based methods therefore need to refine the mesh to resolve this increasingly narrow layer as $\epsilon$ decreases. In contrast, our mesh-free approach does not require mesh refinement: although the probability of sampling points in $\Gamma_\epsilon$ decreases with $\epsilon$, it remains positive under uniform sampling, allowing the localized functional to be evaluated in this region through repeated resampling at each epoch.

\begin{algorithm}[t]
    \caption{Algorithm for solving anisotropic eikonal equation~\eqref{eq:eikonal}}
    \label{alg:algorithm}
    \begin{algorithmic}[1]
        \State Initialize $\theta,\theta'$ and $\epsilon=10^{-k_{\min}}$.
        \For{$n=0,\ldots,N-1$}
            \State Sample $m_\Omega\sim U(\Omega)$.
            \State Evaluate $\widehat{\mathcal E}_{\epsilon,\lambda}$ in~\eqref{eq:discrete_loss}.
            \State Update $\theta,\theta'$ by gradient descent.
            \If{loss stagnates for $N_\epsilon$ epochs}
                \If{$\epsilon=10^{-k_{\max}}$}
                    \State Stop.
                \Else
                    \State $\epsilon\leftarrow\epsilon/10$.
                \EndIf
            \EndIf
        \EndFor
    \end{algorithmic}
\end{algorithm}
\paragraph{Automatic Annealing of the Viscosity Parameter}
Motivated by Proposition~\ref{prop:bound}, we fix $\lambda$ and gradually reduce $\epsilon$ during training. Such a continuation in the viscosity is also used by a finite element method that halves the viscosity coefficient level by level~\cite{yang2022continuous} and by neural solvers for distance functions~\cite{wang2024hotspot,krishnan2026viscoreg}; here it is applied with $\lambda$ fixed and triggered by stagnation of the training loss. We start from $\epsilon=10^{-k_{\min}}$, chosen so that $\Gamma_\epsilon$ in~\eqref{eq:safe_region} is wide enough to resolve the localized functional, and divide $\epsilon$ by $10$ until $\epsilon=10^{-k_{\max}}$, where $k_{\max}\geq k_{\min}$. Rather than prescribing a fixed number of epochs for each $\epsilon$, we reduce $\epsilon$ whenever $\widehat{\mathcal{E}}_{\epsilon,\lambda}$ fails to decrease for $N_\epsilon$ consecutive epochs, so that each viscosity level is optimized sufficiently but not longer than necessary (Algorithm~\ref{alg:algorithm}). A theoretical analysis of this schedule is left for future work, while its limited sensitivity to the scheduling choices is examined in Section~\ref{sec:experiments}.

\section{Numerical Results}\label{sec:experiments}
We examine the numerical properties of the proposed method\footnote{All computations use a single NVIDIA GV100 (TITAN V) GPU.} (Section~\ref{ex:random_pnts}), apply it to further isotropic and anisotropic examples (Section~\ref{ex:exmples}), and study its components in ablations (Section~\ref{sec:ablation}).

Unless stated otherwise, the computational domain is $\Omega=[-1,1]^2$ and we use a fully connected network with five hidden layers of 64 nodes and the softplus activation $\sigma(x)=\beta^{-1}\log(1+e^{\beta x})$ with $\beta=100$, trained by Adam~\cite{kingma2014adam} with initial learning rate $10^{-3}$, decayed by a factor of $0.99$ whenever the loss decreases. Each epoch uses $5{,}000$ collocation points sampled uniformly from $\Omega$, and Algorithm~\ref{alg:algorithm} uses $\lambda=1$, $k_{\min}=0$, $k_{\max}=5$ and $N_\epsilon=500$. For source points $\{\by_i\}_{i\in\mathcal I}$, we take $\phi(\bx)=\prod_{i\in\mathcal I}\|\bx-\by_i\|$ in~\eqref{eq:eik_primary_variable}. On the $N_h$ nodes $\bx_i$ of a uniform grid with spacing $h=1/256$ we report the mean squared errors
\[
E_2(u_\theta)=\frac1{N_h}\sum_{i=1}^{N_h}|u_\theta(\bx_i)-u^*(\bx_i)|^2,\qquad
E_2(\nabla u_\theta)=\frac1{N_h}\sum_{i=1}^{N_h}\|\nabla u_\theta(\bx_i)-\nabla u^*(\bx_i)\|^2,
\]
with $\nabla u_\theta$ computed by automatic differentiation, and the maximum error $E_\infty(u_\theta)=\max_i|u_\theta(\bx_i)-u^*(\bx_i)|$; the figures label $E_2$ as the ``$L^2$ error''. Each training curve shows one run, whereas Table~\ref{tab:all_errors} reports means and standard deviations over five seeds. The curves are plotted on a logarithmic scale, and their background colors mark the successive values of $\epsilon$ chosen by the annealing.

\subsection{Numerical Properties of the Proposed Method}\label{ex:random_pnts}
\begin{figure}[t]
	\centering
	\begin{minipage}[t]{0.49\textwidth}
		\centering
		\subfigure[]{\includegraphics[width=0.48\linewidth]{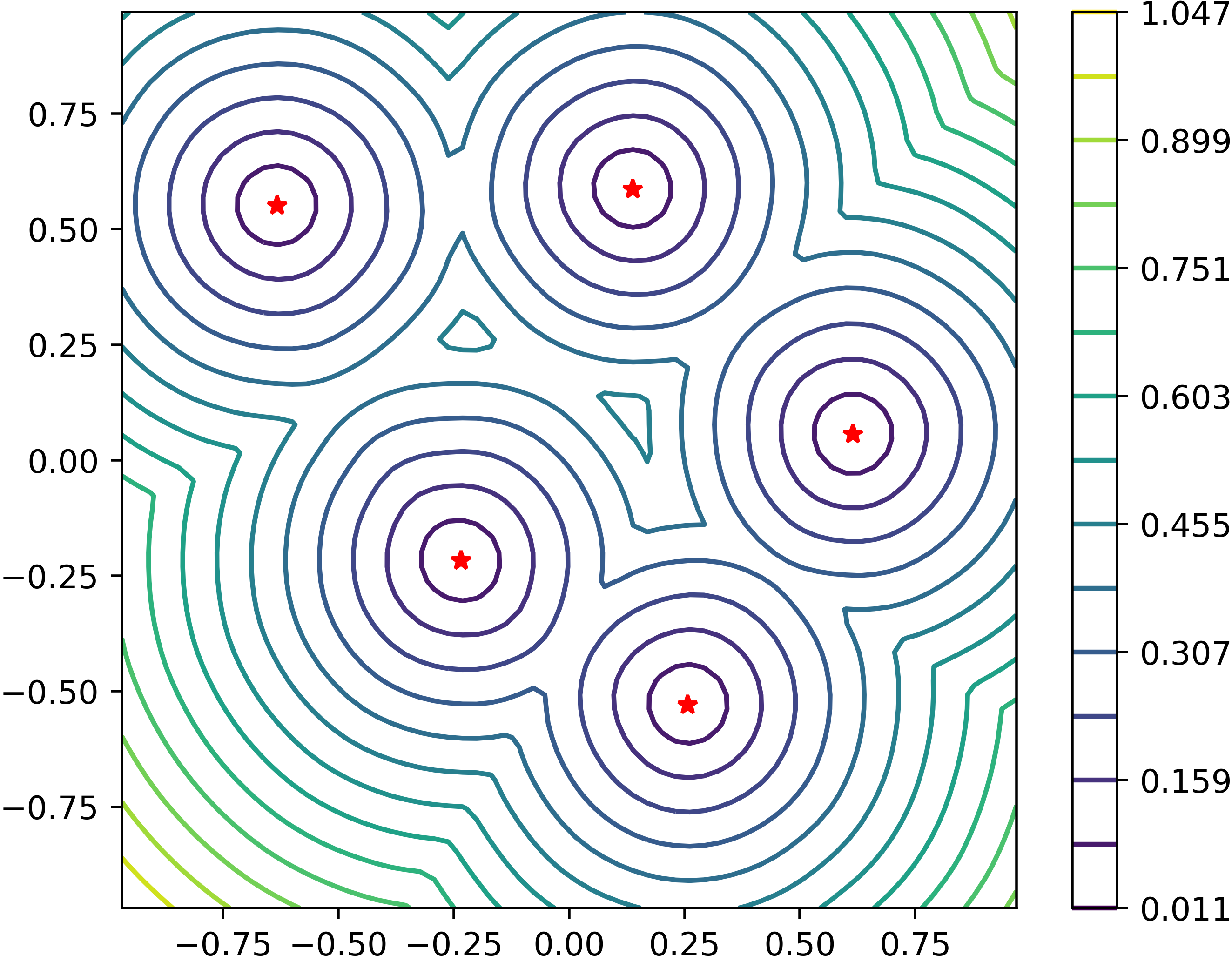}}\hfill
		\subfigure[]{\includegraphics[width=0.48\linewidth]{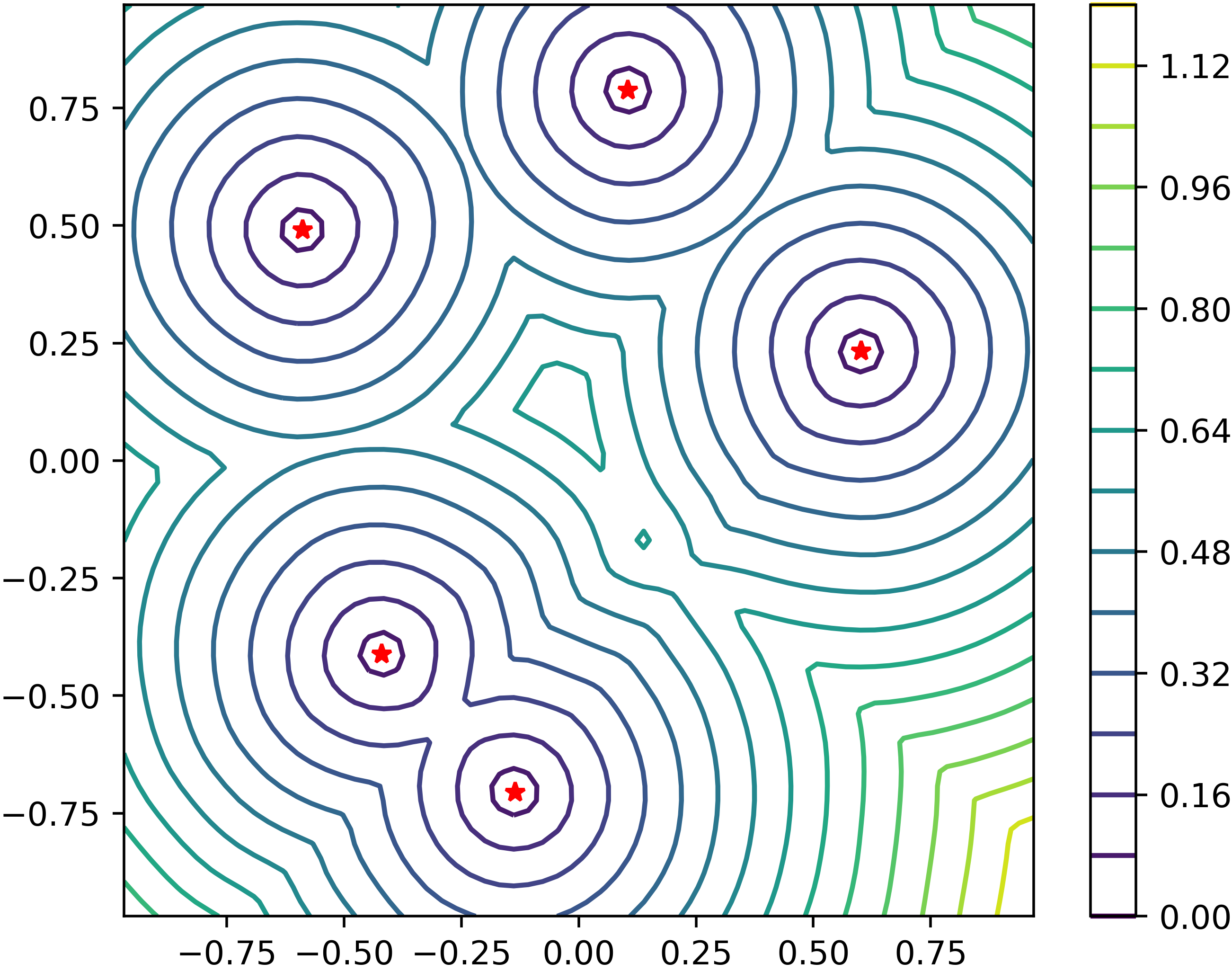}}
		\caption{Proposed method for the two sets of five random source points (red stars) in Section~\ref{ex:random_pnts}.}
		\label{fig:rand}
	\end{minipage}\hfill
	\begin{minipage}[t]{0.49\textwidth}
		\centering
		\subfigure[Ablation 1]{\includegraphics[width=0.48\linewidth]{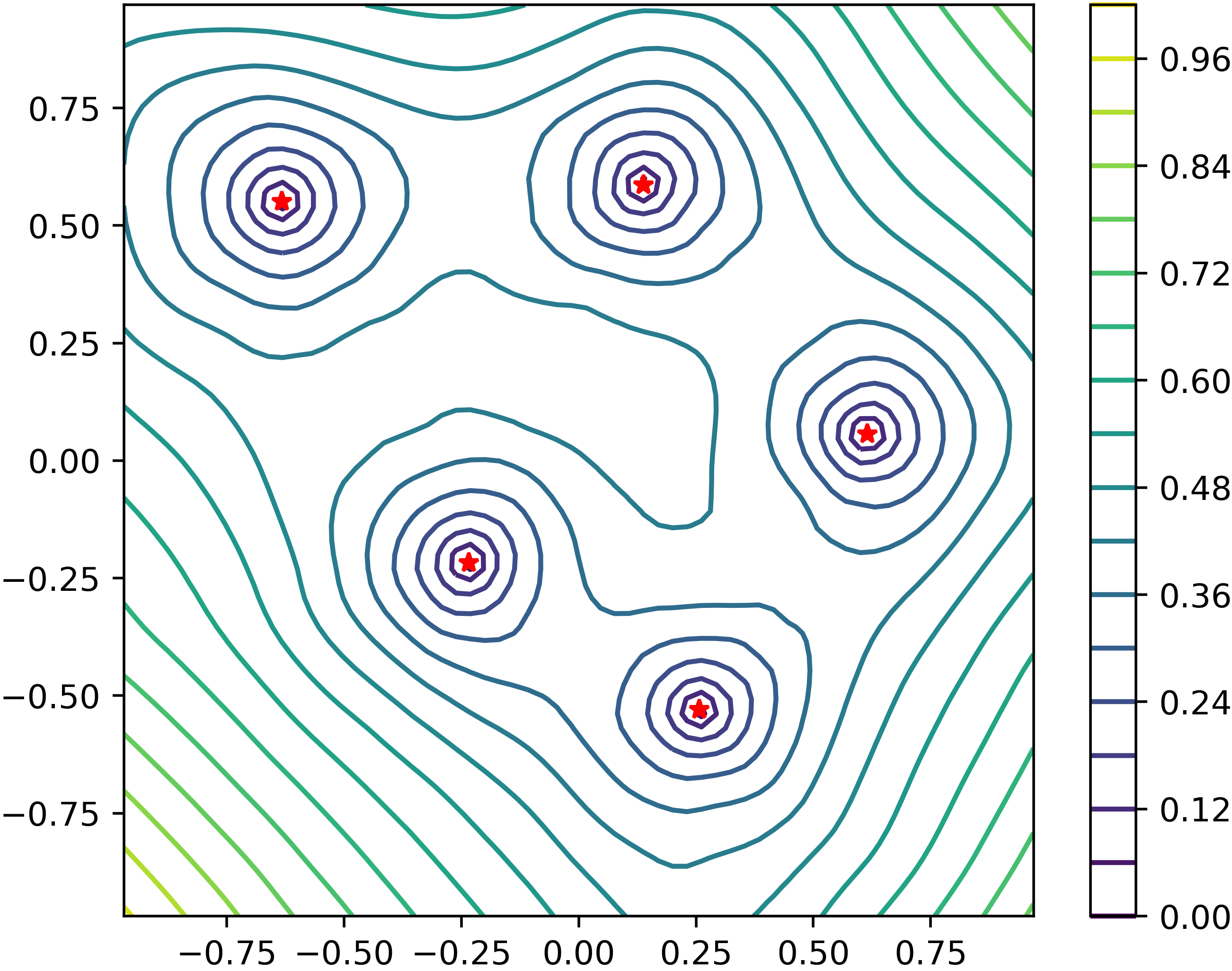}}\hfill
		\subfigure[Ablation 2]{\includegraphics[width=0.48\linewidth]{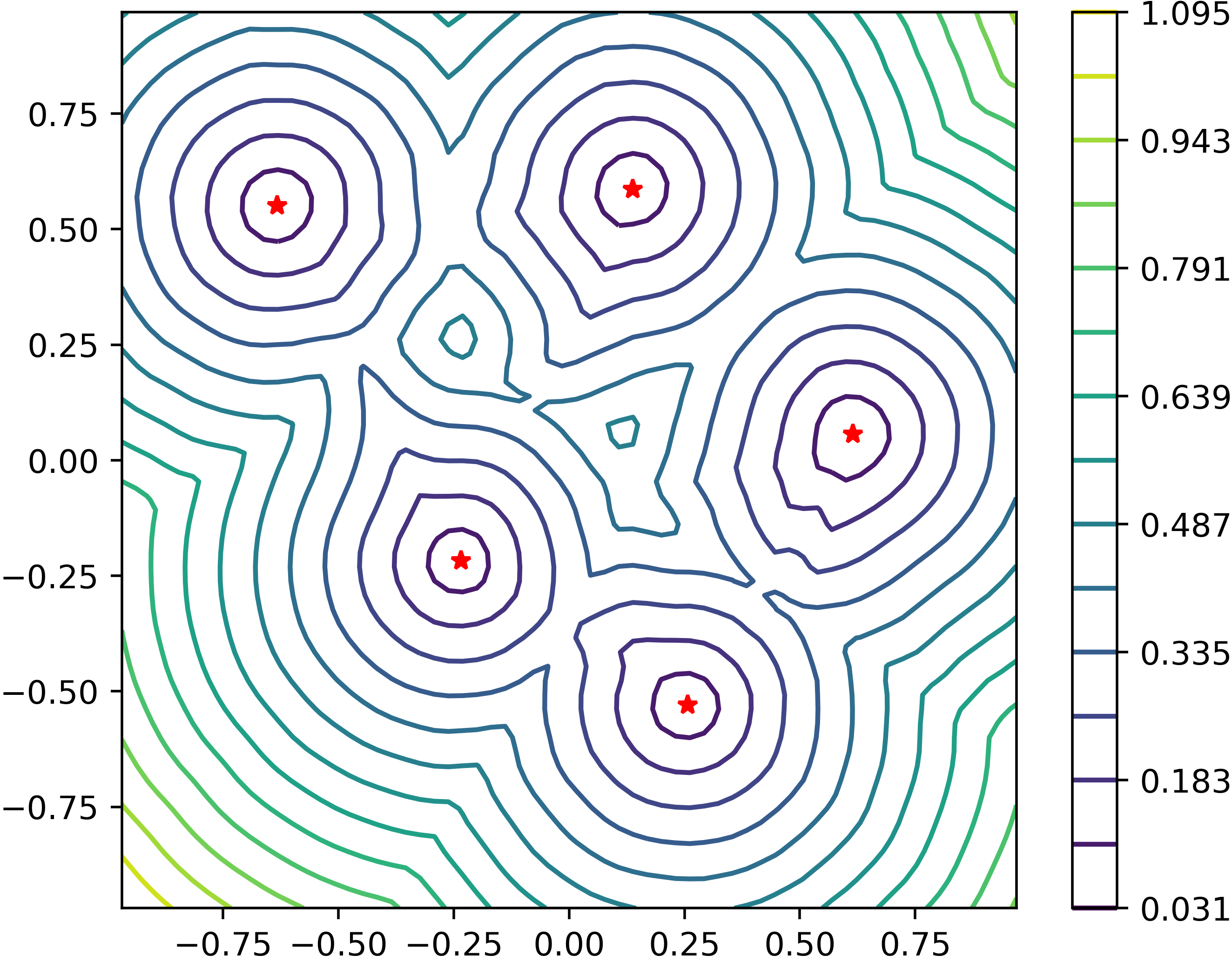}}
		\caption{Ablation studies for the sources of Figure~\ref{fig:rand}(a): minimizing~\eqref{eq:ablation1} (left) and~\eqref{eq:ablation2} (right).}
		\label{fig:rand_ablation}
	\end{minipage}
\end{figure}

\begin{figure}[t]
	\centering
	\subfigure[Automatic annealing with factor $1/2$]{\includegraphics[width=0.99\linewidth]{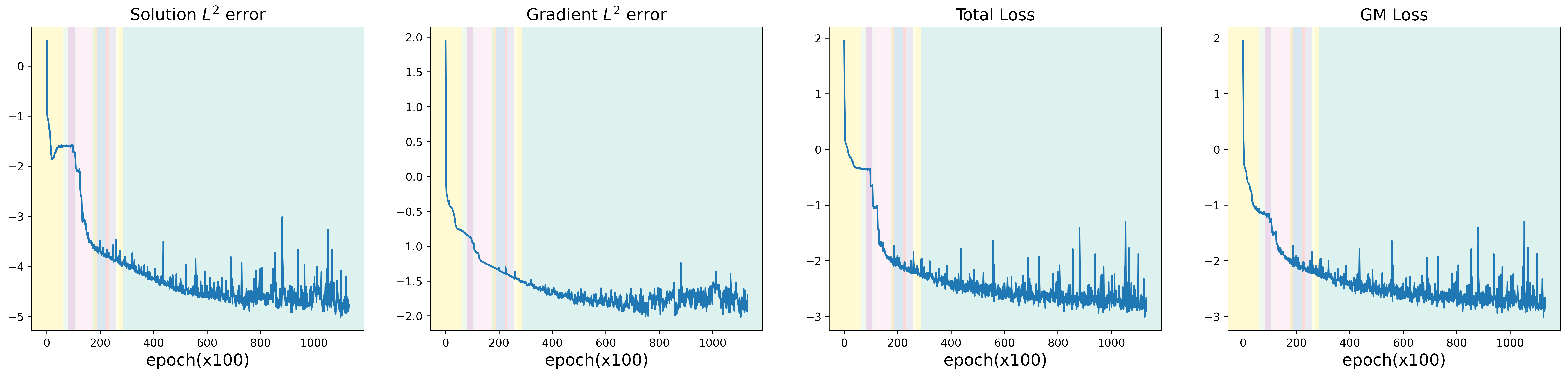}}\\
	\subfigure[Automatic annealing with factor $1/10$]{\includegraphics[width=0.99\linewidth]{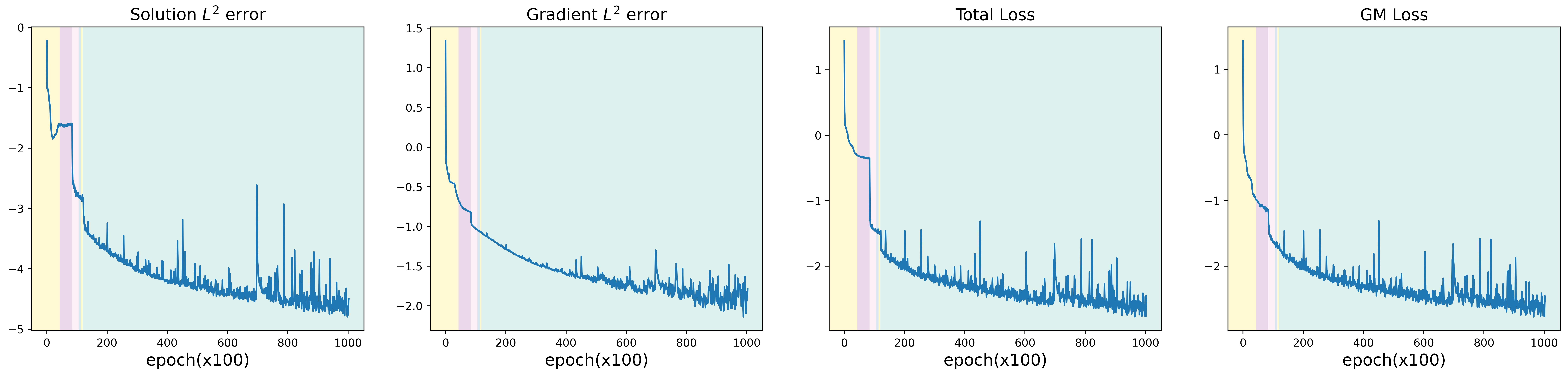}}
	\caption{Solution error, gradient error, total loss and GM loss during training for Figure~\ref{fig:rand}(a), for two annealing schedules of $\epsilon$ (logarithmic $y$-axes; the background colors mark the values of $\epsilon$).}
	\label{fig:rand_eps}
\end{figure}

We first take $M=I$, $g=0$ and two sets of five randomly chosen points as $\Gamma$, so that $u^*$ is the Euclidean distance to $\Gamma$. Figure~\ref{fig:rand} shows the solutions, and Figure~\ref{fig:rand_eps} the errors and losses during training for the first set (for the second set, see Figure~\ref{fig:rand_eps_2}). As $\epsilon$ decreases, the solution and gradient errors decrease, so the network approaches the viscosity solution; the gradient-matching (GM) loss $\cL_{\mathrm{GM}}$ in~\eqref{eq:proposed_formulation} decreases as well, so the constraint $\nabla u=\bp$ is well satisfied. We use the first set (Figure~\ref{fig:rand}(a)) to examine three further properties.

\begin{figure}[t]
	\centering
	\includegraphics[width=0.42\linewidth]{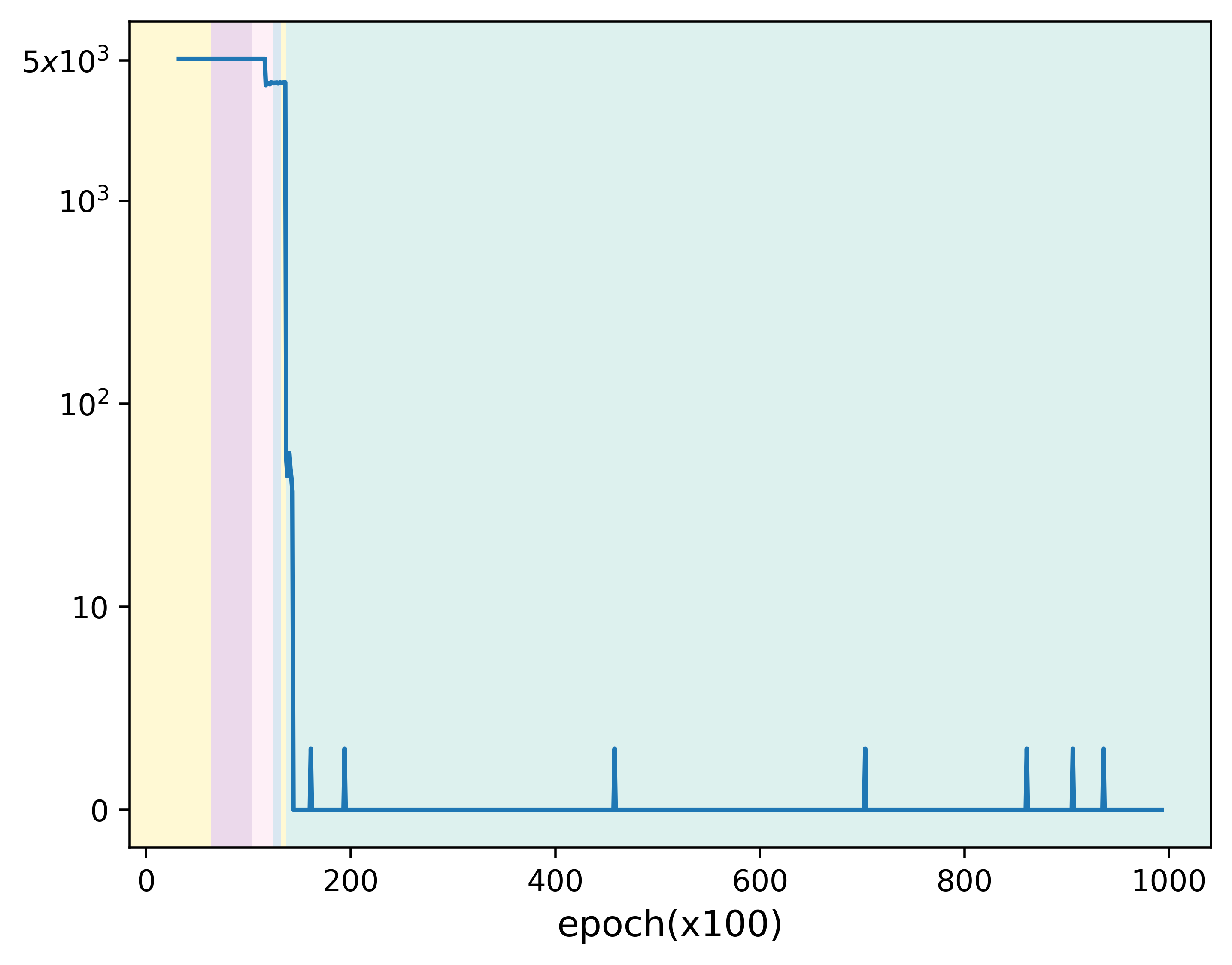}
	\caption{Number of collocation points (out of $5{,}000$ per epoch) in $\Gamma_\epsilon$, where $s_\epsilon^u$ is numerically nonzero, for Figure~\ref{fig:rand}(a) with the annealing factor $1/10$; the background colors mark the values of $\epsilon$.}
	\label{fig:validpnts}
\end{figure}

\paragraph{Sensitivity to the annealing schedule}
Besides the factor $1/10$ of Algorithm~\ref{alg:algorithm}, we decrease $\epsilon$ by a factor of $1/2$: $\epsilon_{n+1}=\epsilon_n/2$ with $\epsilon_0=1$ until $\epsilon_n\le10^{-5}$ (17 halvings, final $\epsilon=2^{-17}\approx7.6\times10^{-6}$), with the same stopping test as in Algorithm~\ref{alg:algorithm}. With both schedules the method learns the viscosity solution as $\epsilon$ decreases (Figure~\ref{fig:rand_eps}), so it is not sensitive to the choice between them. The same holds for Examples~\ref{ex:aniso_ex2} and~\ref{ex:aniso_ex3} (Figures~\ref{fig:lambda_aniso1}(d),(e) and~\ref{fig:lambda_aniso2}(d),(e)) and for the other examples with an exact solution in Section~\ref{ex:exmples} (factor $1/10$ in Figure~\ref{fig:curves_exact}, factor $1/2$ in Figure~\ref{fig:half_exact}). For the examples without an exact solution, the total and GM losses reach about the same final values with both factors (Figure~\ref{fig:loss_curves} and Figure~\ref{fig:half_loss}), although the factor $1/2$ takes more epochs.

\paragraph{Mesh-free resolution of the localized functional}
The integrand $s_\epsilon^u$ of $\mathcal F_\epsilon$ in~\eqref{eq:reduced_vis_functional} is numerically nonzero only in $\Gamma_\epsilon$ in~\eqref{eq:safe_region}. With the precision limit of PyTorch, $\eta=10^{-45}$, and the smallest viscosity $\epsilon=10^{-5}$, the part of $\Gamma_\epsilon$ around one source point has area about $2.68\pi\cdot10^{-7}$. We count the collocation points that fall in $\Gamma_\epsilon$ at each epoch (Figure~\ref{fig:validpnts}): all $5{,}000$ points do so at the first values of $\epsilon$, whereas at $\epsilon=10^{-5}$ the count is zero in most epochs and only occasionally a few points are sampled there, owing to the uniform random resampling. Although these contributions of $s_\epsilon^u$ are sparse in the final stage ($k=k_{\max}$), the solution and gradient errors and the GM loss keep decreasing in it (Figure~\ref{fig:rand_eps}).
In particular, the errors decrease more rapidly at larger values of $\epsilon$, suggesting that the overall solution profile is captured primarily during the early stages of optimization. At the final stage, the further reduction in errors appears to be driven primarily by the decrease in the GM loss, while $\mathcal F_\epsilon$ may still contribute through occasional collocation points sampled within $\Gamma_\epsilon$.

\paragraph{Effect of the unit-norm constraint: HotSpot-type ablations}
We compare the proposed loss with two losses without the unit-norm gradient variable:
\begin{itemize}
	\item \textit{Ablation 1:} the variational problem~\eqref{eq:var_approach} without splitting:
	\begin{equation}\label{eq:ablation1}
		\min_{u\in H^1(\Omega,\mathbb{R})}\ \mathcal{J}_\epsilon(u)+\frac{\lambda_\Gamma}{2}\cL_\Gamma(u).
	\end{equation}
	\item \textit{Ablation 2:} the eikonal equation added to~\eqref{eq:ablation1} as a soft constraint, without splitting:
	\begin{equation}\label{eq:ablation2}
		\min_{u\in H^1(\Omega,\mathbb{R})}\ \mathcal{J}_\epsilon(u)+\frac{\lambda}{2}\int_\Omega\left(\|\nabla u\|_M-1\right)^2\diff\bx+\frac{\lambda_\Gamma}{2}\cL_\Gamma(u).
	\end{equation}
	For $M=I$ and $g=0$, this is the loss of HotSpot~\cite{wang2024hotspot} for learning signed distance functions.
\end{itemize}
In both, $\cL_\Gamma(u)=0$ because, for a fair comparison, the boundary condition is imposed by~\eqref{eq:eik_primary_variable} as in the proposed method. With Ablation~1 (Figure~\ref{fig:rand_ablation}(a)) the network is trained correctly only near the sources: as $\epsilon$ decreases, the loss~\eqref{eq:ablation1} carries no information outside $\Gamma_\epsilon$, and the numerical solution is localized on $\Gamma$. In the proposed loss, the term $\frac{\lambda}{2}\mathcal L_{\mathrm{GM}}$ with $|\mathbf p|_M=1$ acts in all of $\Omega$ and prevents this localization (Figure~\ref{fig:rand}(a)). Adding the eikonal residual (Ablation~2) counteracts this localization and captures the overall shape (Figure~\ref{fig:rand_ablation}(b)), but does not approximate the viscosity solution. The difference from the proposed loss is explicit for $M=I$. For fixed $u$, the minimum of $\mathcal E_{\epsilon,\lambda}(u,\bp)$ in~\eqref{eq:proposed_formulation} over $\|\bp\|=1$ is attained at $\bp=\nabla u/\|\nabla u\|$, and
\begin{align*}
	\min_{\|\bp\|=1}\mathcal E_{\epsilon,\lambda}(u,\bp)&=\mathcal F_\epsilon(u)+\frac\lambda2\int_\Omega\left(\|\nabla u\|-1\right)^2\diff\bx,\\
	\mathcal J_\epsilon(u)&=\mathcal F_\epsilon(u)+\frac12\int_\Omega s_\epsilon^u\left(\|\nabla u\|^2-1\right)\diff\bx .
\end{align*}
Hence~\eqref{eq:ablation2} differs from the proposed loss with $\bp$ eliminated only by $\frac12\int_\Omega s_\epsilon^u(\|\nabla u\|^2-1)\diff\bx$, that is, by the gradient inside the exponential weight. We attribute the inaccuracy of Ablation~2 mainly to this product of the exponential weight and the gradient norm, which the unit-norm constraint eliminates. A second difference concerns training only: as $\bp$ is a separate network output, the network minimizes a loss that is quadratic in $\nabla u$ for given $\bp$ instead of the residual $(\|\nabla u\|-1)^2$, which is nonsmooth at $\nabla u=\mathbf 0$. Our ablations do not separate these two effects.

\subsection{Examples}\label{ex:exmples}
\begin{figure}[t]
	\centering
	\subfigure[Ex.~\ref{ex:disconti1}, proposed]{\includegraphics[width=0.3\textwidth]{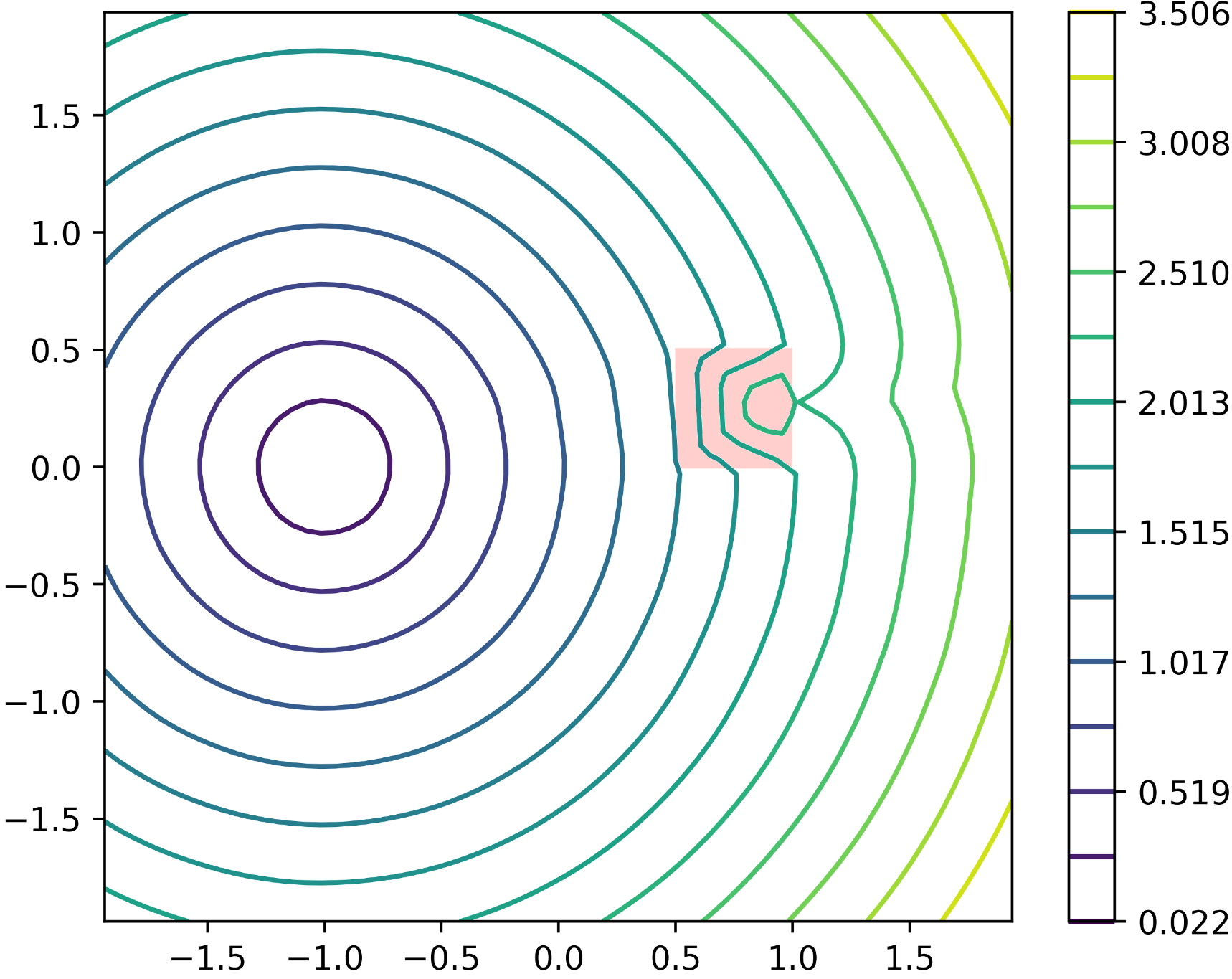}}\hfill
	\subfigure[Ex.~\ref{ex:disconti1}, Ablation 2]{\includegraphics[width=0.3\textwidth]{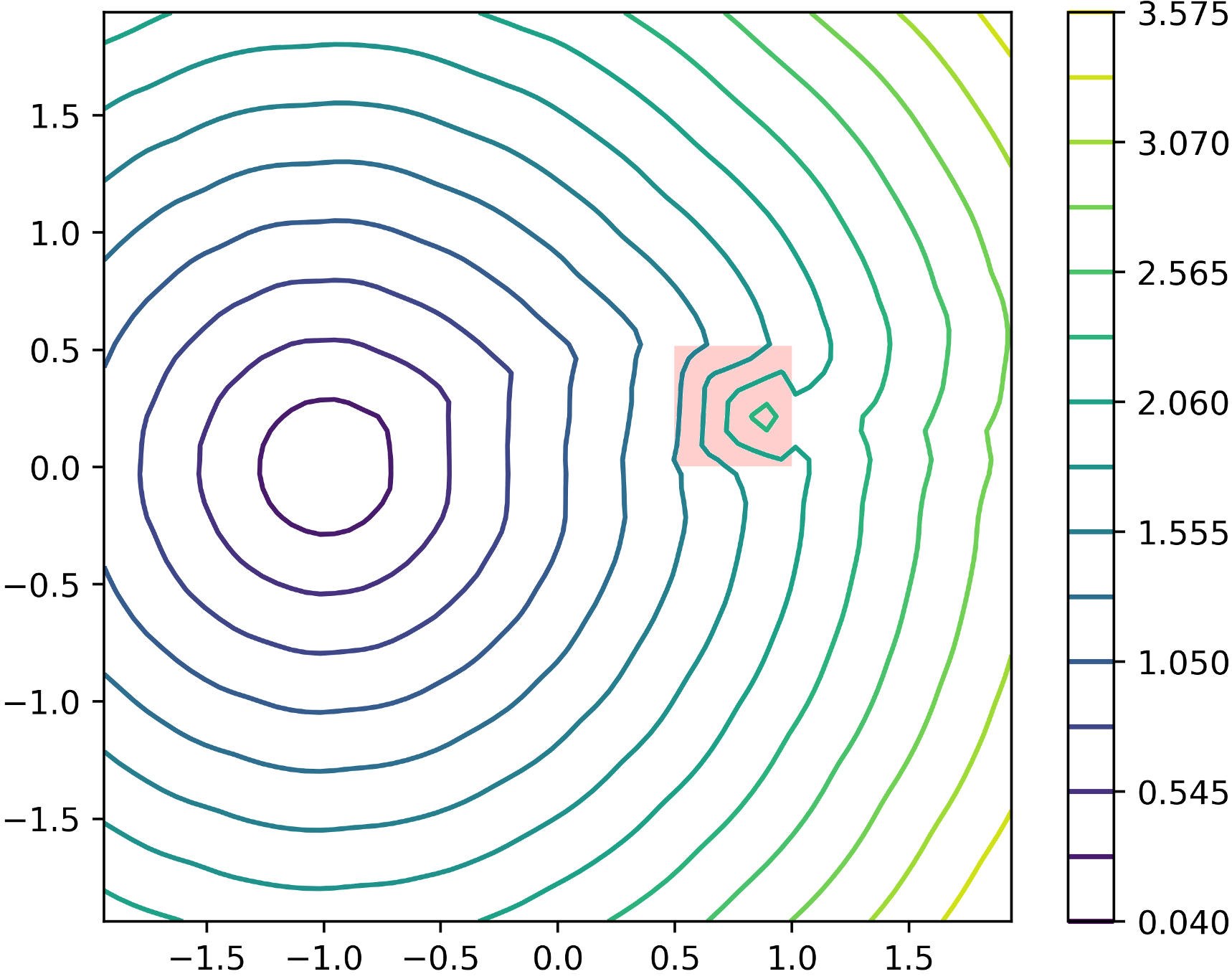}}\hfill
	\subfigure[Ex.~\ref{ex:shape_shading}]{\includegraphics[width=0.3\textwidth]{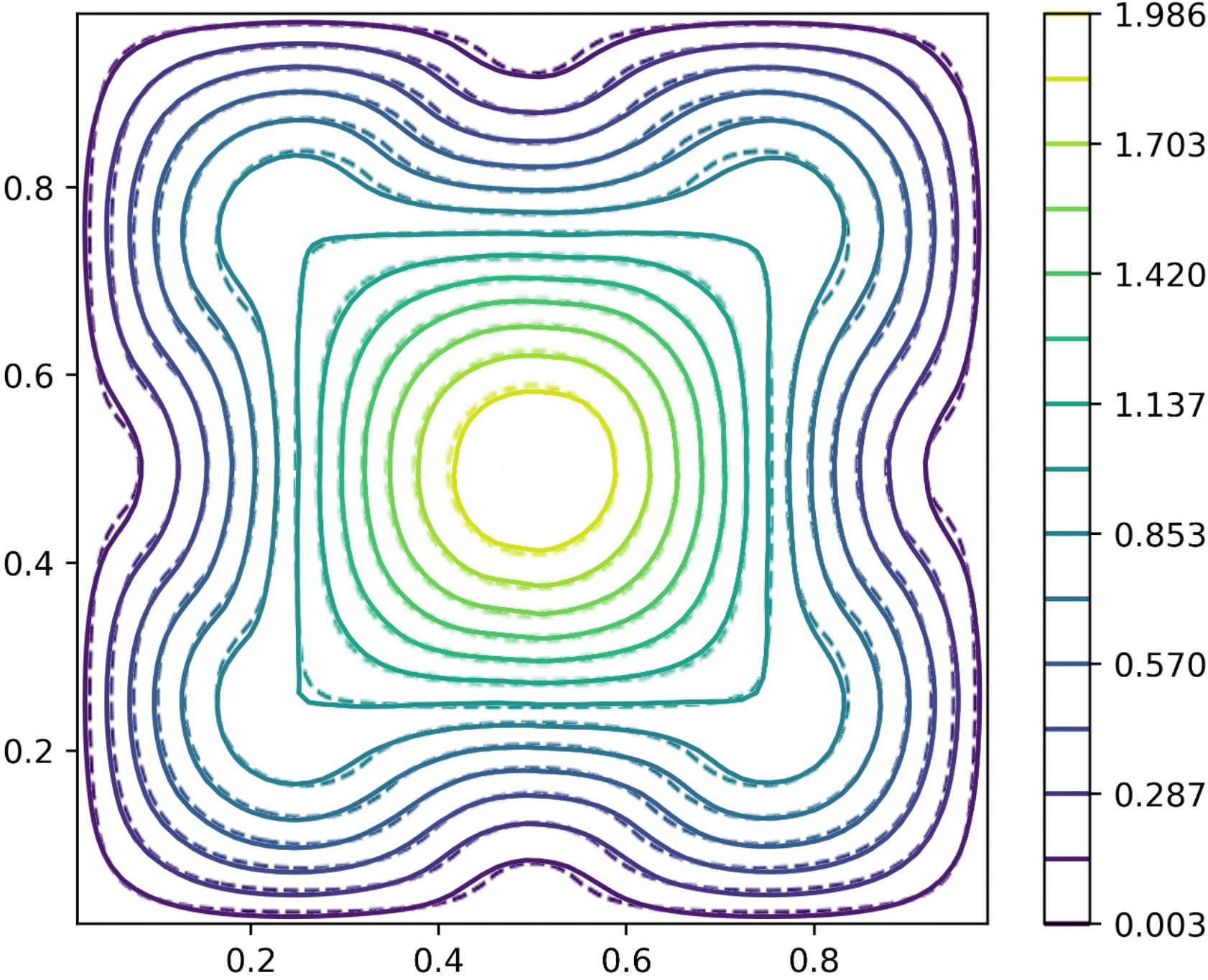}}\\
	\subfigure[Ex.~\ref{ex:sinusoidal}]{\includegraphics[width=0.3\textwidth]{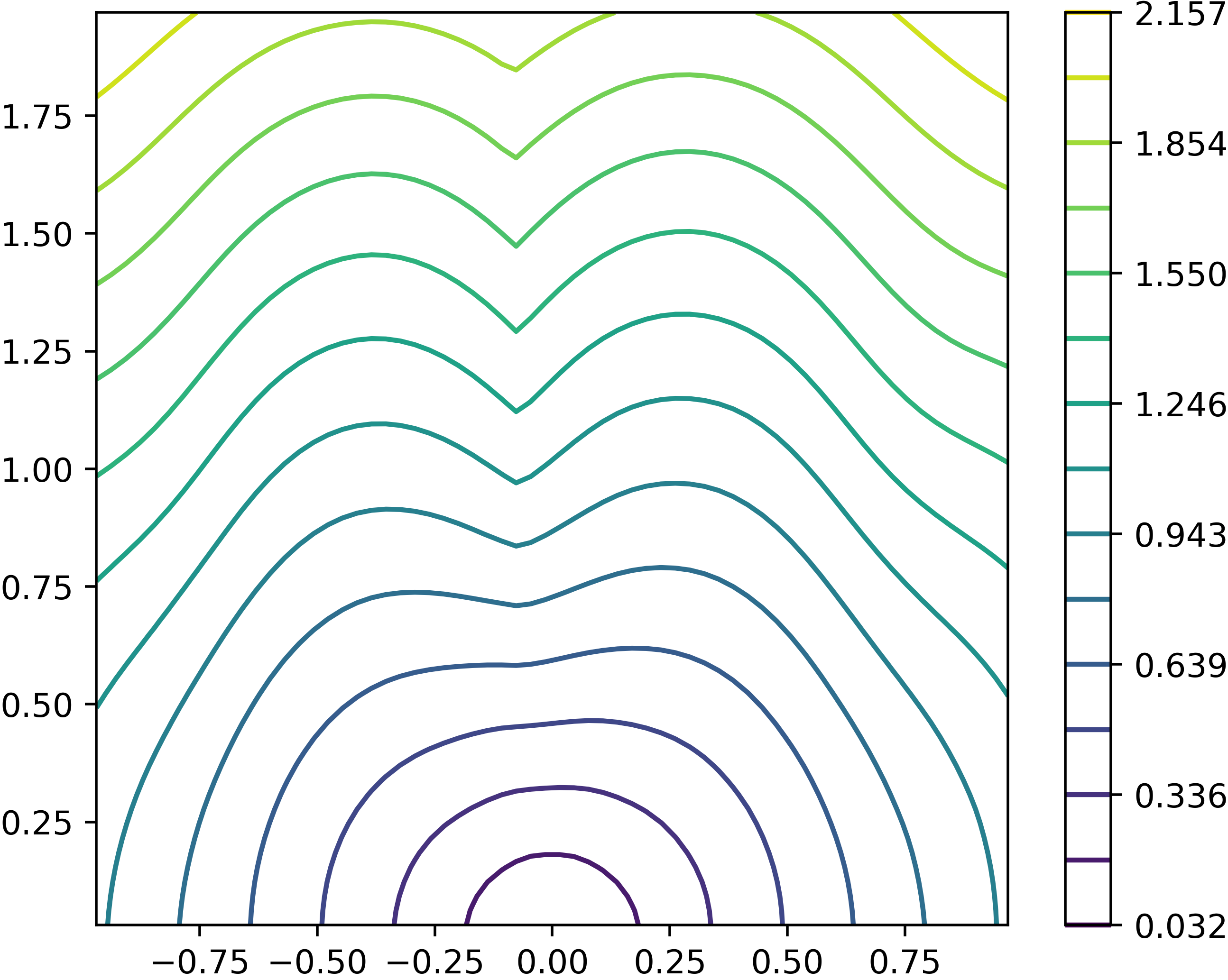}}\hfill
	\subfigure[Ex.~\ref{ex:single_source}]{\includegraphics[width=0.3\textwidth]{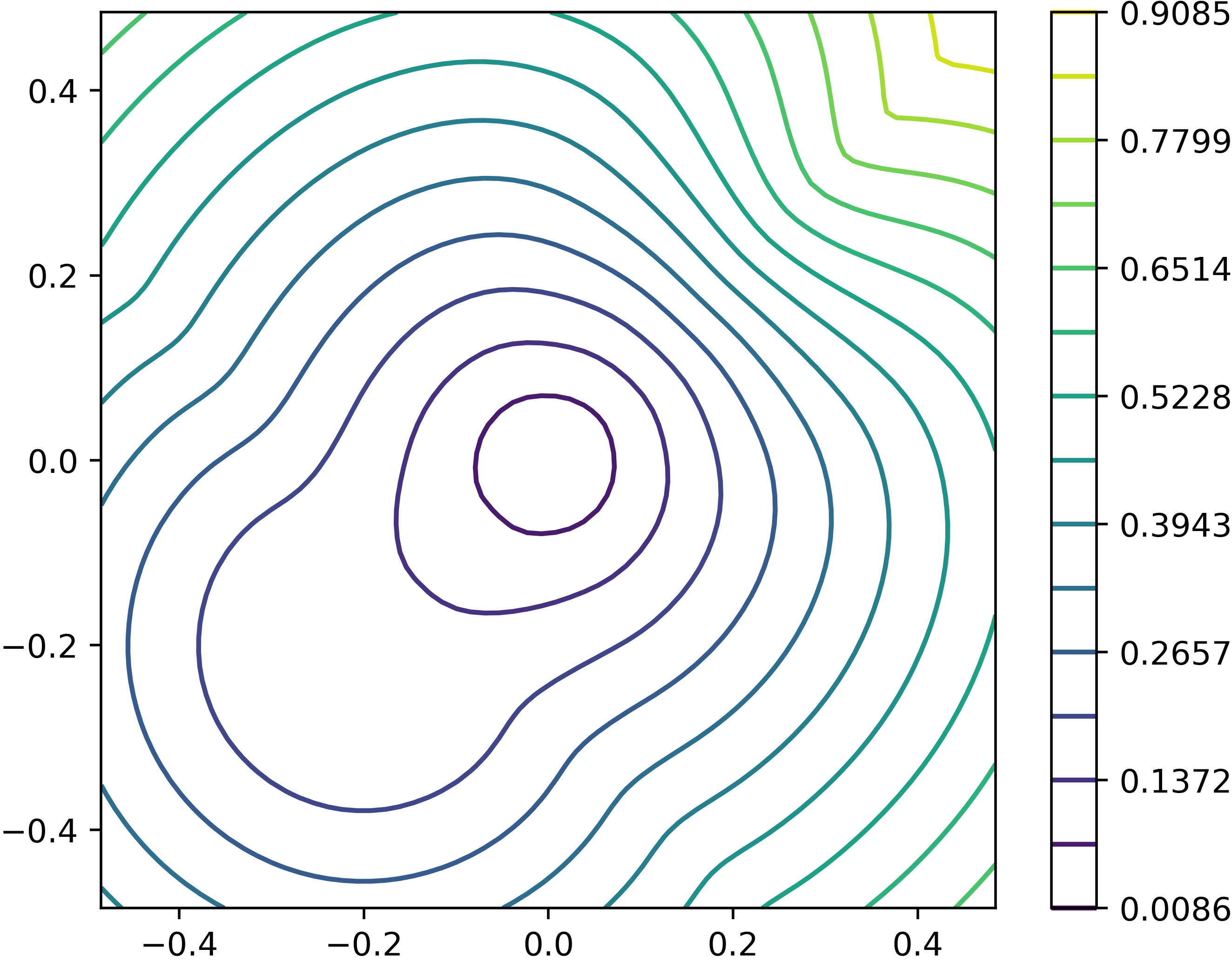}}\hfill
	\subfigure[Ex.~\ref{ex:aniso_ex3_rectangle}]{\includegraphics[width=0.3\textwidth]{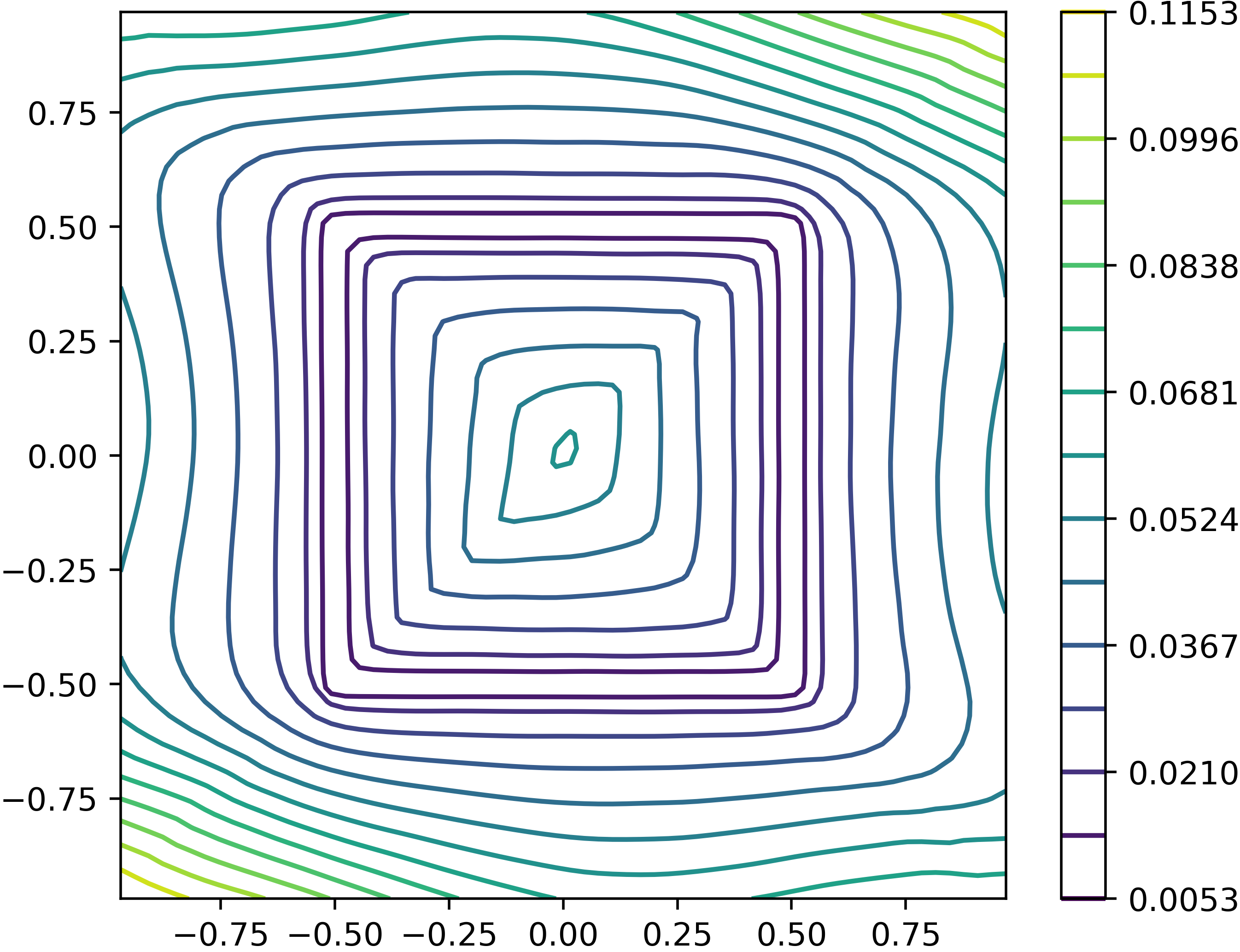}}
	\caption{Contour plots of the predicted solutions. (a), (b) The red-shaded region $[0.5,1]\times[0,0.5]$ is where $f$, and hence $M$, is discontinuous. (c) Predicted (solid) and exact (dashed) solutions.}
	\label{fig:shape_shading}
\end{figure}

\begin{figure}[t]
	\centering
	\subfigure[Example~\ref{ex:sine}]{\label{fig:sine_eps}\includegraphics[width=0.99\linewidth]{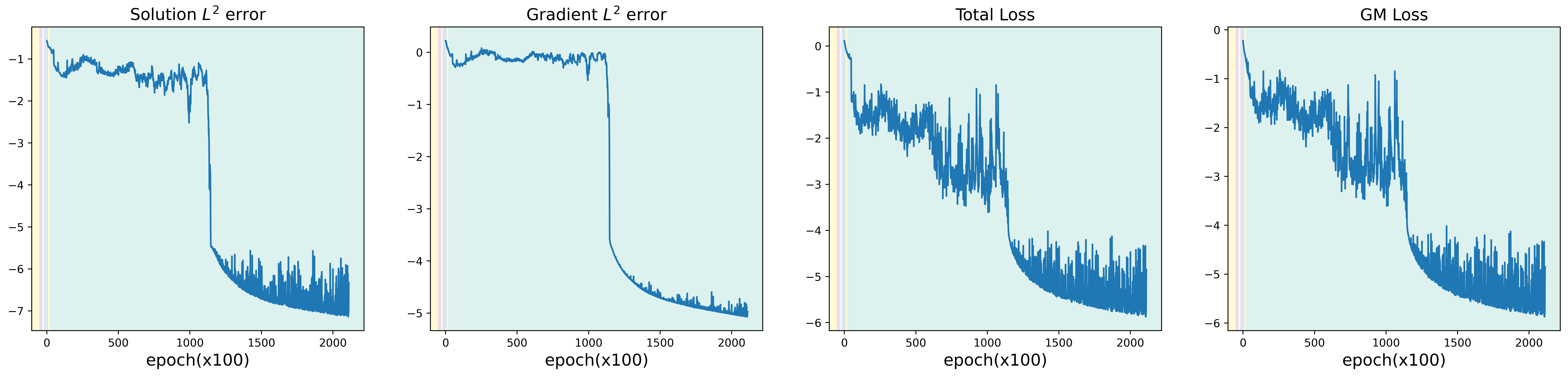}}\\
	\subfigure[Example~\ref{ex:shape_shading}]{\label{fig:err_shape_shading}\includegraphics[width=0.99\linewidth]{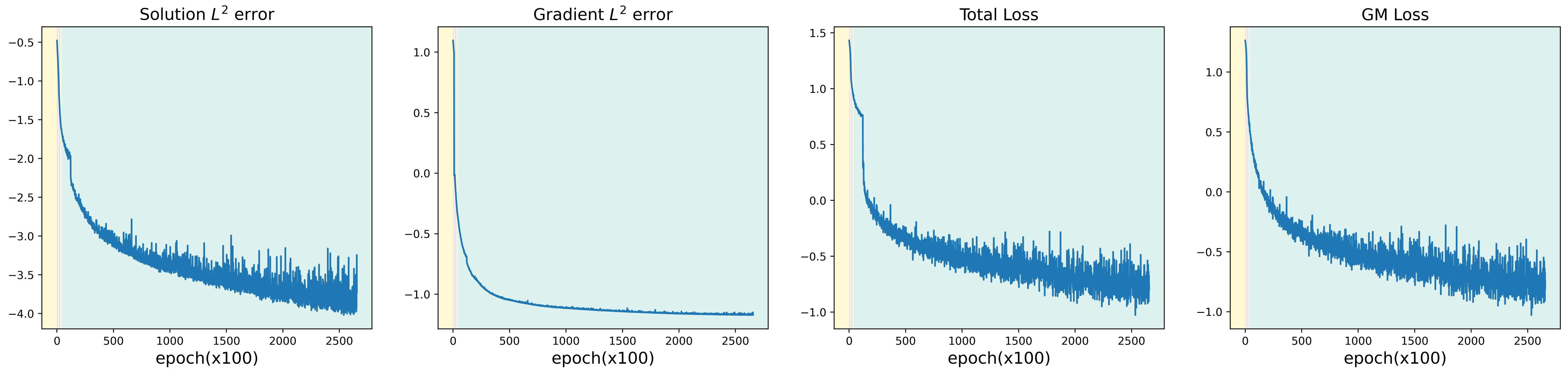}}\\
	\subfigure[Example~\ref{ex:aniso_ex1}]{\label{fig:aniso_ex1}\includegraphics[width=0.99\linewidth]{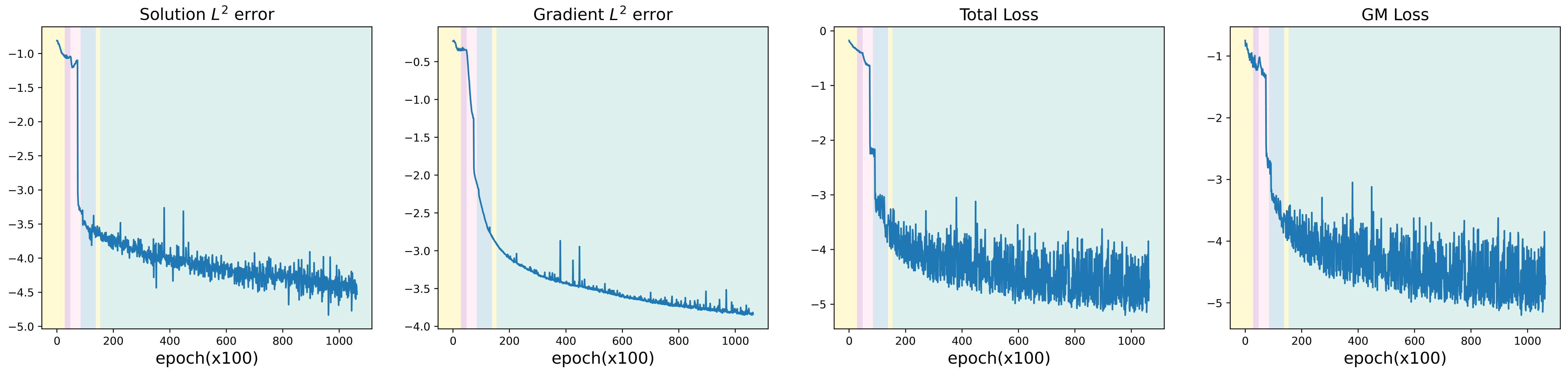}}
	\caption{Solution error, gradient error, total loss and GM loss during training with the annealing factor $1/10$.}
	\label{fig:curves_exact}
\end{figure}

\begin{figure}[t]
	\centering
	\subfigure[Example~\ref{ex:disconti1}]{\label{fig:loss_disconti1}\includegraphics[width=0.48\linewidth]{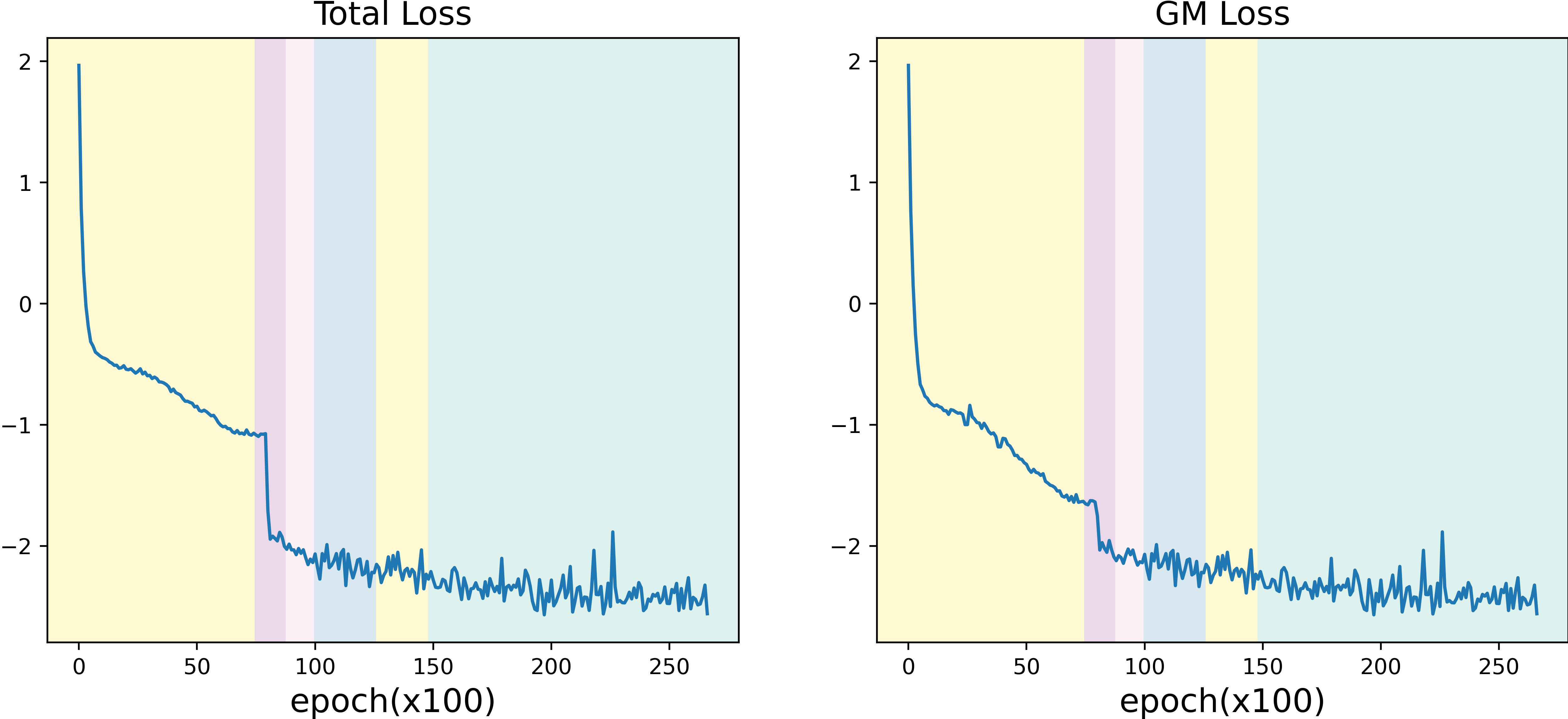}}\hfill
	\subfigure[Example~\ref{ex:sinusoidal}]{\label{fig:loss_sinusoidal}\includegraphics[width=0.48\linewidth]{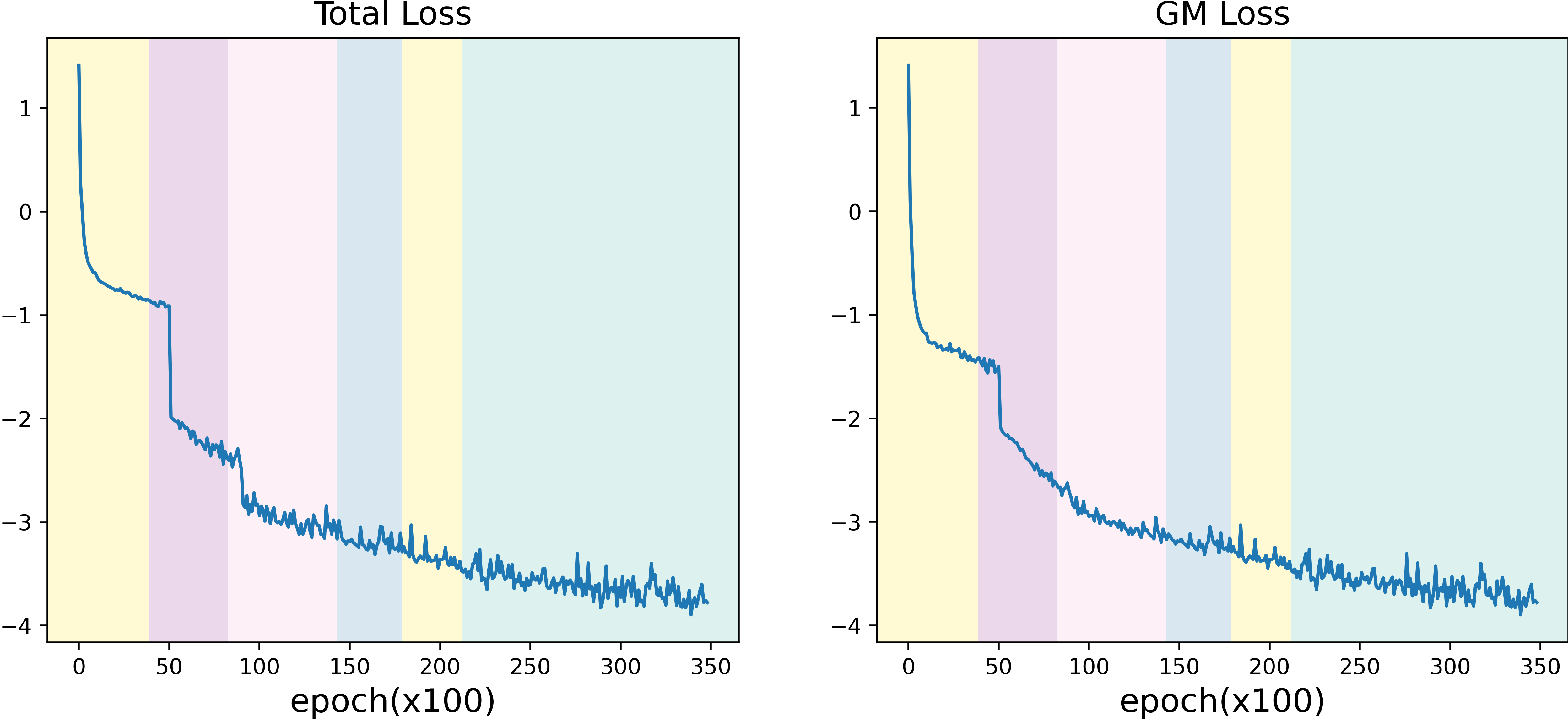}}\\
	\subfigure[Example~\ref{ex:single_source}]{\label{fig:loss_single_source}\includegraphics[width=0.48\linewidth]{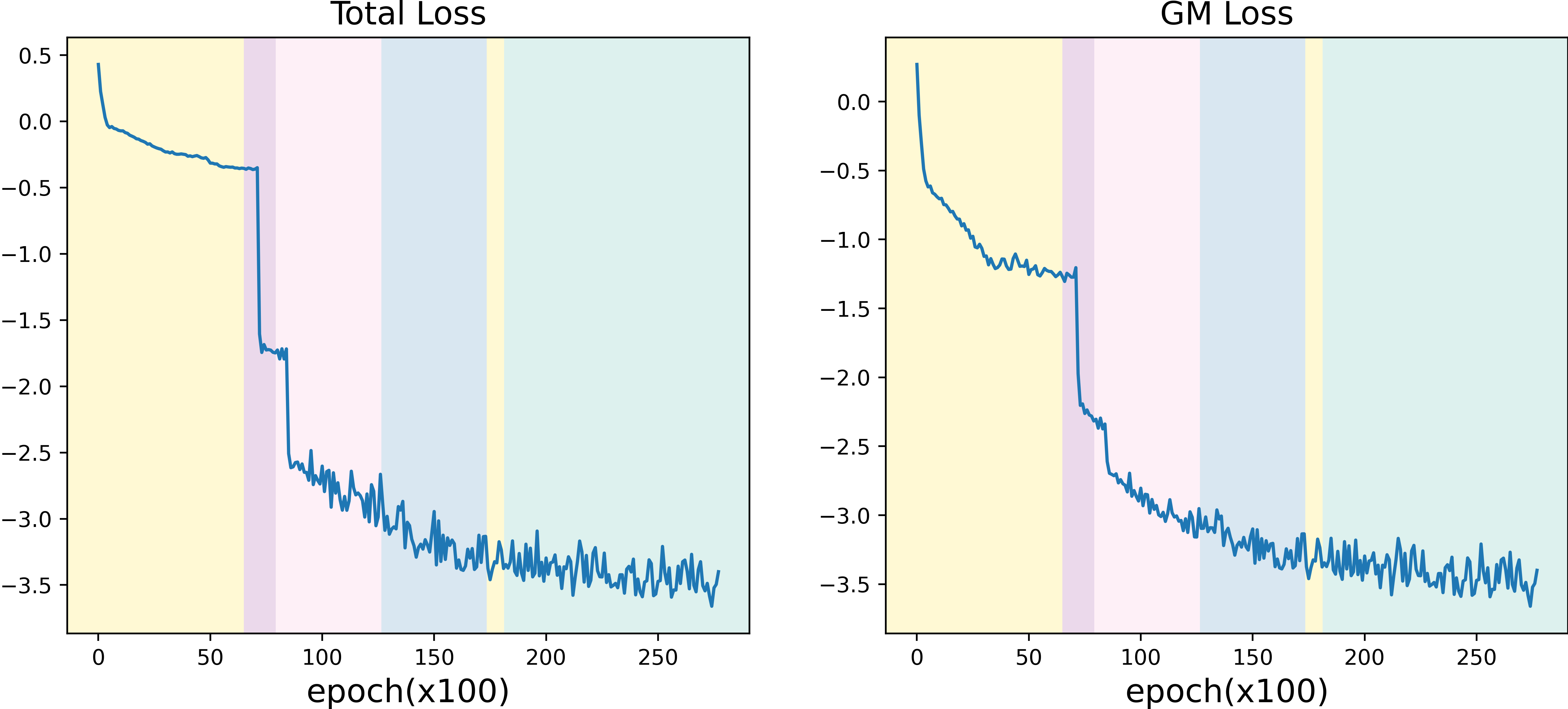}}\hfill
	\subfigure[Example~\ref{ex:aniso_ex3_rectangle}]{\label{fig:loss_aniso_ex3_rectangle}\includegraphics[width=0.48\linewidth]{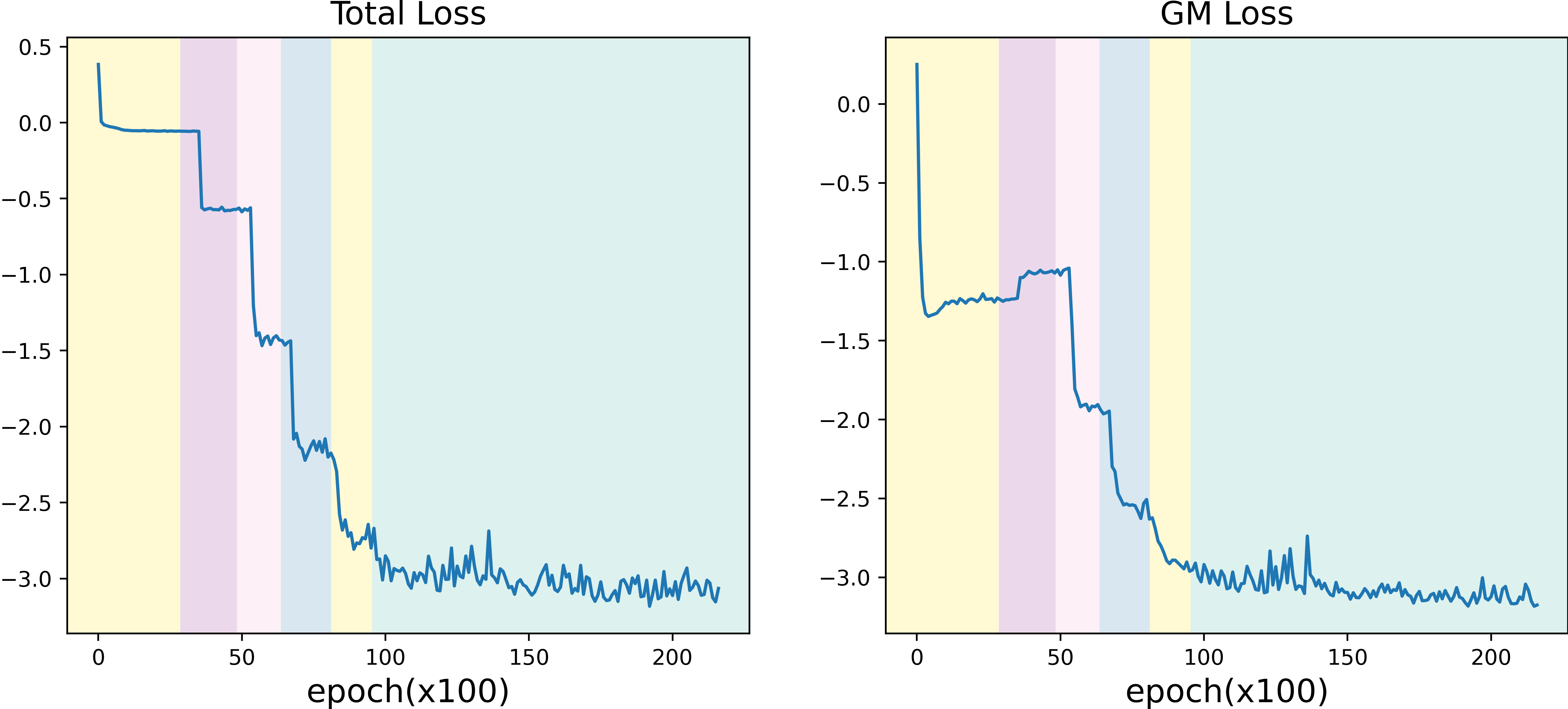}}
	\caption{Total loss and GM loss during training with the annealing factor $1/10$ for the examples without an exact solution.}
	\label{fig:loss_curves}
\end{figure}

\begin{figure}[p]
	\centering
	\subfigure[(A1) with $\lambda=1/\epsilon^2$]{\includegraphics[width=0.99\linewidth]{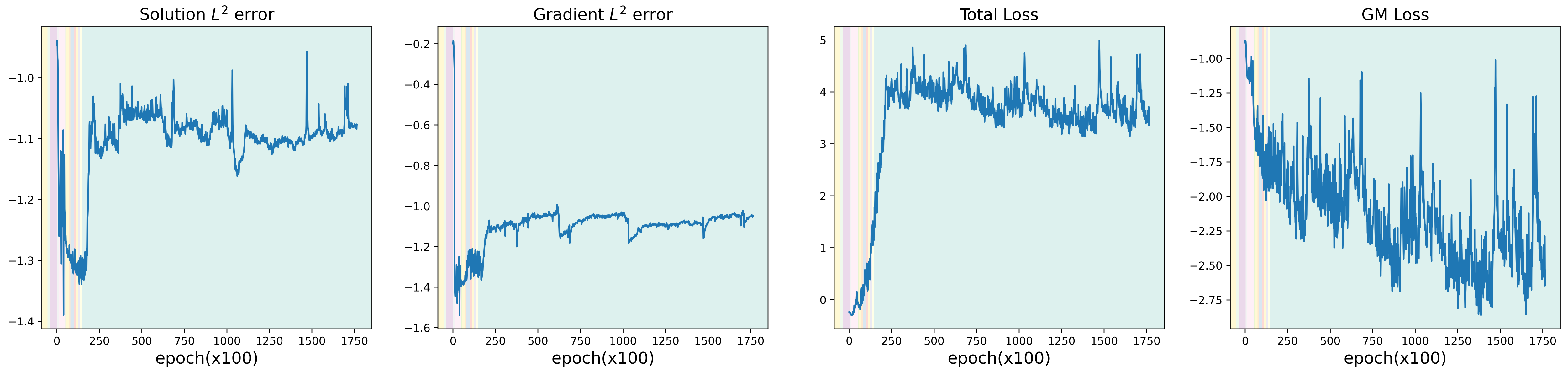}}\\
	\subfigure[(A1) with $\lambda=1/\epsilon$]{\includegraphics[width=0.99\linewidth]{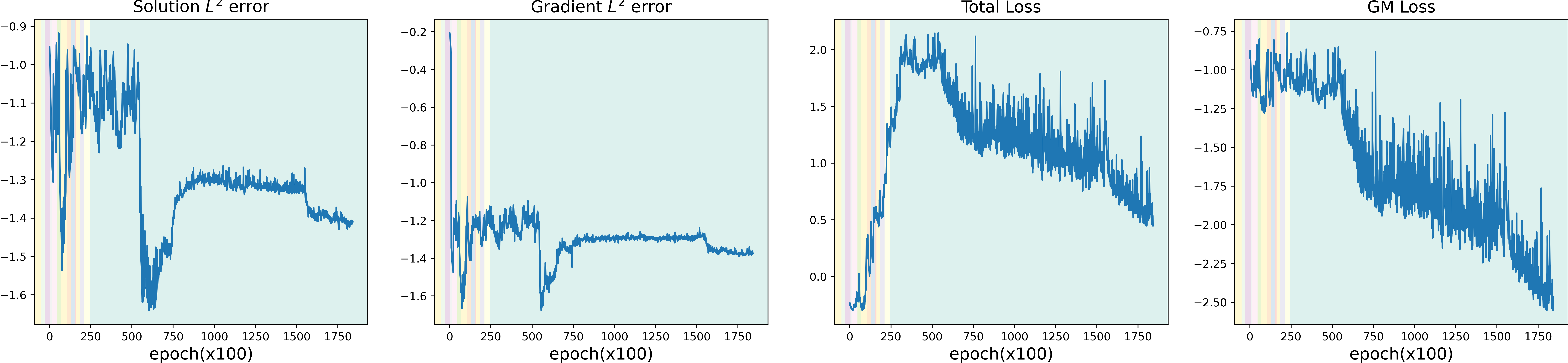}}\\
	\subfigure[(A2): $\epsilon=1$ fixed, $\lambda$ doubled progressively]{\includegraphics[width=0.99\linewidth]{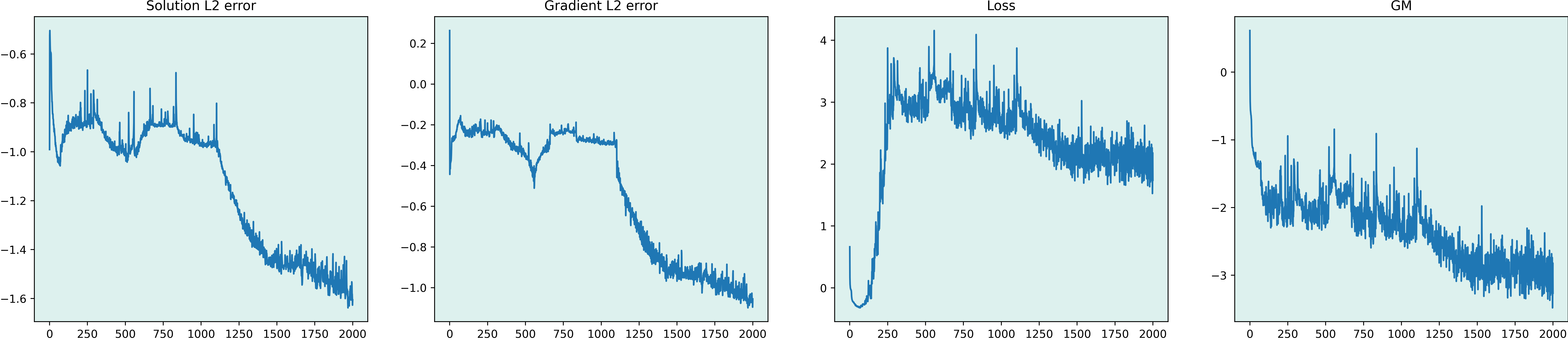}}\\
	\subfigure[(A3): $\lambda=1$ fixed, $\epsilon$ annealed with factor $1/10$]{\includegraphics[width=0.99\linewidth]{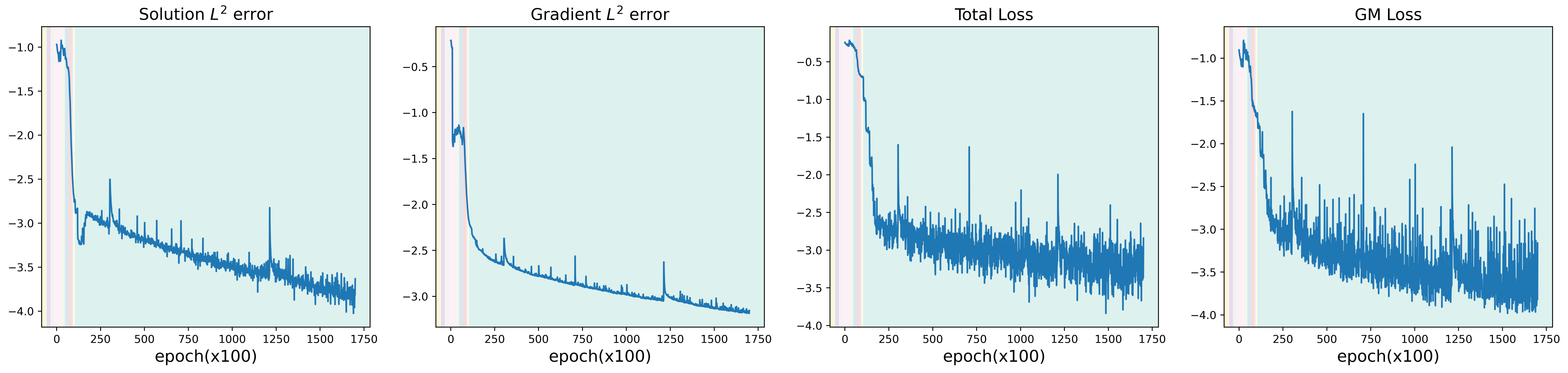}}\\
	\subfigure[(A3): $\lambda=1$ fixed, $\epsilon$ annealed with factor $1/2$]{\includegraphics[width=0.99\linewidth]{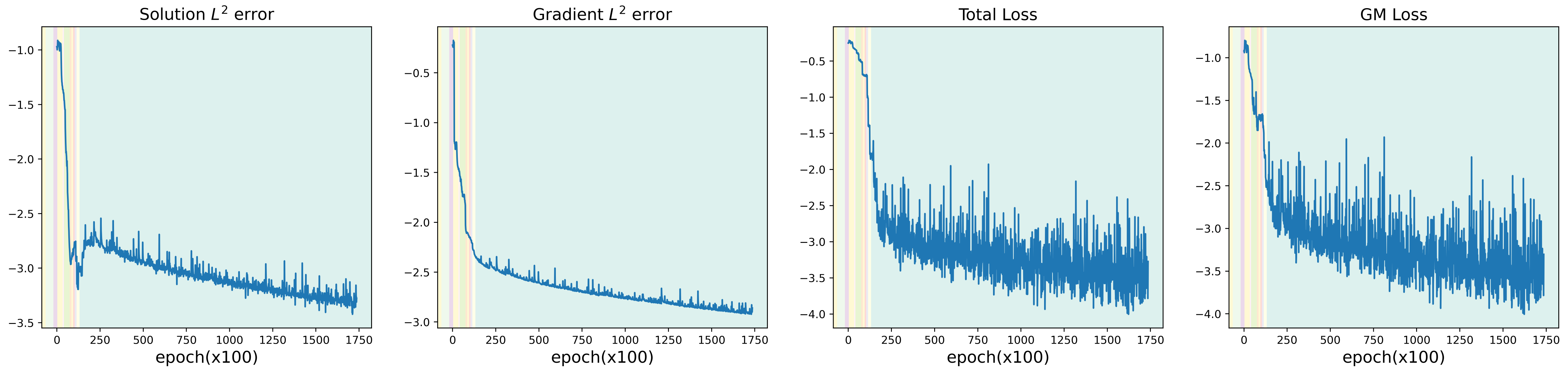}}
	\caption{Comparison of the three asymptotic strategies (A1)--(A3) described in Section~\ref{sec:ablation_lambda} for Example~\ref{ex:aniso_ex2}: solution error, gradient error, total loss, and GM loss during training (logarithmic y-axes). Panels (d) and (e) provide the training curves corresponding to the proposed strategy (A3) considered in Example~\ref{ex:aniso_ex2}, using the two annealing schedules described in Section~\ref{ex:random_pnts}.}
	\label{fig:lambda_aniso1}\label{fig:aniso_ex2}
\end{figure}

\begin{figure}[p]
	\centering
	\subfigure[(A1) with $\lambda=1/\epsilon^2$]{\includegraphics[width=0.99\linewidth]{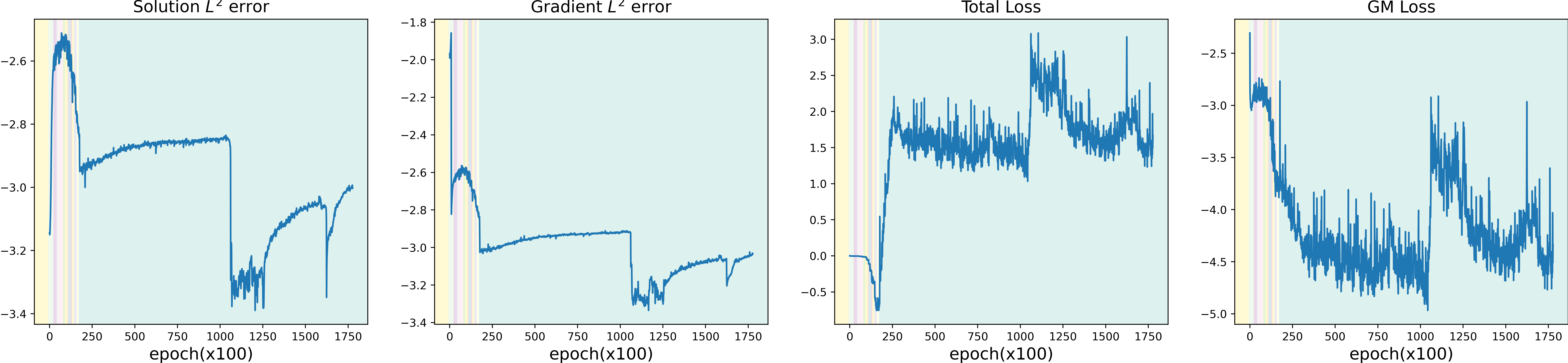}}\\
	\subfigure[(A1) with $\lambda=1/\epsilon$]{\includegraphics[width=0.99\linewidth]{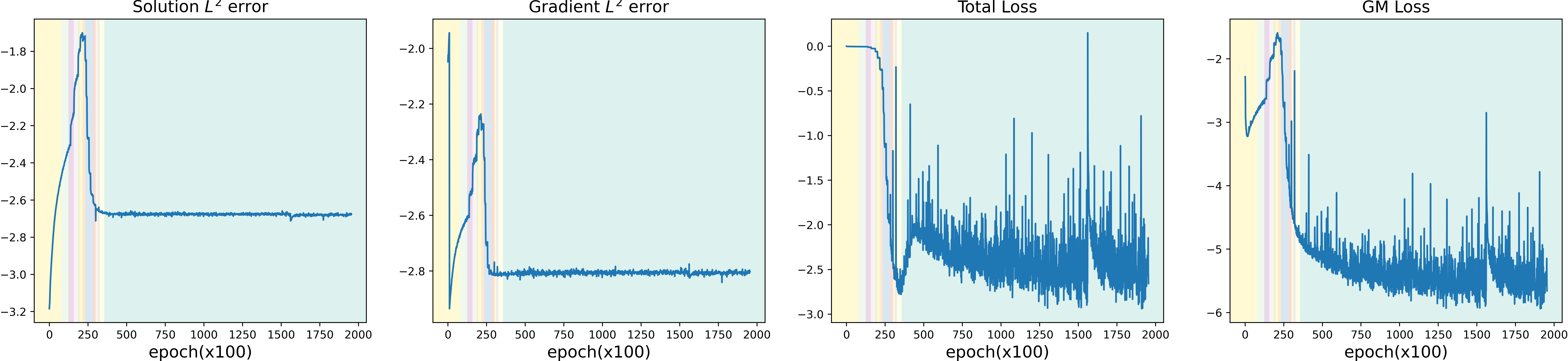}}\\
	\subfigure[(A2): $\epsilon=1$ fixed, $\lambda$ doubled progressively]{\includegraphics[width=0.99\linewidth]{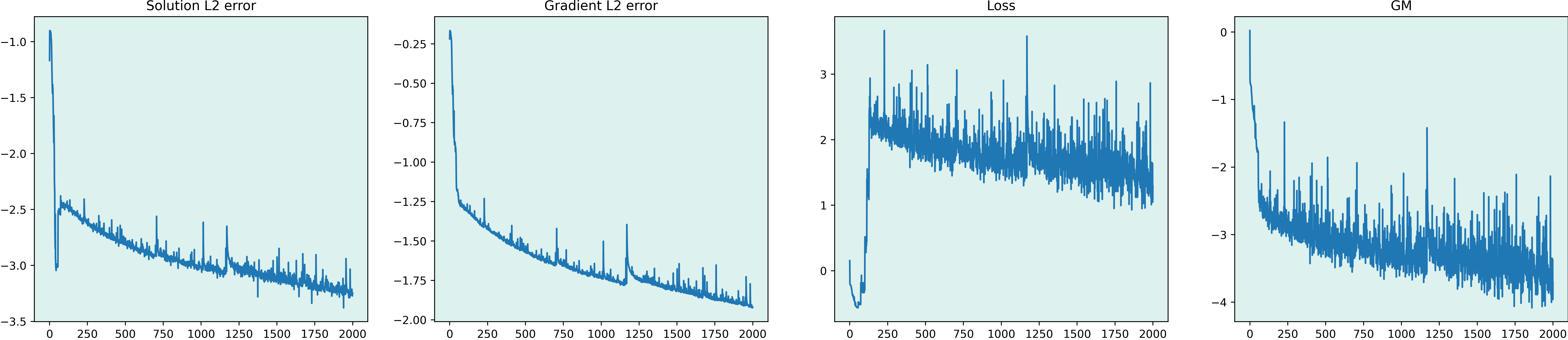}}\\
	\subfigure[(A3): $\lambda=1$ fixed, $\epsilon$ annealed with factor $1/10$]{\includegraphics[width=0.99\linewidth]{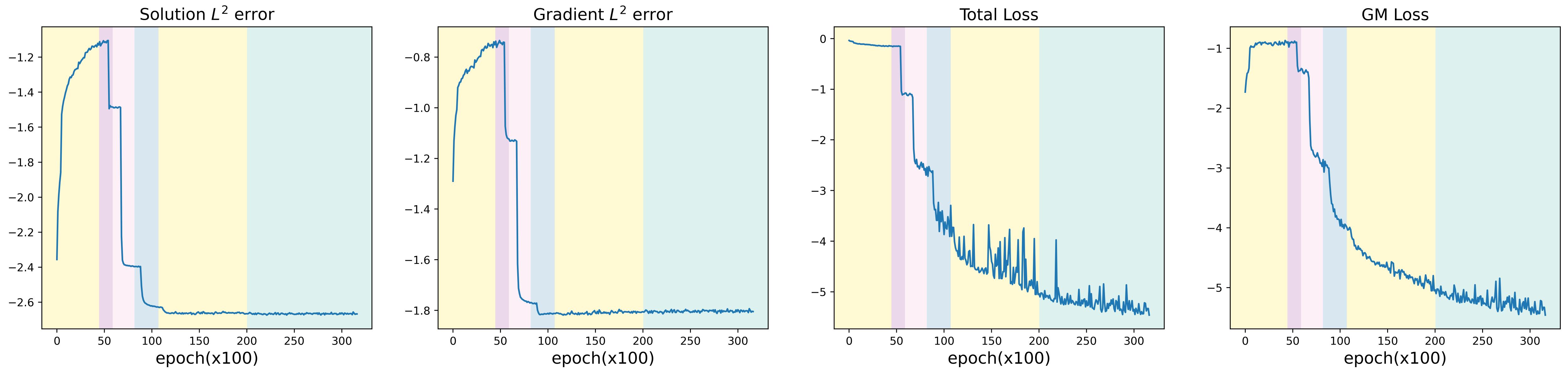}}\\
	\subfigure[(A3): $\lambda=1$ fixed, $\epsilon$ annealed with factor $1/2$]{\includegraphics[width=0.99\linewidth]{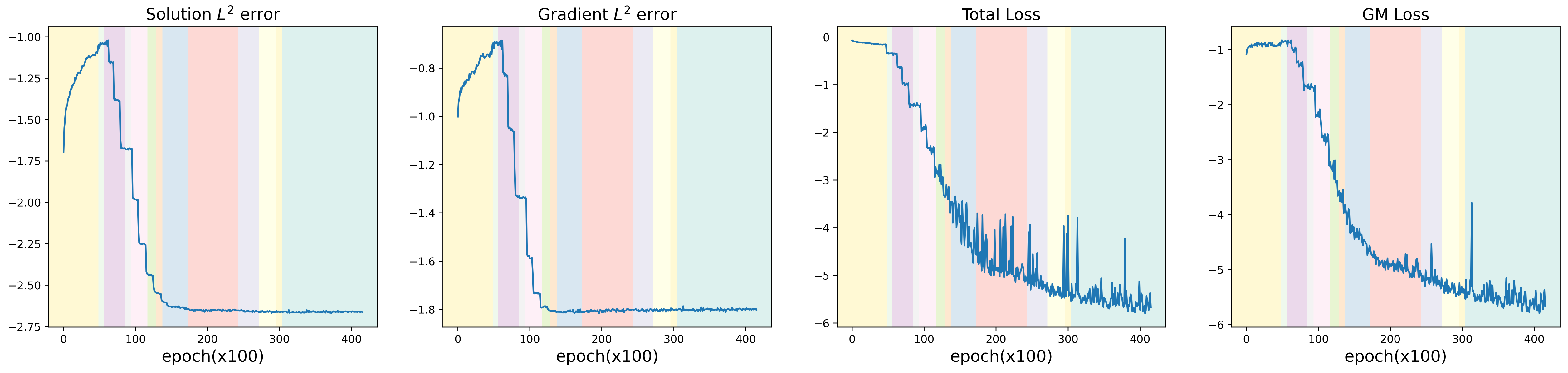}}
	\caption{Comparison of the three asymptotic strategies (A1)--(A3) described in Section~\ref{sec:ablation_lambda} for Example~\ref{ex:aniso_ex3}: solution error, gradient error, total loss and GM loss during training (logarithmic $y$-axes). Panels (d) and (e) provide the training curves corresponding to the proposed strategy (A3) considered in Example~\ref{ex:aniso_ex3}, using the two annealing schedules described in Section~\ref{ex:random_pnts}.}
	\label{fig:lambda_aniso2}\label{fig:aniso_ex3}
\end{figure}

We now vary $M$, $\Gamma$ and $g$ in~\eqref{eq:eikonal}. In the isotropic examples $M=f^{-2}I$, so that~\eqref{eq:eikonal} reads $\|\nabla u\|=f$: $f$ is the slowness, the reciprocal of the speed. Figure~\ref{fig:shape_shading} collects the contour plots, Table~\ref{tab:all_errors} the errors for the examples with an exact or reference solution, and Figures~\ref{fig:curves_exact}, \ref{fig:loss_curves}, \ref{fig:lambda_aniso1} and~\ref{fig:lambda_aniso2} the training curves; Section~\ref{appen:further} adds further curves with the annealing factor $1/2$.

\begin{example}[Discontinuous slowness]\label{ex:disconti1}
Following~\cite{cristiani2007fast}, let $\Omega=[-2,2]^2$, $f=5/2$ in $[0.5,1]\times[0,0.5]$ and $f=1$ otherwise, with $\Gamma=\{(-1,0)\}$ and $g=0$; we also run Ablation~2~\eqref{eq:ablation2}. The proposed method captures the solution across the discontinuity (Figure~\ref{fig:shape_shading}(a)), in agreement with~\cite{cristiani2007fast}. Between the source point and the region of the discontinuity, the solution is the Euclidean distance from the source point; Ablation~2 distorts it there under the influence of the discontinuity (Figure~\ref{fig:shape_shading}(b)), which we attribute to the nonlinearities of~\eqref{eq:ablation2} identified in Section~\ref{ex:random_pnts}.
\end{example}

\begin{example}[Periodic slowness]\label{ex:sine}
For the low-frequency, symmetric periodic slowness
$f(x,y)=\frac{\pi}{2}\sqrt{\sin^2(\pi+\frac{\pi x}{2})+\sin^2(\pi+\frac{\pi y}{2})}$, $\Gamma=\{(0,0)\}$ and $g=0$, the exact solution is $u(x,y)=\cos(\pi+\frac{\pi x}{2})+\cos(\pi+\frac{\pi y}{2})+2$~\cite{hu2024high}. The network approaches it as $\epsilon$ decreases (Figure~\ref{fig:sine_eps}).
\end{example}

\begin{example}[Shape-from-shading]\label{ex:shape_shading}
Shape-from-shading is the inverse problem of recovering the depth of a surface from the brightness variations of an image. We solve the problem of~\cite{rouy1992viscosity}, in the form of~\cite{hu2024high}, modeled by~\eqref{eq:eikonal} on $\Omega=[0,1]^2$ with the symmetric, highly oscillatory periodic slowness
\begin{equation*}
	f(x,y)=2\pi\sqrt{\left[\cos(2\pi x)\sin(2\pi y)\right]^2+\left[\cos(2\pi y)\sin(2\pi x)\right]^2}.
\end{equation*}
Here $\Gamma$ consists of $\partial\Omega$, where $g=0$, the four points $(\frac14,\frac14)$, $(\frac34,\frac34)$, $(\frac34,\frac14)$, $(\frac14,\frac34)$, where $g=1$, and the point $(\frac12,\frac12)$, where $g=2$. The boundary condition is imposed softly through $\cL_\Gamma$ in~\eqref{eq:loss_boundary}. The slowness vanishes at some points, where the speed is infinite and the equation is singular. The exact solution, the shape under vertical lighting, is $u=\max\{u_1,u_2\}$ if $|x+y-1|<\frac12$ and $|x-y|<\frac12$, and $u=u_1$ otherwise, where $u_1(x,y)=|\sin(2\pi x)\sin(2\pi y)|$ and $u_2(x,y)=1+\cos(2\pi x)\cos(2\pi y)$. The prediction with the annealing factor $1/10$ captures its non-smooth structure (Figure~\ref{fig:shape_shading}(c)) and approaches the viscosity solution as $\epsilon$ decreases (Figure~\ref{fig:err_shape_shading}).
\end{example}

\begin{example}[Sinusoidal model]\label{ex:sinusoidal}
To test a singular source point and a non-smooth solution, we solve the isotropic wave front of~\cite{hu2024high} with $f(x,y)=1/\{1+0.2\sin(0.5\pi y)\sin(3\pi(x+0.55))\}$, an asymmetric periodic slowness that oscillates at high frequency in $x$ and varies at low frequency in $y$, with $\Gamma=\{(0,0)\}$ and $g=0$. The singular source point and the non-smoothness are captured (Figure~\ref{fig:shape_shading}(d)).
\end{example}

\begin{example}[Exponential decay slowness]\label{ex:single_source}
Following~\cite{zhao2005fast}, let $\Omega=[-0.5,0.5]^2$ and $f(x,y)=1+e^{-40((x+0.2)^2+(y+0.2)^2)}-e^{-40((x-0.2)^2+(y-0.2)^2)}$, which varies exponentially around two points distinct from the source point, $g=0$ on $\Gamma=\{(0,0)\}$. The result (Figure~\ref{fig:shape_shading}(e)) is consistent with~\cite{zhao2005fast}.
\end{example}

\begin{example}[Constant anisotropic]\label{ex:aniso_ex1}
For the constant metric $M=\operatorname{diag}(a,b)$ with $a=2$ and $b=1$, $\Gamma=\{(0,0)\}$ and $g=0$, the exact solution is $u(x,y)=\sqrt{x^2/a+y^2/b}$. As $\epsilon$ decreases the network approaches it (Figure~\ref{fig:aniso_ex1}).
\end{example}

\begin{example}[Function-scaled constant anisotropic]\label{ex:aniso_ex2}
Following~\cite{luo2012fast}, let
\begin{equation*}
	M(x,y)=\frac{1}{f(x,y)^2}\begin{pmatrix}1&-1\\-1&2\end{pmatrix},\qquad
	f(x,y)=\exp\left(-\sqrt{2\bar x^2+2\bar x\bar y+\bar y^2}\right),
\end{equation*}
with $\bar x=x-x_0$, $\bar y=y-y_0$ and $x_0=y_0=0.5$ in $\Omega=[0,1]^2$, so that the magnitude varies in space while the directional bias is fixed, and $g=0$ on $\Gamma=\{(x_0,y_0)\}$. The exact solution is $u=1-f$. The method learns the viscosity solution as $\epsilon$ decreases, with stable training (Figure~\ref{fig:lambda_aniso1}(d),(e)).
\end{example}

\begin{example}[Fully anisotropic]\label{ex:aniso_ex3}
Following~\cite{luo2012fast}, let $M=(m_{ij})$ with
\begin{equation}\label{eq:aniso_ex_coeff}
	\begin{aligned}
	m_{11}&=150.25\left(1+\mu_1\sin^2(\pi xy)\right), \qquad m_{22}=50.75\left(1+\mu_2\cos^2(\pi xy)\right),\\
	m_{12}&=m_{21}=-86.16953\left(1-\mu_3\sin^2(\pi xy)\right),
	\end{aligned}
\end{equation}
with $\mu_1=1$, $\mu_2=1$ and $\mu_3=0.125$ in $\Omega=[0,1]^2$. Unlike in Examples~\ref{ex:aniso_ex1} and~\ref{ex:aniso_ex2}, both the magnitude and the orientation of the anisotropy vary in space. We set $g=0$ on $\Gamma=\{(x_0,y_0)\}=\{(0.5,0.5)\}$. As the reference solution we use the numerical solution of~\cite{luo2012fast} with the additive factor
\begin{equation*}
	\sqrt{\frac{50.75(1+0.5\mu_2)\bar{x}^2+86.16953(2-\mu_3)\bar{x}\bar{y}+150.25(1+0.5\mu_1)\bar{y}^2}{150.25(1+0.5\mu_1)\,50.75(1+0.5\mu_2)-\bigl(86.16953(1-0.5\mu_3)\bigr)^2}},
\end{equation*}
with $\bar x=x-x_0$ and $\bar y=y-y_0$. The network approaches the viscosity solution as $\epsilon$ decreases (Figure~\ref{fig:lambda_aniso2}(d),(e)).
\end{example}

\begin{example}[Fully anisotropic, square interface]\label{ex:aniso_ex3_rectangle}
With $M$ in~\eqref{eq:aniso_ex_coeff}, let $\Gamma$ be the boundary of the square of side length~$1$ centered at the origin, with $g=0$ and $\phi(x,y)=\left|\max\{|x|-\frac12,|y|-\frac12\}\right|$ in~\eqref{eq:eik_primary_variable}. The prediction (Figure~\ref{fig:shape_shading}(f)) agrees with the result of~\cite{qian2007fast}.
\end{example}

\begin{example}[Surface reconstruction from 3D point clouds]\label{ex:surface_recon}
Reconstructing open surfaces from point clouds is a challenging task in computer vision. An open surface can be represented implicitly by the zero level set of the distance function to the point cloud, which the proposed method computes without modification: $M=I$, $g=0$, and the given points form the set $m_\Gamma$ of boundary collocation points for $\cL_\Gamma$ in~\eqref{eq:loss_boundary}, with $\lambda_\Gamma=1$. The point cloud is centered at its bounding box center and uniformly scaled such that its largest dimension fits within $\Omega=[-1,1]^3$ \cite{sitzmann2020implicit}. No mesh of the domain or of the surface is needed for training; a triangular mesh is extracted from the learned distance function afterwards by the algorithm of~\cite{park2025nsp}. We use six hidden layers of $512$ nodes, $|m_\Omega|=|m_\Gamma|=20{,}000$ collocation points per epoch, $100{,}000$ epochs, $N_\epsilon=1000$, and a learning-rate decay by a factor of $0.99$ every $2{,}000$ epochs. As surface reconstruction is not our focus, we test only two car models from ShapeNet~\cite{shapenet2015}, whose multi-layered, non-closed shapes make the point set $\Gamma$ intricate. In Figure~\ref{fig:car}, the wheels, side mirrors and front bumper are well reconstructed, the surface is closely aligned with the point cloud, and the interior view shows that the seats and the steering wheel are captured. This shows that the method applies to 3D problems and to surface reconstruction and highlights the advantages of its mesh-free formulation in handling complex geometries.

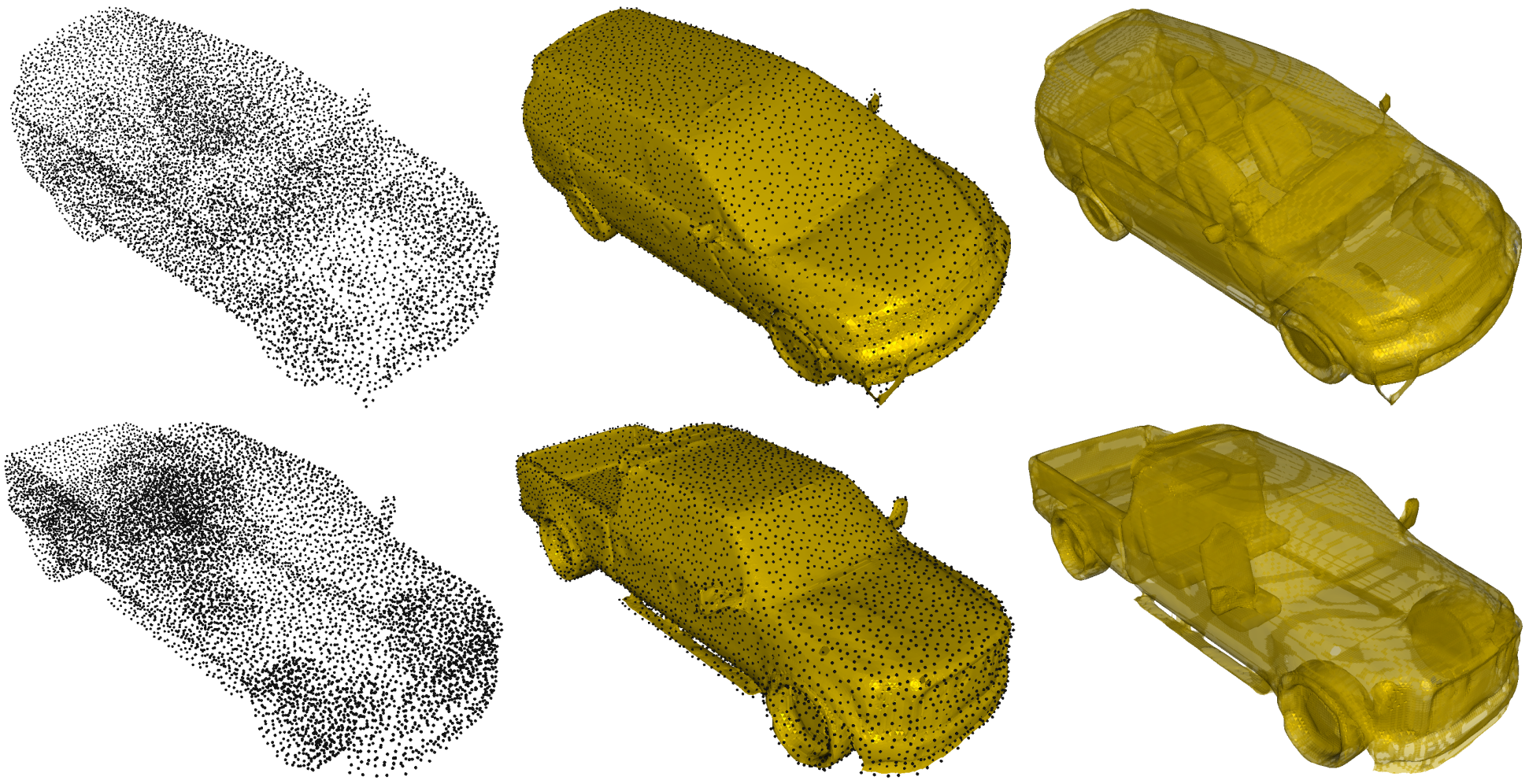
\begin{figure}[t]
	\begin{tikzpicture}[every node/.style={font=\normalsize}]
		\node[anchor=south west,inner sep=0] (image1) at (0,0)
		{\includegraphics[width=0.33\linewidth]{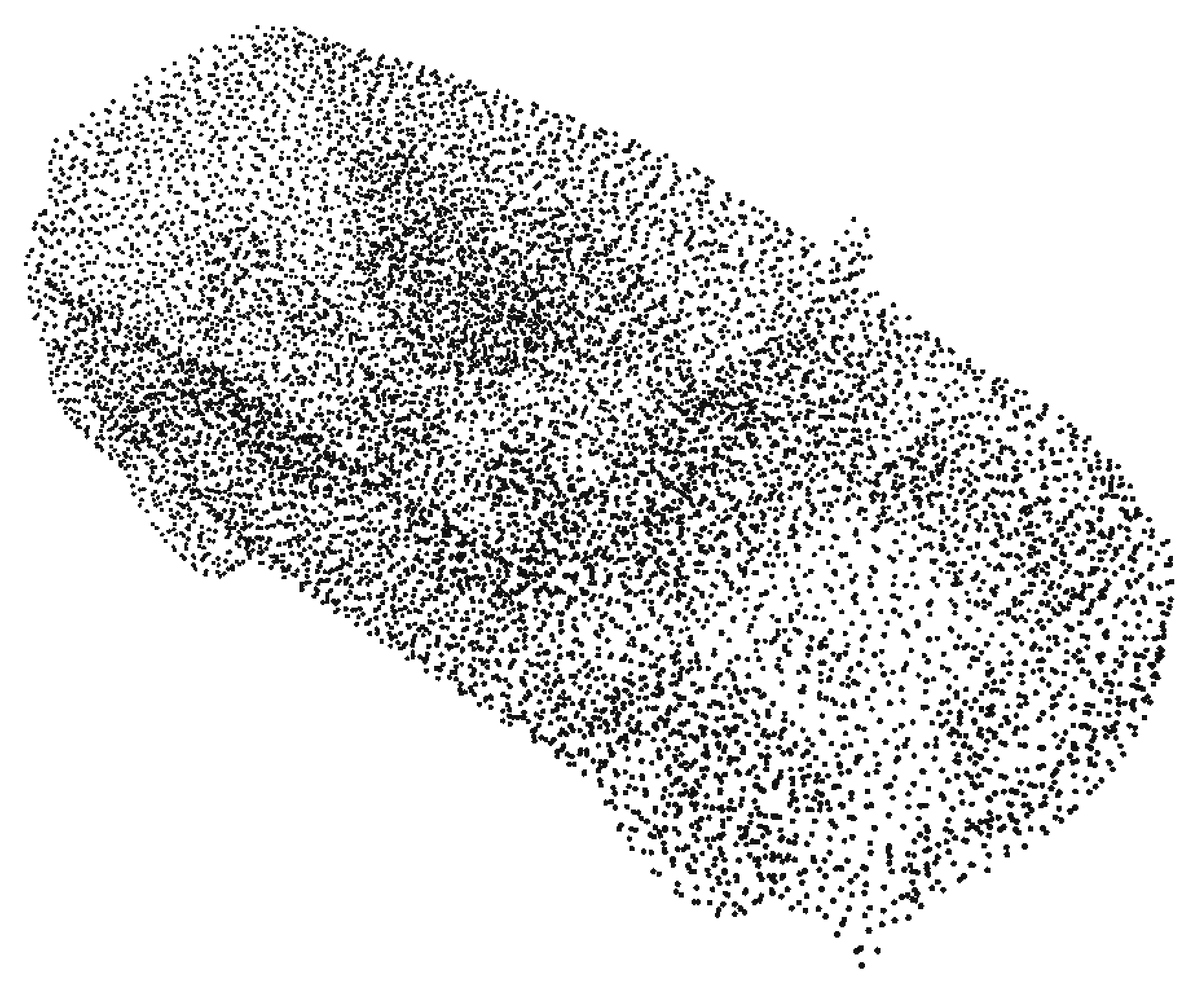}};
		\node[anchor=south west,inner sep=0] (image2) at (0.33\linewidth,0)
		{\includegraphics[width=0.33\linewidth]{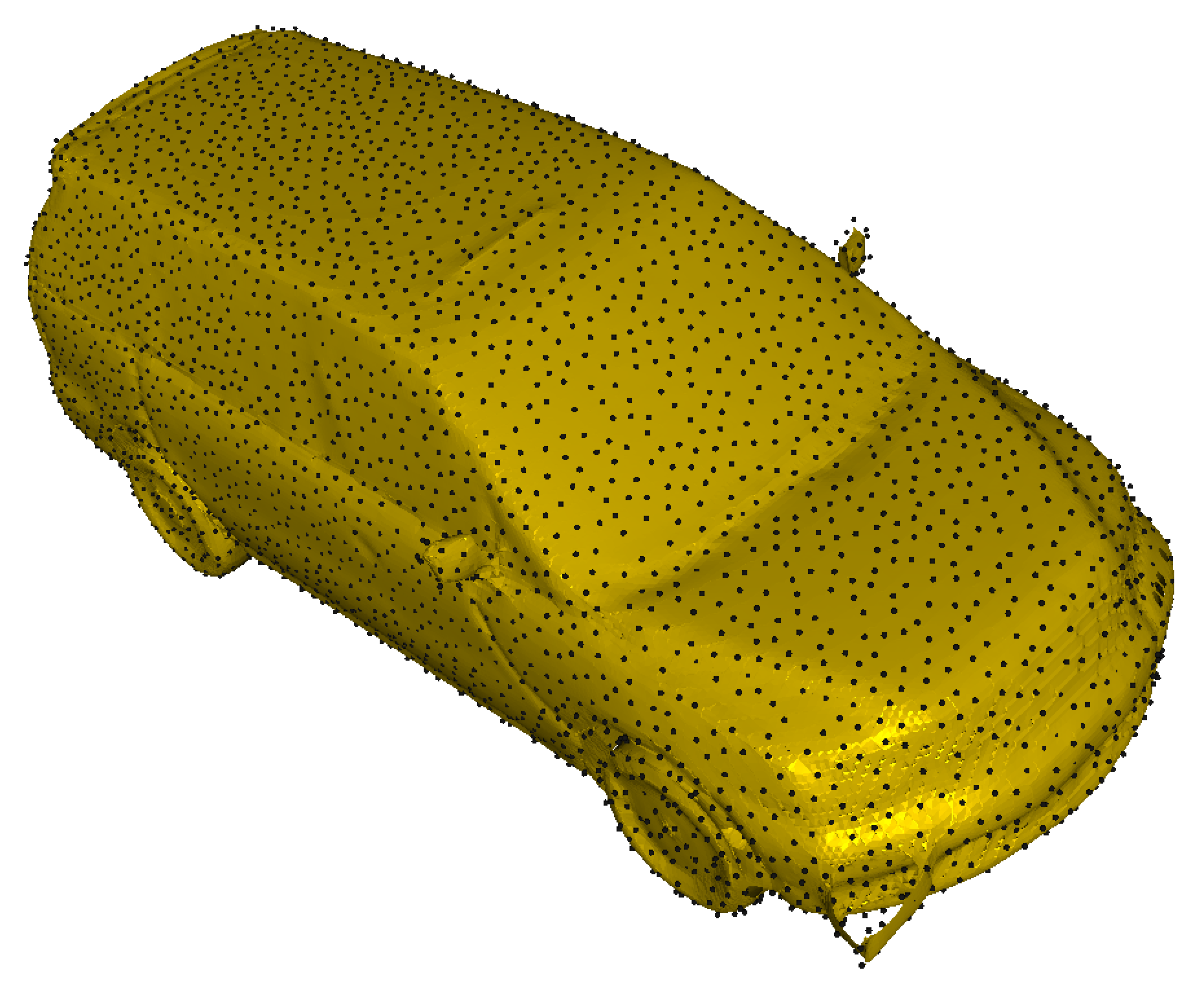}};
		\node[anchor=south west,inner sep=0] (image3) at (0.66\linewidth,0)
		{\includegraphics[width=0.33\linewidth]{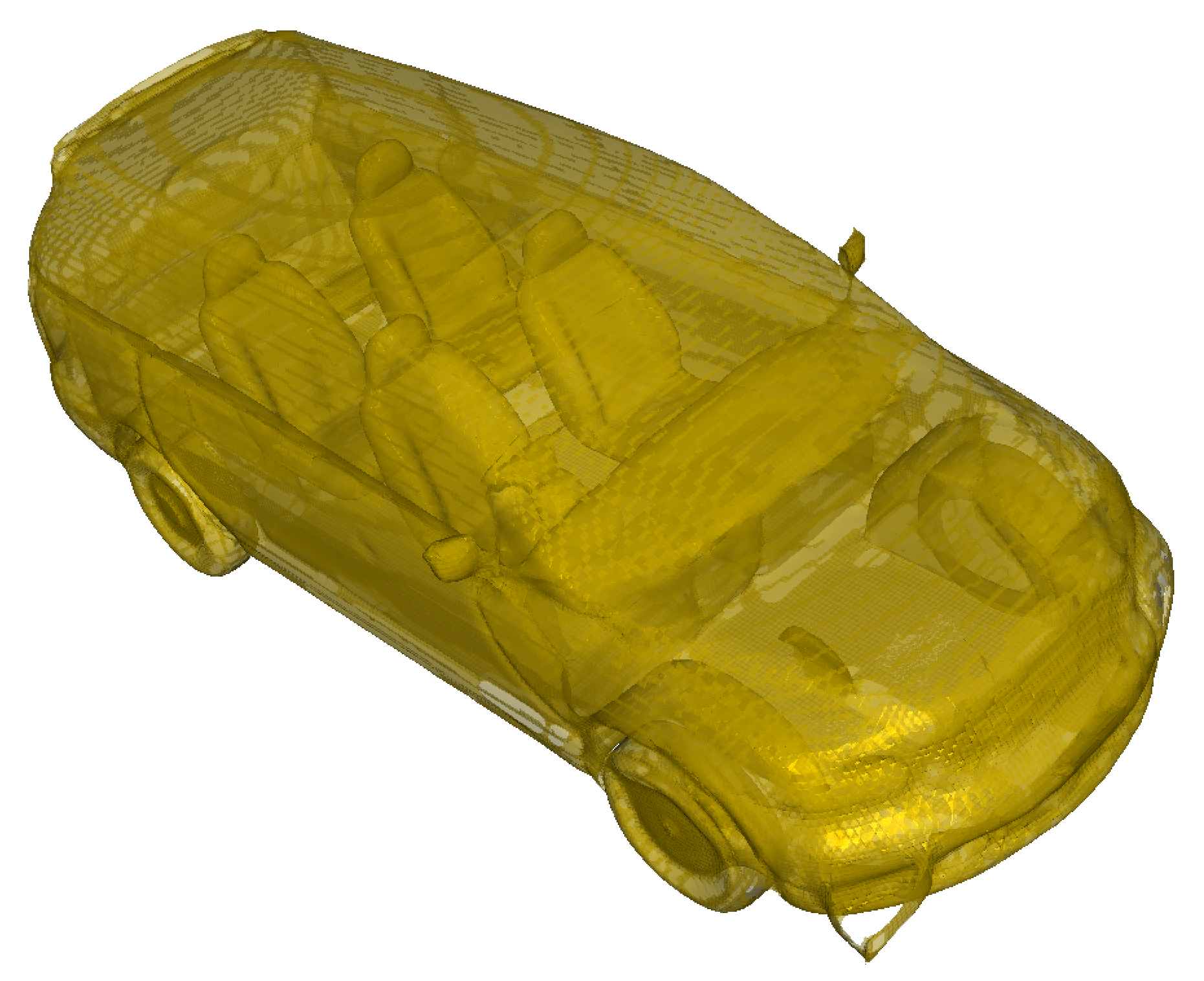}};
		\node[anchor=north west,inner sep=0] (image4) at (image1.south west)
		{\includegraphics[width=0.33\linewidth]{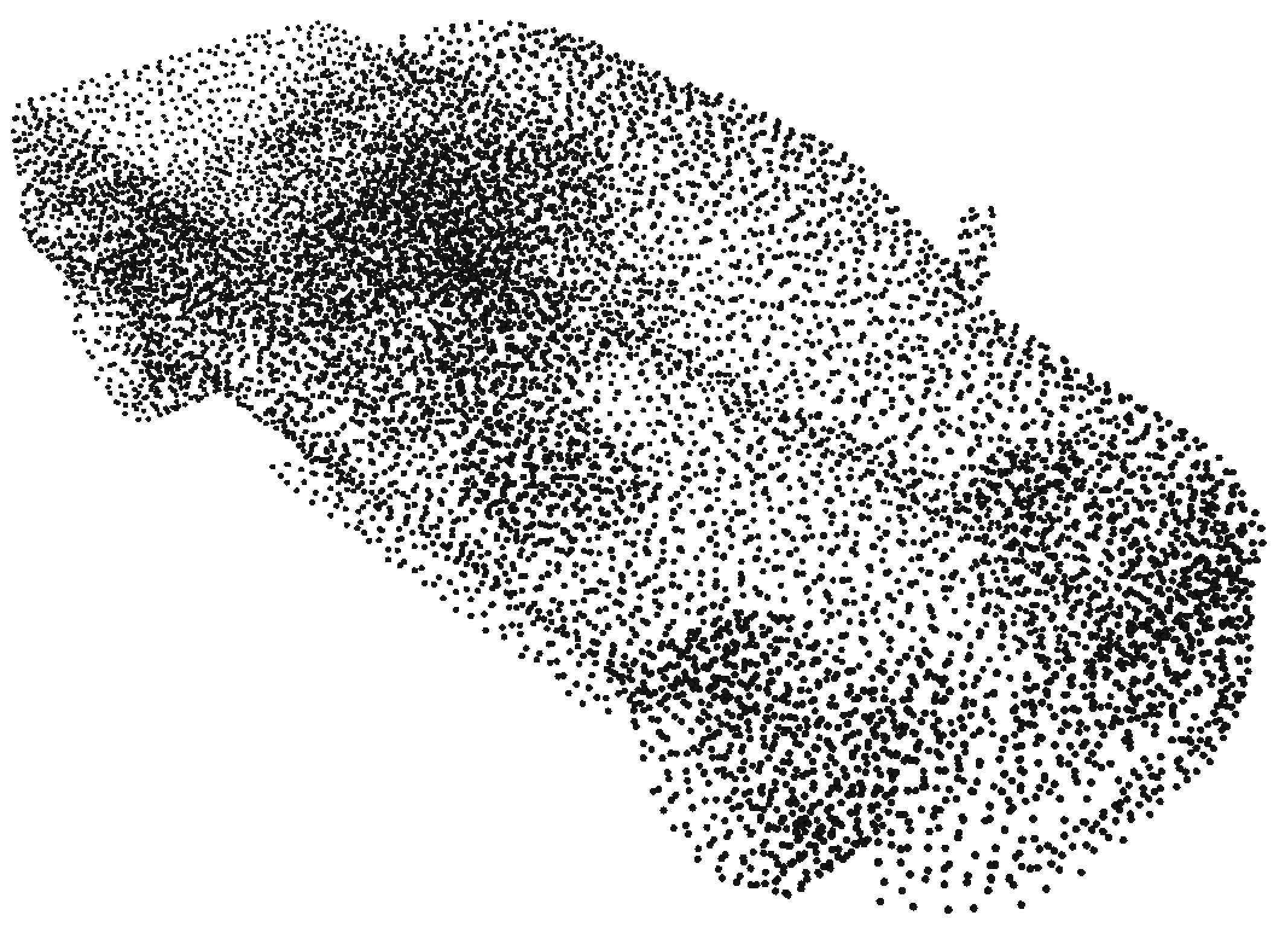}};
		\node[anchor=north west,inner sep=0] (image5) at (image2.south west)
		{\includegraphics[width=0.33\linewidth]{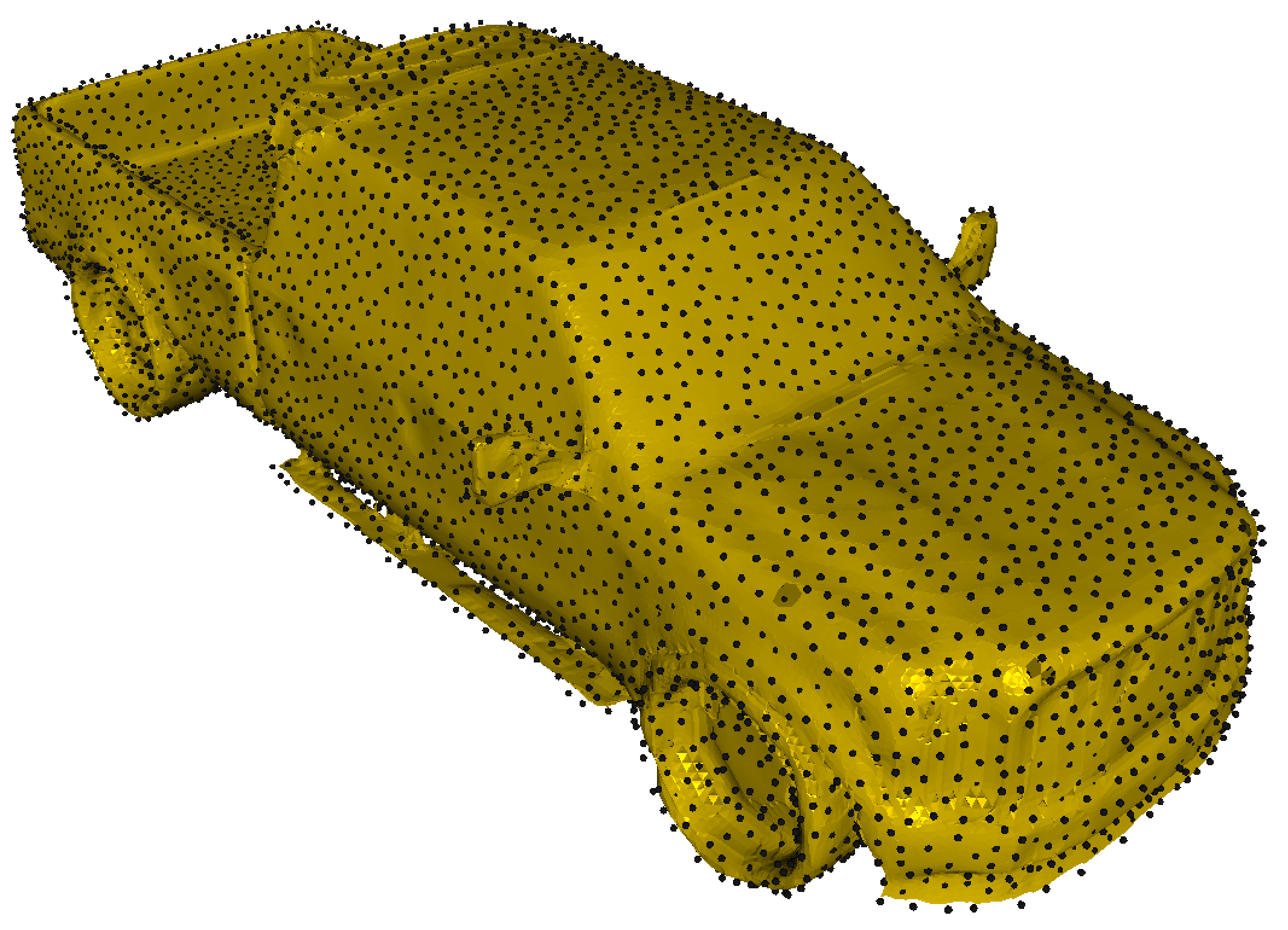}};
		\node[anchor=north west,inner sep=0] (image6) at (image3.south west)
		{\includegraphics[width=0.33\linewidth]{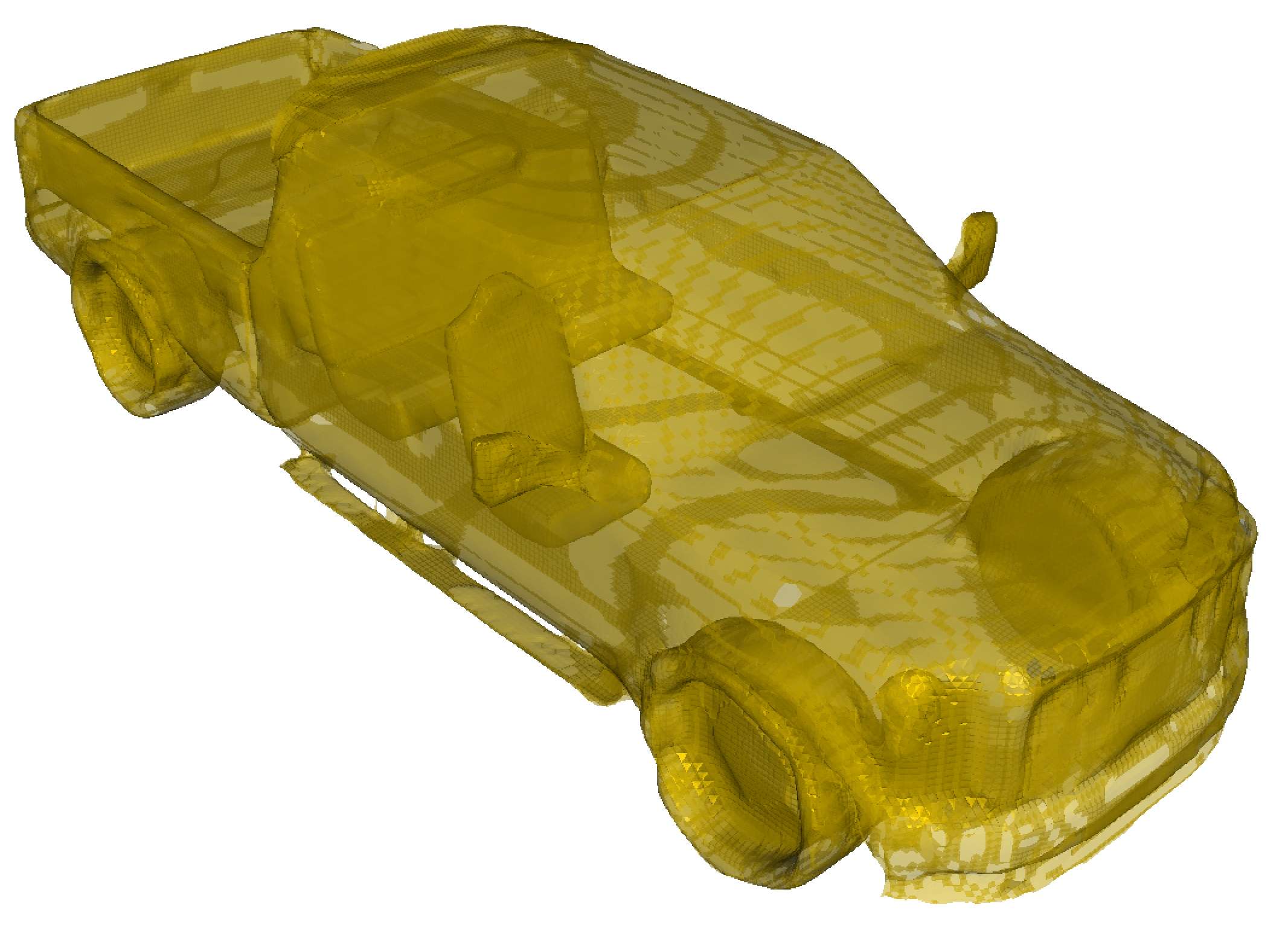}};
		\node[anchor=north,inner sep=1pt] at (image4.south) {(a)};
		\node[anchor=north,inner sep=1pt] at (image5.south) {(b)};
		\node[anchor=north,inner sep=1pt] at (image6.south) {(c)};
	\end{tikzpicture}
	\caption{Example~\ref{ex:surface_recon} for car1 (top) and car2 (bottom): given point cloud (a), reconstructed surface with the point cloud (b), and interior view (c).}
	\label{fig:car}
\end{figure}
\end{example}

\subsection{Model Ablation}\label{sec:ablation}

\subsubsection{Robustness to stochasticity}
Network initialization and collocation sampling are random. Table~\ref{tab:all_errors} reports the mean and standard deviation of the errors over five random seeds, including the gradient error. The standard deviations are small relative to the means for most examples; for Examples~\ref{ex:sine} and~\ref{ex:aniso_ex1} they are comparable to the means, but the errors themselves remain small. Overall, the performance is not highly sensitive to the randomness of network training.

\begin{table}[t]
\centering
\caption{Mean and standard deviation of the errors over five independent runs with different random seeds.}
\label{tab:all_errors}
\small
\setlength{\tabcolsep}{3.5pt}
\begin{tabular}{lccc}
\toprule
\textbf{Example}
& $\boldsymbol{E_2(u_\theta)}$
& $\boldsymbol{E_2(\nabla u_\theta)}$
& $\boldsymbol{E_\infty(u_\theta)}$ \\
\midrule
Sec.~\ref{ex:random_pnts} (rand 1)
& $(8.280 \pm 0.591)\times 10^{-5}$
& $(3.850 \pm 0.235)\times 10^{-2}$
& $(2.150 \pm 0.195)\times 10^{-1}$ \\
Sec.~\ref{ex:random_pnts} (rand 2)
& $(6.347 \pm 0.419)\times 10^{-4}$
& $(3.695 \pm 0.494)\times 10^{-2}$
& $(1.793 \pm 0.457)\times 10^{-1}$ \\
Example~\ref{ex:sine}
& $(3.926 \pm 3.324)\times 10^{-8}$
& $(9.846 \pm 5.766)\times 10^{-6}$
& $(7.115 \pm 2.810)\times 10^{-4}$ \\
Example~\ref{ex:shape_shading}
& $(2.752 \pm 0.181)\times 10^{-4}$
& $(7.335 \pm 0.164)\times 10^{-2}$
& $(4.743 \pm 0.437)\times 10^{-2}$ \\
Example~\ref{ex:aniso_ex1}
& $(2.836 \pm 2.460)\times 10^{-4}$
& $(1.392 \pm 1.888)\times 10^{-3}$
& $(2.566 \pm 1.836)\times 10^{-2}$ \\
Example~\ref{ex:aniso_ex2}
& $(6.635 \pm 0.406)\times 10^{-4}$
& $(1.383 \pm 0.644)\times 10^{-2}$
& $(3.333 \pm 0.856)\times 10^{-2}$ \\
Example~\ref{ex:aniso_ex3}
& $(2.178 \pm 0.064)\times 10^{-3}$
& $(1.600 \pm 0.043)\times 10^{-2}$
& $(1.003 \pm 0.012)\times 10^{-1}$ \\
\bottomrule
\end{tabular}
\end{table}

\subsubsection{Collocation points and minimal viscosity}\label{sec:ablation_collocation}
\begin{figure}[t]
	\centering
	\subfigure[Number $N$ of collocation points]{\includegraphics[width=0.75\linewidth]{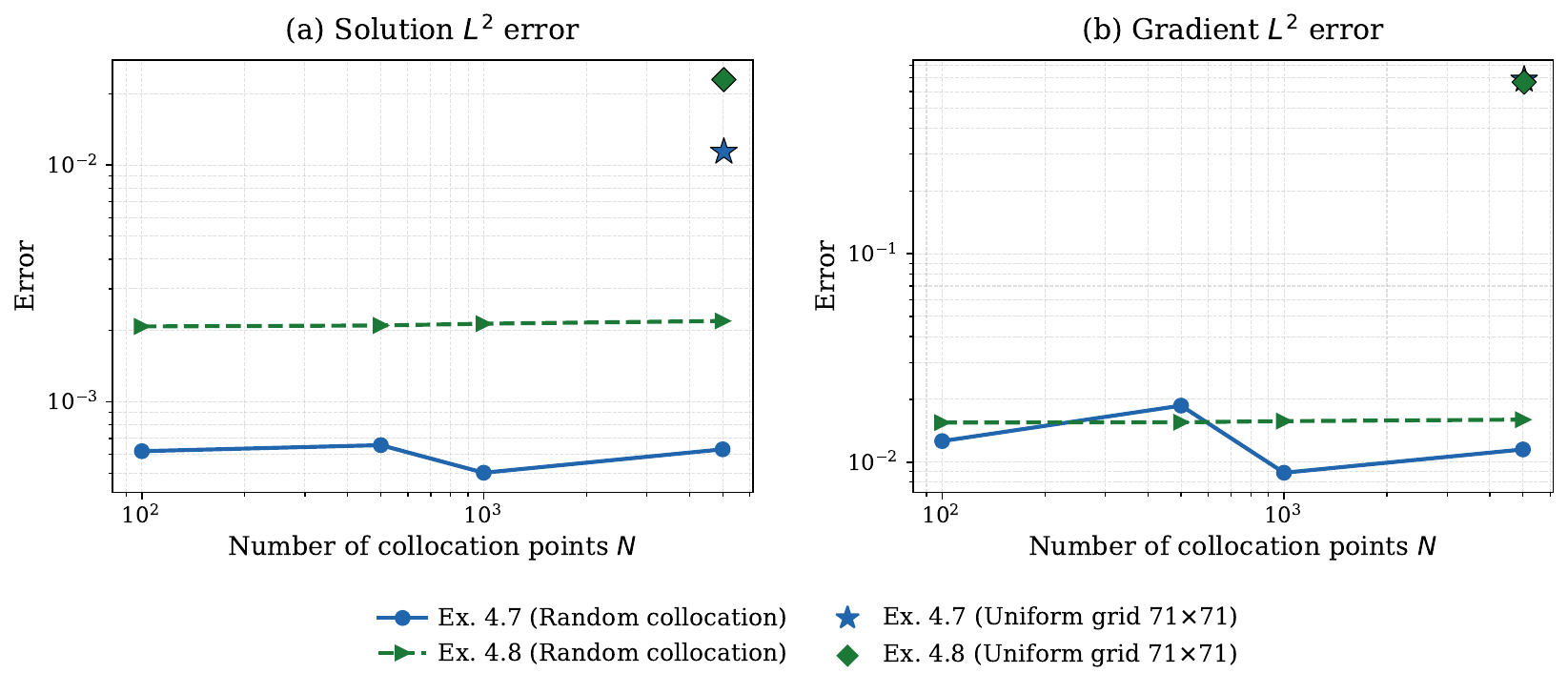}\label{fig:collocation_ablation}}\\
	\subfigure[Minimal viscosity $\epsilon_{\min}$]{\includegraphics[width=0.75\linewidth]{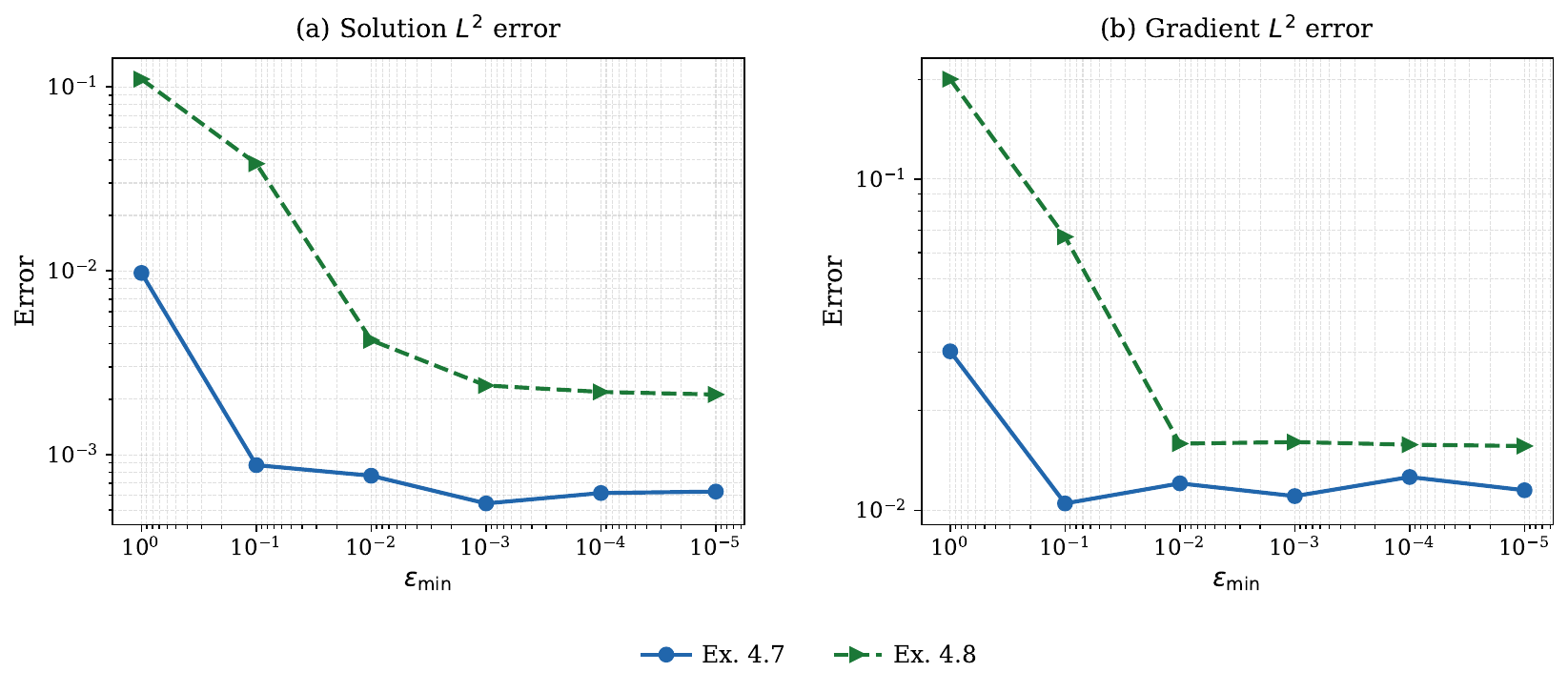}\label{fig:eps_min_ablation}}
	\caption{$E_2(u_\theta)$ (left panel) and $E_2(\nabla u_\theta)$ (right panel) for Examples~\ref{ex:aniso_ex2} (solid blue) and~\ref{ex:aniso_ex3} (dashed green): (a) against $N$ with random resampling, where the blue star and green diamond are the results on a fixed uniform $71\times71$ grid; (b) against $\epsilon_{\min}$.}
	\label{fig:sampling}
\end{figure}

To assess the effects of the number of collocation points and random resampling, we consider the anisotropic Examples~\ref{ex:aniso_ex2} and~\ref{ex:aniso_ex3}, varying the number of points $N\in\{100,500,1000,5000\}$ resampled at each epoch and comparing random sampling with a fixed uniform $71\times71$ grid (5041 points). With random resampling the errors are almost independent of $N$ (Figure~\ref{fig:collocation_ablation}): $E_2(u_\theta)$ stays near $6\times10^{-4}$ for Example~\ref{ex:aniso_ex2} and $2\times10^{-3}$ for Example~\ref{ex:aniso_ex3}, and $E_2(\nabla u_\theta)$ below $2\times10^{-2}$. The fixed grid gives values of $E_2(u_\theta)$ one order of magnitude larger, about $10^{-2}$ and $2\times10^{-2}$, and $E_2(\nabla u_\theta)\approx0.7$ for both examples. As $\epsilon$ decreases, the region where $s_\epsilon^u$ is numerically nonzero shrinks; a fixed grid may place too few points in it, whereas resampling at every epoch keeps reaching it, so random collocation is important for capturing the localized structure of the singular perturbation.

For the same two examples, we investigate how small the minimum viscosity $\epsilon_{\min}=10^{-k_{\max}}$ must be to accurately recover the viscosity solution. Decreasing $\epsilon_{\min}$ over $\{1,10^{-1},\ldots,10^{-5}\}$ reduces the solution and gradient errors most from $1$ to $10^{-1}$ and $10^{-2}$, after which they change little (Figure~\ref{fig:eps_min_ablation}). A sufficiently small $\epsilon_{\min}$ is thus needed, but the method is not sensitive to its precise value, and avoiding unnecessarily small values saves continuation steps.

\subsubsection{Comparison of asymptotic strategies}\label{sec:ablation_lambda}
To test the decoupling of the two limits, we compare three asymptotic strategies:
\begin{enumerate}
	\renewcommand{\labelenumi}{(S\arabic{enumi})}
	\item decrease $\epsilon$ and increase $\lambda$ simultaneously, with $\lambda=1/\epsilon^2$ or $\lambda=1/\epsilon$;
	\item fix $\epsilon=1$ and double $\lambda$ progressively, the strategy suggested by Proposition~\ref{prop:penalty_convergence};
	\item fix $\lambda$ and decrease $\epsilon$, the strategy suggested by Proposition~\ref{prop:bound} and used in all the examples above. 
\end{enumerate}
Figures~\ref{fig:lambda_aniso1} and~\ref{fig:lambda_aniso2} show the three strategies for Examples~\ref{ex:aniso_ex2} and~\ref{ex:aniso_ex3}. Both propositions assume $\Gamma=\partial\Omega$, and Proposition~\ref{prop:bound} also $M=I$ and $g=0$; these anisotropic point-source examples therefore lie beyond their hypotheses. Under (A1), as $\epsilon$ decreases and $\lambda$ increases, the training loss oscillates increasingly, and the two couplings end at different errors. Although $\cL_{\mathrm{GM}}$ decreases on the whole, its fluctuations grow, so the increasing penalty makes the optimization more sensitive and unstable. (A2) is substantially more stable than (A1), and, consistently with the behavior predicted by Proposition~\ref{prop:penalty_convergence}, the stronger penalty drives the solution toward the viscosity solution; for large $\lambda$, however, the total loss increases or fluctuates strongly although $\cL_{\mathrm{GM}}$ decreases, reflecting the deteriorating conditioning. (A3) trains stably and converges toward the viscosity solution with either annealing factor (panels (d) and (e)), consistently with the behavior predicted by Proposition~\ref{prop:bound}. For Example~\ref{ex:aniso_ex2} it is also far more accurate: its $E_2(u_\theta)$ ends below $10^{-3}$, whereas those of (A1) and (A2) end between $2\times10^{-2}$ and $10^{-1}$. For Example~\ref{ex:aniso_ex3} the losses oscillate in the same way under (A1) and the total loss fluctuates under (A2), but the final errors of (A1) and (A2) are comparable to or smaller than those of (A3) in Table~\ref{tab:all_errors}; in particular, (A1) ends with $E_2(\nabla u_\theta)$ between $10^{-3}$ and $1.6\times10^{-3}$, about ten times smaller than the $1.6\times10^{-2}$ of (A3). For this example the advantage of (A3) is the stability of training, not accuracy.
Together with Propositions~\ref{prop:penalty_convergence} and~\ref{prop:bound}, these results indicate that the unit-norm gradient variable permits a stable strategy that avoids controlling the two asymptotic parameters simultaneously, and that, of the two decoupled routes, (A3) trains more stably than (A2).

\section{Conclusion}\label{sec:conclusion}
We computed the viscosity solution of the anisotropic eikonal equation through a variational formulation of its viscous approximation with a unit-norm gradient variable. The unit-norm constraint removes the gradient from the exponentially weighted part of the functional and decouples the viscosity and penalty limits. We proved convergence to the viscosity solution as $\lambda\to\infty$ for each fixed $\epsilon>0$, including strong $L^2$ convergence of the gradient, and, for the distance function to the boundary of a convex or smooth domain, as $\epsilon\to0$ for any fixed $\lambda>0$ with an $O(\epsilon/\lambda)$ error bound. The network satisfies the unit-norm constraint exactly at every point and evaluates the localized functional at sampled points without viscosity-dependent mesh refinement. In isotropic and anisotropic experiments, including point sources beyond the analysis, a fixed $\lambda$ with decreasing $\epsilon$ trains stably, whereas coupling $\lambda$ to $\epsilon$ makes training oscillate. Compared with mesh-based solvers, the method avoids mesh generation; demonstrating this advantage on a domain that requires an unstructured mesh is left to future work. Open problems are the fixed-$\lambda$ analysis for general $M$, $g$ and point sources, and further applications of the framework.

\appendix
\section{Proof of Proposition \ref{prop:penalty_convergence}}\label{appen:proof_fix_epsilon}
\begin{proof}
Fix $\epsilon>0$ and, for each $\lambda>0$, let $(u_\lambda,\mathbf p_\lambda)=(u_{\epsilon,\lambda},\mathbf p_{\epsilon,\lambda})$ be any global minimizer. Here $\mathcal F_\epsilon$ is given by~\eqref{eq:reduced_vis_functional}. Since $M$ is symmetric, continuous and uniformly positive definite on $\overline\Omega$, there is $m_M>0$ with $m_M\|\boldsymbol\xi\|^2\le\|\boldsymbol\xi\|_{M(\bx)}^2=\|M^{1/2}(\bx)\boldsymbol\xi\|^2$ for all $\bx\in\overline\Omega$, $\boldsymbol\xi\in\bR^d$, and $M^{1/2}$ is bounded on $\overline\Omega$. As $\Omega$ is a bounded Lipschitz domain, every function in $W^{1,\infty}(\Omega)$ has a Lipschitz representative on $\overline\Omega$ whose boundary values coincide with its trace, and we
consistently identify each such function with this representative.

\textit{Step 1: uniform bounds.} The Lipschitz function $u^*$ is differentiable a.e.\ and satisfies $\|\nabla u^*\|_M=1$ at every point of differentiability~\cite{evans2022partial}, so $(u^*,\nabla u^*)\in\mathcal A$ and $\mathcal L_{\mathrm{GM}}(u^*,\nabla u^*)=0$. The minimality of $(u_\lambda,\mathbf p_\lambda)$ gives
\begin{equation}\label{eq:objective_inequality_ubar}
\mathcal F_\epsilon(u_\lambda)+\frac\lambda2\|\nabla u_\lambda-\mathbf p_\lambda\|_{L^2(\Omega)}^2\le\mathcal F_\epsilon(u^*)<\infty .
\end{equation}
Since $\mathcal F_\epsilon\ge0$, $\|\nabla u_\lambda-\mathbf p_\lambda\|_{L^2(\Omega)}\le(2\mathcal F_\epsilon(u^*)/\lambda)^{1/2}\to0$ as $\lambda\to\infty$. Moreover, $m_M\|\mathbf p_\lambda\|^2\le\|\mathbf p_\lambda\|_M^2=1$ a.e., so $\|\mathbf p_\lambda\|_{L^2(\Omega)}\le(|\Omega|/m_M)^{1/2}$, and by the triangle inequality $\|\nabla u_\lambda\|_{L^2(\Omega)}\le(2\mathcal F_\epsilon(u^*))^{1/2}+(|\Omega|/m_M)^{1/2}$ for $\lambda\ge1$. As $u_\lambda-u^*\in H_0^1(\Omega)$, the Poincar\'e inequality bounds $\|u_\lambda\|_{H^1(\Omega)}$ uniformly in $\lambda\ge1$.

\textit{Step 2: every limit is a subsolution.} Let $\lambda_n\to\infty$. By Step~1, weak compactness and the Rellich--Kondrachov theorem, a subsequence (not relabelled) satisfies $u_{\lambda_n}\rightharpoonup\bar u$ weakly in $H^1(\Omega)$, $u_{\lambda_n}\to\bar u$ in $L^2(\Omega)$ and a.e.\ in $\Omega$, and $\mathbf p_{\lambda_n}\rightharpoonup\bar{\mathbf p}$ weakly in $L^2(\Omega;\bR^d)$. Since $\nabla u_{\lambda_n}\rightharpoonup\nabla\bar u$ and $\nabla u_{\lambda_n}-\mathbf p_{\lambda_n}\to0$ strongly, $\bar{\mathbf p}=\nabla\bar u$. The set $\{\mathbf q\in L^2(\Omega;\bR^d):\|\mathbf q\|_M\le1\text{ a.e.\ in }\Omega\}$ is convex, because $\|\cdot\|_{M(\bx)}$ is a norm, and closed, because an $L^2$-convergent sequence has an a.e.\ convergent subsequence; hence it is weakly closed. As it contains every $\mathbf p_{\lambda_n}$, $\|\nabla\bar u\|_M\le1$ a.e.\ in $\Omega$; hence $\nabla\bar u\in L^\infty(\Omega;\bR^d)$ and, by the Sobolev embeddings and Morrey's inequality, $\bar u\in W^{1,\infty}(\Omega)$. The trace operator is linear and continuous on $H^1(\Omega)$, hence weakly continuous, so $\bar u=g$ on $\Gamma$. The viscosity solution is the maximal subsolution~\cite{lions1982generalized,ennaji2021augmented}: since $\boldsymbol\xi\mapsto\|\boldsymbol\xi\|_{M(\bx)}$ is convex, every $v\in W^{1,\infty}(\Omega)$ with $\|\nabla v\|_M\le1$ a.e.\ in $\Omega$ is a viscosity subsolution of~\eqref{eq:eikonal}, and if $v=g$ on $\Gamma$ the comparison principle gives $v\le u^*$ in $\overline\Omega$. Therefore
\begin{equation}\label{eq:subsolution_minimizer_u}
\bar u\le u^*\qquad\text{in }\Omega.
\end{equation}

\textit{Step 3: identification of the limit.} By Fatou's lemma and~\eqref{eq:objective_inequality_ubar},
\[
\int_\Omega\bigl(\mathrm e^{-2\bar u/\epsilon}-\mathrm e^{-2u^*/\epsilon}\bigr)\diff\bx=\mathcal F_\epsilon(\bar u)-\mathcal F_\epsilon(u^*)\le\liminf_{n\to\infty}\mathcal F_\epsilon(u_{\lambda_n})-\mathcal F_\epsilon(u^*)\le0 .
\]
By~\eqref{eq:subsolution_minimizer_u}, the integrand is nonnegative, so it vanishes a.e.; since $t\mapsto\mathrm e^{-2t/\epsilon}$ is strictly decreasing, $\bar u=u^*$ a.e.\ in $\Omega$, and $\bar{\mathbf p}=\nabla u^*$.

\textit{Step 4: strong convergence of $\mathbf p$ and $\nabla u$.} Multiplication by $M^{1/2}$ is a bounded linear operator on $L^2(\Omega;\bR^d)$, so $M^{1/2}\mathbf p_{\lambda_n}\rightharpoonup M^{1/2}\nabla u^*$ weakly. Because $\|\mathbf p_{\lambda_n}\|_M=\|\nabla u^*\|_M=1$ a.e., $\|M^{1/2}\mathbf p_{\lambda_n}\|_{L^2(\Omega)}^2=\|M^{1/2}\nabla u^*\|_{L^2(\Omega)}^2=|\Omega|$ for every $n$, so
\[
\|M^{1/2}(\mathbf p_{\lambda_n}-\nabla u^*)\|_{L^2(\Omega)}^2=2|\Omega|-2\bigl(M^{1/2}\mathbf p_{\lambda_n},M^{1/2}\nabla u^*\bigr)_{L^2(\Omega)}\longrightarrow2|\Omega|-2|\Omega|=0 .
\]
Therefore, $m_M\|\mathbf p_{\lambda_n}-\nabla u^*\|_{L^2(\Omega)}^2\le\|M^{1/2}(\mathbf p_{\lambda_n}-\nabla u^*)\|_{L^2(\Omega)}^2\to0$, and
\[
\|\nabla u_{\lambda_n}-\nabla u^*\|_{L^2(\Omega)}\le\|\nabla u_{\lambda_n}-\mathbf p_{\lambda_n}\|_{L^2(\Omega)}+\|\mathbf p_{\lambda_n}-\nabla u^*\|_{L^2(\Omega)}\longrightarrow0 .
\]
Every sequence $\lambda_n\to\infty$ thus has a subsequence along which $u_{\lambda_n}\to u^*$ strongly in $H^1(\Omega)$, that is, $u_{\lambda_n}\to u^*$ and $\nabla u_{\lambda_n}\to\nabla u^*$ in $L^2$, and $\mathbf p_{\lambda_n}\to\nabla u^*$ strongly in $L^2(\Omega;\bR^d)$. The limits do not depend on the subsequence, so these convergences hold as $\lambda\to\infty$.
\end{proof}

\section{Proof of Proposition~\ref{prop:bound}}\label{appen:proof_fix_lambda}

\paragraph{Idea}
The proof relies on a convex relaxation of $\mathcal E_{\epsilon,\lambda}$ and a comparison principle for its minimizer. The relaxation yields a convex variational problem whose Euler--Lagrange equation admits comparison arguments. We first establish that its minimizer coincides with the minimizer of the original functional. We then construct suitable sub- and supersolutions to obtain pointwise bounds relative to the distance function $u^*=\operatorname{dist}(\cdot,\partial\Omega)$. Finally, exploiting the concentration of the source term near the boundary, we derive quantitative estimates for the penalization energy and the gradient.

\paragraph{Setting}
Throughout the proof, we take $M=I$, $g=0$, and $\Gamma=\partial\Omega$, so that $u^*=\operatorname{dist}(\cdot,\partial\Omega)$ \cite{lions1982generalized}. In particular, $u^*\in H_0^1(\Omega)$, $u^*\ge0$, and $\|\nabla u^*\|=1$ a.e. The admissible class is $\mathcal A=\{(u,\mathbf p)\in H_0^1(\Omega)\times L^2(\Omega;\mathbb R^d):\|\mathbf p\|=1\text{ a.e.}\}$.

Assume that $\Delta u^*\le K$ in $\mathcal D'(\Omega)$
for some $K\ge0$. Since $\nabla u^*\in L^\infty(\Omega)$ and nonnegative functions in $C_c^\infty(\Omega)$ are dense in the nonnegative functions in $H_0^1(\Omega)$,
this is equivalent to
\begin{equation}\label{eq:app_K}
\int_\Omega\nabla u^*\cdot\nabla\zeta\diff\bx
\ge -K\int_\Omega\zeta\diff\bx
\qquad
\text{for all }\zeta\in H_0^1(\Omega),\ \zeta\ge0.
\end{equation}
This condition holds with $K=0$ if $\Omega$ is convex, since $u^*$ is then concave, and with $K=(d-1)/r$ under a uniform exterior sphere condition of radius $r$, in particular for $C^{1,1}$ domains. Indeed, if $\mathcal C_r$ denotes the set of centres of exterior tangent balls of radius $r$, then
\[
u^*(\bx)=\inf_{\mathbf a\in\mathcal C_r} \bigl(\|\bx-\mathbf a\|-r\bigr),
\]
and $\Delta\|\bx-\mathbf a\|\le(d-1)/r$ in $\Omega$.
Consequently, $u^*-\frac{K}{2d}\|\bx\|^2$ is superharmonic as a continuous infimum of superharmonic functions. The condition may fail near a re-entrant corner $\mathbf c$ of a planar domain, where $u^*(\bx)=\|\bx-\mathbf c\|$ in a sector and $\Delta u^*(\bx)=1/\|\bx-\mathbf c\|$.

\paragraph{Preliminaries}
We repeatedly use the following properties of Sobolev functions~\cite{gilbarg2001elliptic}: for $v\in H^1(\Omega)$,
\[
v_+\in H^1(\Omega),\qquad \nabla v_+=\mathbf 1_{\{v>0\}}\nabla v,\qquad v_+\in H_0^1(\Omega),
\]
whenever $v\in H_0^1(\Omega)$. Moreover, if $v\in W^{1,1}_{\mathrm{loc}}(\Omega)$, then $\nabla v=0$ a.e.\ on every level set of $v$.

\paragraph{Convex relaxation}
Let $B\subset\mathbb R^d$ denote the closed unit ball, and let $P_B$ be the Euclidean projection onto $B$. Define
\[
\boldsymbol\sigma(\boldsymbol\xi) := \boldsymbol\xi- P_B\boldsymbol\xi = (\|\boldsymbol\xi\|-1)_+ \frac{\boldsymbol\xi}{\|\boldsymbol\xi\|}, 
\qquad \boldsymbol\sigma(\mathbf 0):=\mathbf 0.
\]
This map is the gradient of the convex function $\frac12\operatorname{dist}(\cdot,B)^2$ and satisfies $\|\boldsymbol\sigma(\boldsymbol\xi)\|\le\|\boldsymbol\xi\|$. Moreover, the firm nonexpansiveness of $P_B$ yields 
\begin{equation}\label{eq:app_firm}
\bigl(\boldsymbol\sigma(\mathbf a)-\boldsymbol\sigma(\mathbf b)\bigr)\cdot(\mathbf a-\mathbf b)-\|\boldsymbol\sigma(\mathbf a)-\boldsymbol\sigma(\mathbf b)\|^2=(P_B\mathbf a-P_B\mathbf b)\cdot(\mathbf a-\mathbf b)-\|P_B\mathbf a-P_B\mathbf b\|^2\ \ge\ 0 .
\end{equation}
for all $\mathbf a,\mathbf b\in\mathbb R^d$. We consider the convex relaxation
\[
\mathcal G_{\epsilon,\lambda}(u) := \mathcal F_\epsilon(u) + \frac\lambda2\int_\Omega \operatorname{dist}(\nabla u,B)^2\diff\bx, \qquad u\in H_0^1(\Omega).
\]

\begin{lemma}[Convex relaxation]\label{lem:app_relax}
Let $\epsilon,\lambda>0$. (i) $\mathcal G_{\epsilon,\lambda}$ has a unique minimizer $u$; it satisfies $u\ge0$ and
\begin{equation}\label{eq:app_EL}
\lambda\int_\Omega\boldsymbol\sigma(\nabla u)\cdot\nabla\varphi\diff\bx=\frac2\epsilon\int_\Omega s_\epsilon^{u}\varphi\diff\bx\qquad\text{for all }\varphi\in H_0^1(\Omega).
\end{equation}
(ii) $\boldsymbol\sigma(\nabla u)\in H^1_{\mathrm{loc}}(\Omega;\mathbb R^d)$. (iii) $\|\nabla u\|>1$ a.e.\ in $\Omega$. (iv) $(u,\nabla u/\|\nabla u\|)$ is the unique global minimizer of $\mathcal E_{\epsilon,\lambda}$ on $\mathcal A$.
\end{lemma}

\begin{proof}
(i) The functional $\mathcal G_{\epsilon,\lambda}$ is convex, finite at $0$, lower semicontinuous on $H_0^1(\Omega)$. Indeed, lower semicontinuity of $\mathcal F_\epsilon$ follows from Fatou's lemma along subsequences converging a.e., while the second term is continuous, since $|\operatorname{dist}(\mathbf a,B)^2-\operatorname{dist}(\mathbf b,B)^2|\le(\|\mathbf a\|+\|\mathbf b\|)\|\mathbf a-\mathbf b\|$). Thus, $\mathcal G_{\epsilon,\lambda}$ is also weakly lower
semicontinuous. Moreover, it is coercive by
$\operatorname{dist}(\boldsymbol\xi,B)^2\ge\frac12\|\boldsymbol\xi\|^2-1$ and the Poincar\'e inequality. A minimizer $u$ therefore exists, and it is unique because $t\mapsto\mathrm e^{-2t/\epsilon}$ is strictly convex. Since $s_\epsilon^{u_+}\le s_\epsilon^u$ and $\operatorname{dist}(\nabla u_+,B)\le\operatorname{dist}(\nabla u,B)$, $u_+$ is also a minimizer. Uniqueness therefore gives
$u=u_+\ge0$. 

For $\varphi\in C_c^\infty(\Omega)$ and $|t|\le1$, we have 
\[
s_\epsilon^{u+t\varphi}\le\mathrm e^{2\|\varphi\|_\infty/\epsilon},\quad  |\boldsymbol\sigma(\nabla u+t\nabla\varphi)\cdot\nabla\varphi|\le(\|\nabla u\|+\|\nabla\varphi\|)\|\nabla\varphi\|.
\]
Differentiation under the integral is therefore justified
by dominated convergence, and the resulting stationarity
condition is~\eqref{eq:app_EL}. Since both sides of
\eqref{eq:app_EL} are continuous in
$\varphi\in H_0^1(\Omega)$, the identity extends to all
such $\varphi$.

(ii) We show that
$\boldsymbol\sigma(\nabla u)\in H^1_{\mathrm{loc}}
(\Omega;\mathbb R^d)$ using difference quotients~\cite{evans2022partial}. Let $\omega\Subset\Omega$ and choose $\chi\in C_c^\infty(\Omega)$ with $0\le\chi\le1$, $\chi=1$ on $\omega$. For $0<\|\mathbf h\|<\operatorname{dist}(\operatorname{supp}\chi,\partial\Omega)$, set $\delta u=u(\cdot+\mathbf h)-u$ and $\delta\boldsymbol\sigma=\boldsymbol\sigma(\nabla u(\cdot+\mathbf h))-\boldsymbol\sigma(\nabla u)$. Testing~\eqref{eq:app_EL} with $(\chi^2\delta u)(\cdot-\mathbf h)$ and with $\chi^2\delta u$, and subtracting the resulting identities, gives 
\[
\frac2\epsilon\int(s_\epsilon^{u(\cdot+\mathbf h)}-s_\epsilon^u)\chi^2\delta u\diff\bx\le0.
\]
Using $\nabla(\chi^2\delta u)=\chi^2\nabla\delta u+2\chi\,\delta u\,\nabla\chi$ and \eqref{eq:app_firm}, we obtain $\|\chi\delta\boldsymbol\sigma\|_{L^2}^2\le2\|\chi\delta\boldsymbol\sigma\|_{L^2}\|\delta u\nabla\chi\|_{L^2}$. Consequently, 
\[
\|\delta\boldsymbol\sigma\|_{L^2(\omega)}\le2\|\nabla\chi\|_\infty\|\mathbf h\|\|\nabla u\|_{L^2(\Omega)},
\]
which yields the claimed local $H^1$ regularity~\cite{evans2022partial}.

(iii) By (ii) and~\eqref{eq:app_EL}, $\nabla\cdot\boldsymbol\sigma(\nabla u)=-\frac{2}{\lambda\epsilon}s_\epsilon^u<0$ a.e. On $Z=\{\|\nabla u\|\le1\}$ every component of $\boldsymbol\sigma(\nabla u)\in H^1_{\mathrm{loc}}$ vanishes. Its weak gradient
therefore vanishes a.e. on $Z$, by the level-set property
of Sobolev functions. Hence, the divergence also vanishes
a.e. on $Z$, which implies $|Z|=0$.

(iv) For any $(v,\mathbf q)\in\mathcal A$, $\|\nabla v-\mathbf q\|\ge|\|\nabla v\|-1|\ge\operatorname{dist}(\nabla v,B)$, and thus 
\[
\mathcal E_{\epsilon,\lambda}(v,\mathbf q)\ge\mathcal G_{\epsilon,\lambda}(v)\ge\mathcal G_{\epsilon,\lambda}(u).
\]
By (iii), $\mathbf p=\nabla u/\|\nabla u\|$ is defined a.e.\ and $\|\nabla u-\mathbf p\|=\|\nabla u\|-1=\operatorname{dist}(\nabla u,B)$.
Therefore, $\mathcal E_{\epsilon,\lambda}(u,\mathbf p)=\mathcal G_{\epsilon,\lambda}(u)$,
so $(u,\mathbf p)$ is another global minimizer, then $v=u$ by uniqueness in (i). Moreover, equality in the first inequality above requires $\nabla u\cdot\mathbf q=\|\nabla u\|$, which, since $\|\mathbf q\|=1$, gives $\mathbf q=\mathbf p$ a.e.
\end{proof}

\begin{lemma}[Weak comparison]\label{lem:app_comparison}
Let $u$ be as in Lemma~\ref{lem:app_relax} and $v\in H_0^1(\Omega)$ with $v\ge0$. If $v$ is a supersolution, $\lambda\int_\Omega\boldsymbol\sigma(\nabla v)\cdot\nabla\varphi\diff\bx\ge\frac2\epsilon\int_\Omega s_\epsilon^v\varphi\diff\bx$ for all $\varphi\in H_0^1(\Omega)$ with $\varphi\ge0$, then $u\le v$ a.e. If $v$ is a subsolution (reverse inequality), then $v\le u$ a.e.
\end{lemma}

\begin{proof}
Testing~\eqref{eq:app_EL} and the supersolution inequality with $\varphi=(u-v)_+$, and subtracting, gives
\[
0\le\lambda\int_{\{u>v\}}(\boldsymbol\sigma(\nabla u)-\boldsymbol\sigma(\nabla v))\cdot\nabla(u-v)\diff\bx\le\frac2\epsilon\int_\Omega(s_\epsilon^u-s_\epsilon^v)(u-v)_+\diff\bx,
\]
where the first inequality follows from~\eqref{eq:app_firm}. Since $s_\epsilon^u<s_\epsilon^v$ on $\{u>v\}$, the right-hand side is nonpositive and is strictly negative if $|\{u>v\}|>0$. Hence $u\le v$ a.e. The subsolution case follows analogously by testing with $\varphi=(v-u)_+$, which gives $v\le u$ a.e.
\end{proof}

\begin{proof}[Proof of Proposition~\ref{prop:bound}]
Write $(u,\mathbf p)=(u_{\epsilon,\lambda},\mathbf p_{\epsilon,\lambda})$. By Lemma~\ref{lem:app_relax}, the global minimizer exists and is unique, $u$ minimizes $\mathcal G_{\epsilon,\lambda}$, and $\mathbf p=\nabla u/\|\nabla u\|$ a.e. Moreover, setting $\boldsymbol\sigma:=\boldsymbol\sigma(\nabla u) = \nabla u-\mathbf p$, we have
\[
\|\nabla u\|=1+\|\boldsymbol\sigma\|\quad\text{a.e.},\qquad 
\mathcal L_{\mathrm{GM}}(u,\mathbf p) = \int_\Omega\|\boldsymbol\sigma\|^2\diff\bx.
\]

\textit{Lower bound.} Since $\|\nabla u^*\|=1$ a.e., $\boldsymbol\sigma(\nabla u^*)=\mathbf 0$ a.e., so $u^*$ is a subsolution, and Lemma~\ref{lem:app_comparison} gives $u^*\le u$ a.e.

\textit{Upper bound.} Let $0<\epsilon<2/K$ and define $c=(1-K\epsilon/2)^{-1}$ and $\psi(t)=\frac{c\epsilon}{2\lambda}(1-\mathrm e^{-2t/\epsilon})$. Then $\lambda\psi'(t)=c\,\mathrm e^{-2t/\epsilon}$ and $\psi''=-\frac2\epsilon\psi'$. 
As $\psi\in C^\infty([0,\infty))$ with $\psi(0)=0$ and $\psi',\psi''$ bounded, and $u^*\in W^{1,\infty}(\Omega)\cap H_0^1(\Omega)$, both $u_\psi=u^*+\psi(u^*)$ and, for $\varphi\in H_0^1(\Omega)$ with $\varphi\ge0$, $\zeta=\psi'(u^*)\varphi\ge0$ lie in $H_0^1(\Omega)$; also $u_\psi\ge u^*\ge0$ and $\boldsymbol\sigma(\nabla u_\psi)=\psi'(u^*)\nabla u^*$ a.e. Since $\nabla\zeta=\psi''(u^*)\varphi\nabla u^*+\psi'(u^*)\nabla\varphi$ and $\|\nabla u^*\|=1$, \eqref{eq:app_K} gives
\[
\begin{aligned}
\lambda\int_\Omega\boldsymbol\sigma(\nabla u_\psi)\cdot\nabla\varphi\diff\bx&=\lambda\int_\Omega\bigl(\nabla u^*\cdot\nabla\zeta-\psi''(u^*)\varphi\bigr)\diff\bx\ge\lambda\Bigl(\frac2\epsilon-K\Bigr)\int_\Omega\psi'(u^*)\varphi\diff\bx\\
&=\frac2\epsilon\int_\Omega s_\epsilon^{u^*}\varphi\diff\bx\ge\frac2\epsilon\int_\Omega s_\epsilon^{u_\psi}\varphi\diff\bx ,
\end{aligned}
\]
by $c(\frac2\epsilon-K)=\frac2\epsilon$ and $u_\psi\ge u^*$. Hence, $u_\psi$ is a supersolution.  By Lemma~\ref{lem:app_comparison}, $u\le u_\psi$ a.e. Combining this with the lower bound gives
$0\le u-u^*\le\psi(u^*)\le\frac{\epsilon}{2\lambda(1-K\epsilon/2)}$ a.e. Finally, if $\epsilon\le1/K$, then $(1-K\epsilon/2)^{-1}\le1+K\epsilon$, which yields~\eqref{eq:err_bound_fix_lambda}.

\textit{Boundary layer estimate.} Let $C_\Omega=\sup_{t>0}|\{u^*<t\}|/t$. Then, for every $\epsilon>0$,
\[
\mathcal F_\epsilon(u^*)=\int_0^\infty\frac2\epsilon\mathrm e^{-2t/\epsilon}\,|\{u^*<t\}|\diff t\le C_\Omega\int_0^\infty\frac{2t}\epsilon\mathrm e^{-2t/\epsilon}\diff t=\frac{C_\Omega}2\,\epsilon .
\]
We next show that $C_\Omega<\infty$. Since $\Omega$ is Lipschitz, its boundary can be covered by finitely many cylinders $C_i$ in which, after a rotation of coordinates, $\Omega\cap C_i=\{x_d>\gamma_i(\bx')\}$, with $\gamma_i$ being $L$-Lipschitz on the base $C_i'$. Choose $\rho_0>0$ such that every ball of radius $\rho_0$ centred on $\partial\Omega$ is contained in one of these cylinders. If $u^*(\bx)<t<\rho_0$ and $\mathbf y\in\partial\Omega$ is the nearest boundary point to $\bx$, then $\bx,\mathbf y\in C_i$ for some $i$, and $0<x_d-\gamma_i(\bx')\le|x_d-y_d|+L\|\bx'-\mathbf y'\|<\sqrt{1+L^2}\,t$.
Consequently, $|\{u^*<t\}|\le\sqrt{1+L^2}\,t\sum_i|C_i'|$. For $t\ge\rho_0$, we simply have $|\{u^*< t\}|/t\le|\Omega|/\rho_0$. These estimates establish $C_\Omega<\infty$.

\textit{Bounds on $\mathcal L_{\mathrm{GM}}$ and the gradient.} Testing~\eqref{eq:app_EL} with $\varphi=u$ and using $\boldsymbol\sigma\cdot\nabla u = \|\boldsymbol\sigma\|(1+\|\boldsymbol\sigma\|)$, we obtain
\[
\lambda\int_\Omega \bigl(\|\boldsymbol\sigma\|+\|\boldsymbol\sigma\|^2\bigr) \diff\bx = \frac2\epsilon\int_\Omega u\,\mathrm e^{-2u/\epsilon}\diff\bx.
\]
Since $t\mathrm e^{-2t/\epsilon} \le \frac{\epsilon}{\mathrm e}\mathrm e^{-t/\epsilon}$ for $t\ge0$ and $u\ge u^*$, it follows that
\[
\lambda\int_\Omega \bigl(\|\boldsymbol\sigma\|+\|\boldsymbol\sigma\|^2\bigr) \diff\bx\le\frac2{\mathrm e} \int_\Omega\mathrm e^{-u^*/\epsilon}\diff\bx=\frac2{\mathrm e}\mathcal F_{2\epsilon}(u^*) \le\frac{2C_\Omega}{\mathrm e}\epsilon,
\]
where we used the preceding estimate with $2\epsilon$
in place of $\epsilon$. In particular, $\mathcal L_{\mathrm{GM}}(u,\mathbf p)=\int_\Omega\|\boldsymbol\sigma\|^2\diff\bx\le\frac{2C_\Omega}{\mathrm e}\frac\epsilon\lambda$ for all $\epsilon,\lambda>0$. Neither this estimate nor the lower bound requires~\eqref{eq:app_K}. Together with Lemma~\ref{lem:app_relax}, this proves~(i). 

To estimate the gradient, let $\zeta=u-u^*\in H_0^1(\Omega)$, so that $\zeta\ge0$. Using $\|\nabla u^*\|=1$ and $\|\nabla u\|^2=1+2\|\boldsymbol\sigma\|+\|\boldsymbol\sigma\|^2$, we have  $\|\nabla u-\nabla u^*\|^2=2\|\boldsymbol\sigma\|+\|\boldsymbol\sigma\|^2-2\nabla\zeta\cdot\nabla u^*$ a.e. 
By~\eqref{eq:app_K} and~\eqref{eq:err_bound_fix_lambda}, $-2\int_\Omega\nabla\zeta\cdot\nabla u^*\diff\bx\le2K\int_\Omega\zeta\diff\bx\le K|\Omega|(1+K\epsilon)\frac\epsilon\lambda$.
Therefore, for $0<\epsilon\le1/K$,
\[
\|\nabla u-\nabla u^*\|_{L^2(\Omega)}^2\le2\int_\Omega\bigl(\|\boldsymbol\sigma\|+\|\boldsymbol\sigma\|^2\bigr)\diff\bx+2K\int_\Omega\zeta\diff\bx\le\Bigl(\frac{4C_\Omega}{\mathrm e}+2K|\Omega|\Bigr)\frac\epsilon\lambda .
\]
Finally, $\nabla u=(1+\|\boldsymbol\sigma\|)\mathbf p$ gives $\|\nabla u-\nabla u^*\|^2=\|\boldsymbol\sigma\|^2+(1+\|\boldsymbol\sigma\|)\|\mathbf p-\nabla u^*\|^2\ge\|\mathbf p-\nabla u^*\|^2$ a.e., which proves~\eqref{eq:grad_bound_fix_lambda}.
\end{proof}

While the proposed deep learning methodology involves additional sources of error, including neural network approximation, generalization from the empirical objective to its population counterpart, and optimization errors, these aspects are beyond the scope of the present analysis. Instead, our results establish the convergence of the underlying variational problem to the desired viscosity solution through a suitable separation of the viscosity and penalty parameters. In this way, we provide a theoretical justification for the proposed variational
formulation and establish a mathematical foundation for the proposed methodology, independently of the additional errors introduced by its neural network implementation.
A rigorous convergence analysis of the proposed deep learning methodology is left for future work.

\section{Further training curves}
\begin{figure}[ht]
	\centering
	\subfigure[Automatic annealing with factor $1/2$]{\includegraphics[width=0.99\linewidth]{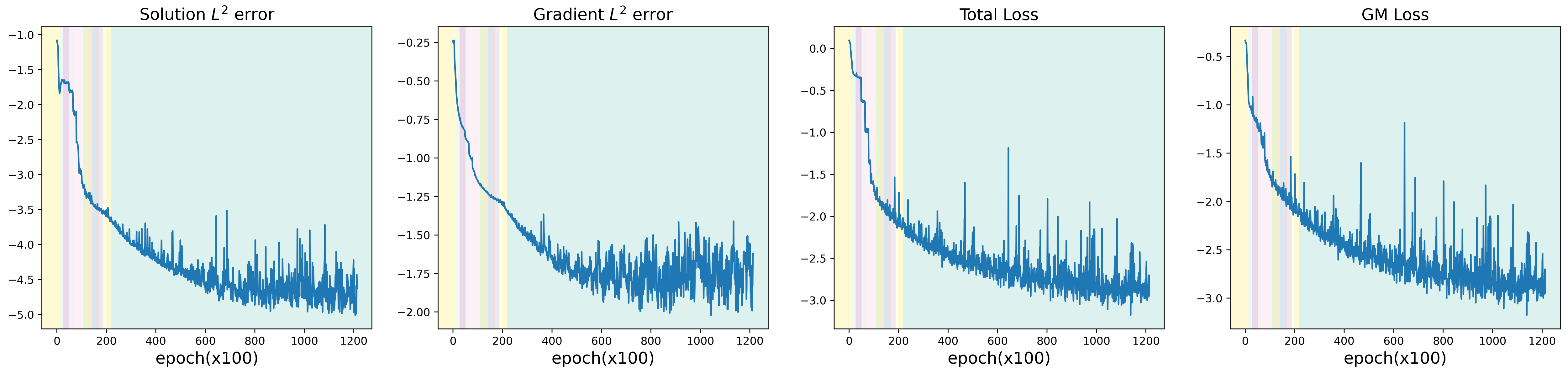}}\\
	\subfigure[Automatic annealing with factor $1/10$]{\includegraphics[width=0.99\linewidth]{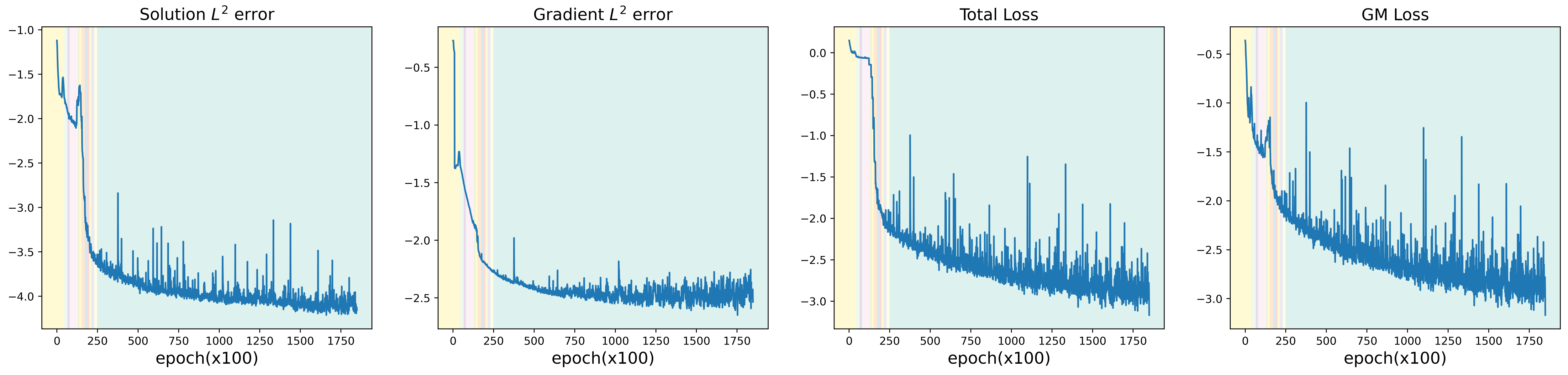}}
	\caption{Errors and losses during training for the second set of random source points in Section~\ref{ex:random_pnts} (Figure~\ref{fig:rand}(b)).}
	\label{fig:rand_eps_2}
\end{figure}

\begin{figure}[ht]
	\centering
	\subfigure[Example~\ref{ex:sine}]{\includegraphics[width=0.99\linewidth]{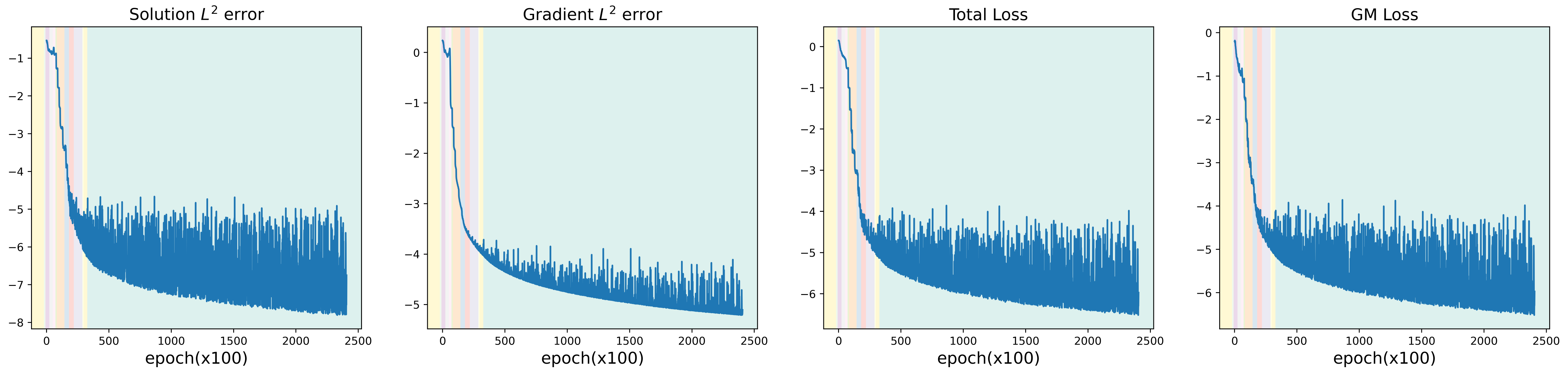}}\\
	\subfigure[Example~\ref{ex:shape_shading}]{\includegraphics[width=0.99\linewidth]{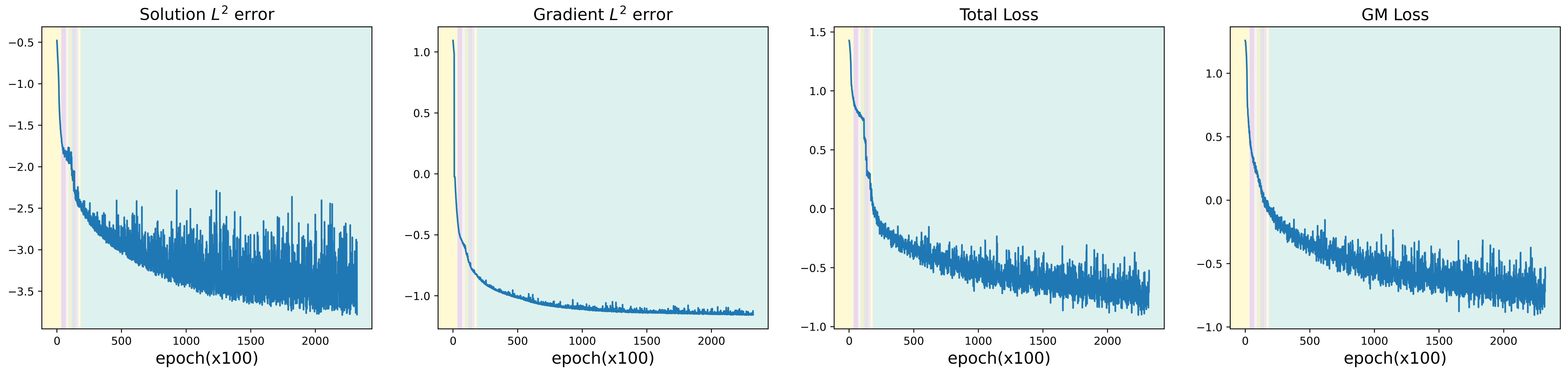}}\\
	\subfigure[Example~\ref{ex:aniso_ex1}]{\includegraphics[width=0.99\linewidth]{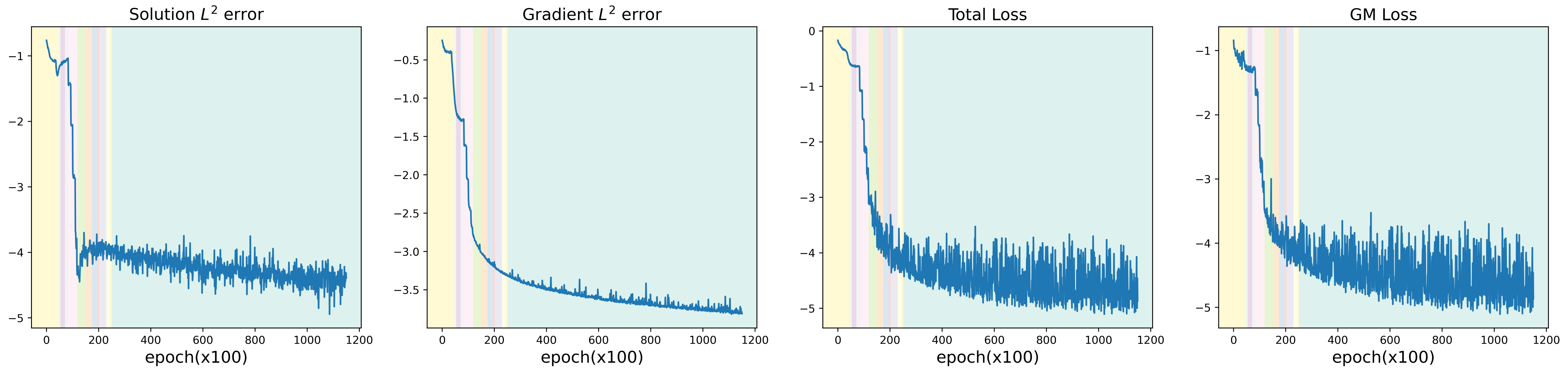}}
	\caption{Errors and losses during training with the annealing factor $1/2$ (for the factor $1/10$, see Figure~\ref{fig:curves_exact}).}
	\label{fig:half_exact}
\end{figure}

\begin{figure}[ht]
	\centering
	\subfigure[Example~\ref{ex:disconti1}]{\includegraphics[width=0.48\linewidth]{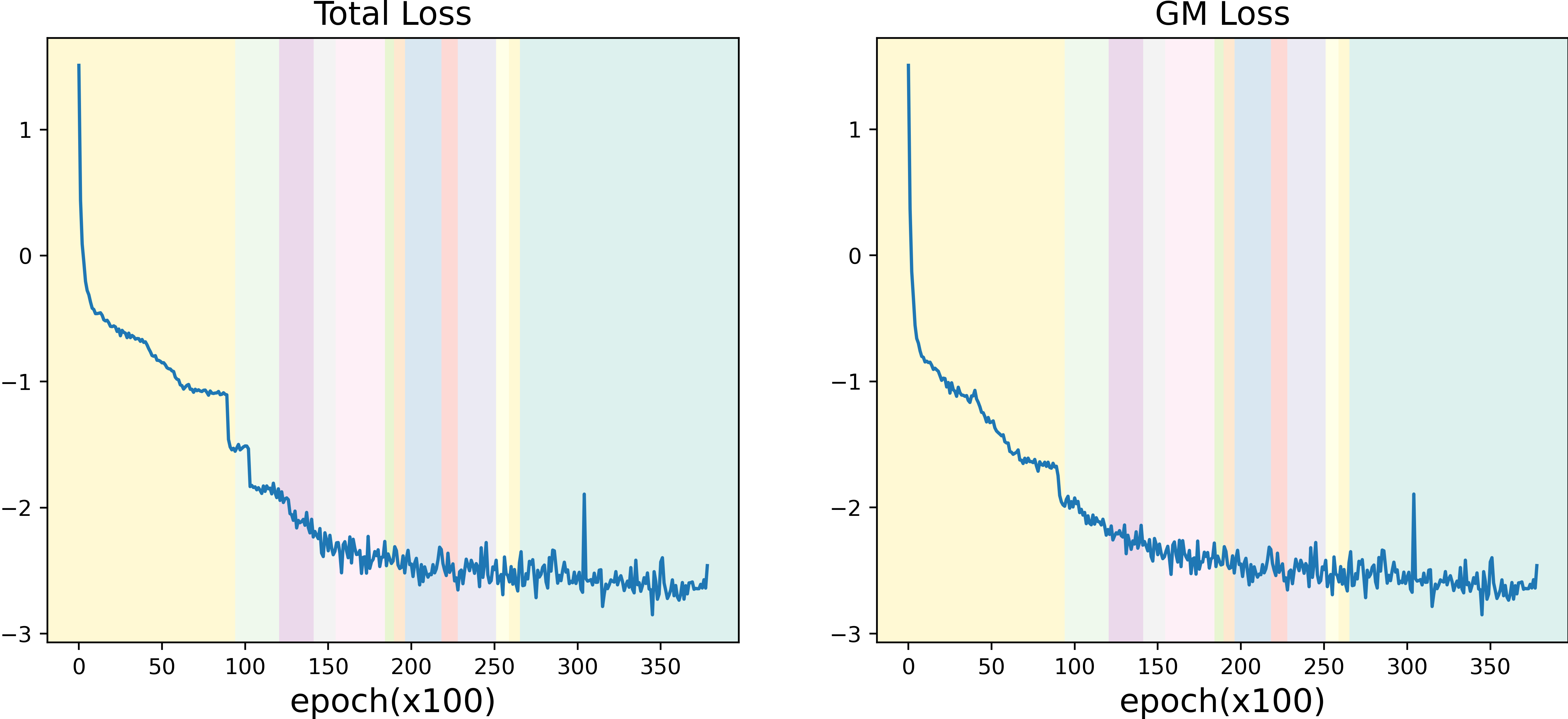}}\hfill
	\subfigure[Example~\ref{ex:sinusoidal}]{\includegraphics[width=0.48\linewidth]{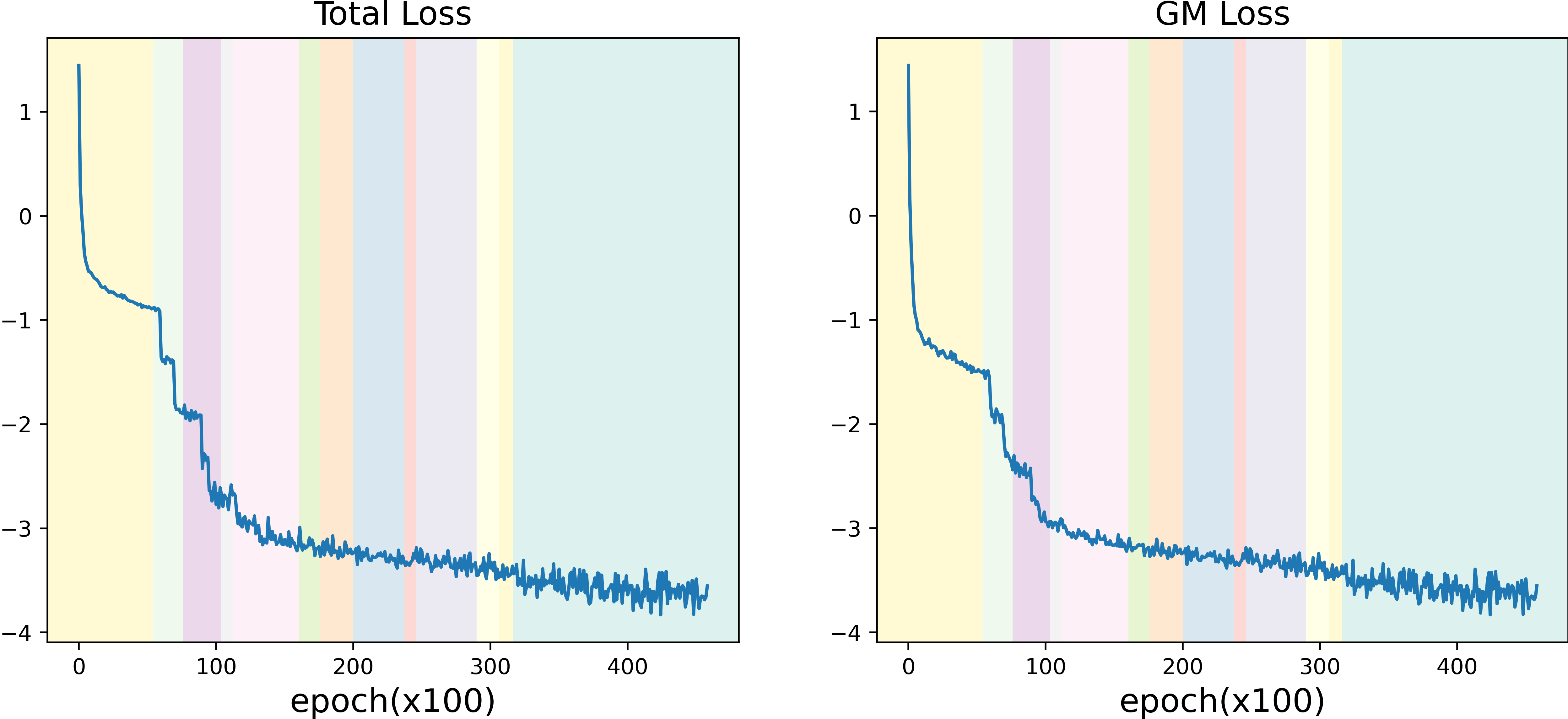}}\\
	\subfigure[Example~\ref{ex:single_source}]{\includegraphics[width=0.48\linewidth]{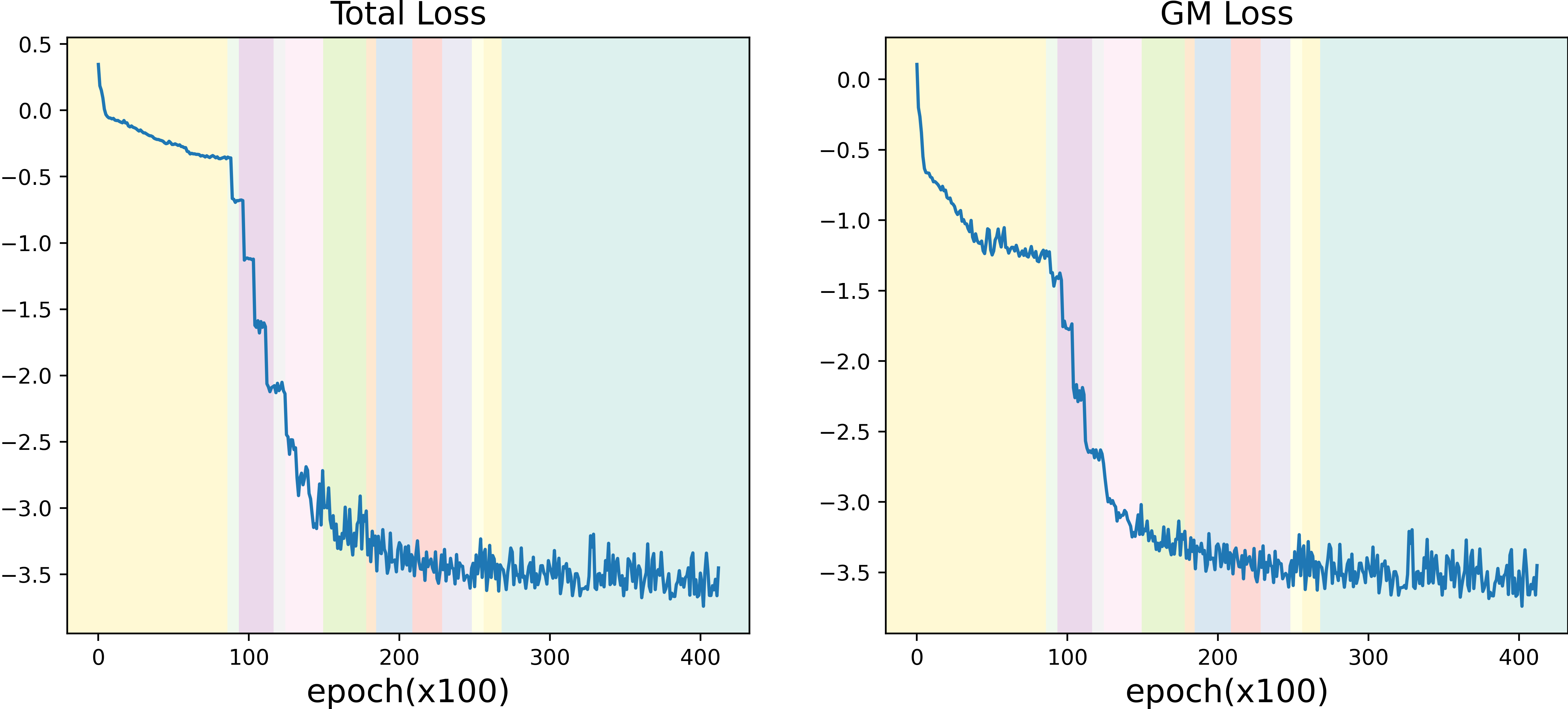}}\hfill
	\subfigure[Example~\ref{ex:aniso_ex3_rectangle}]{\includegraphics[width=0.48\linewidth]{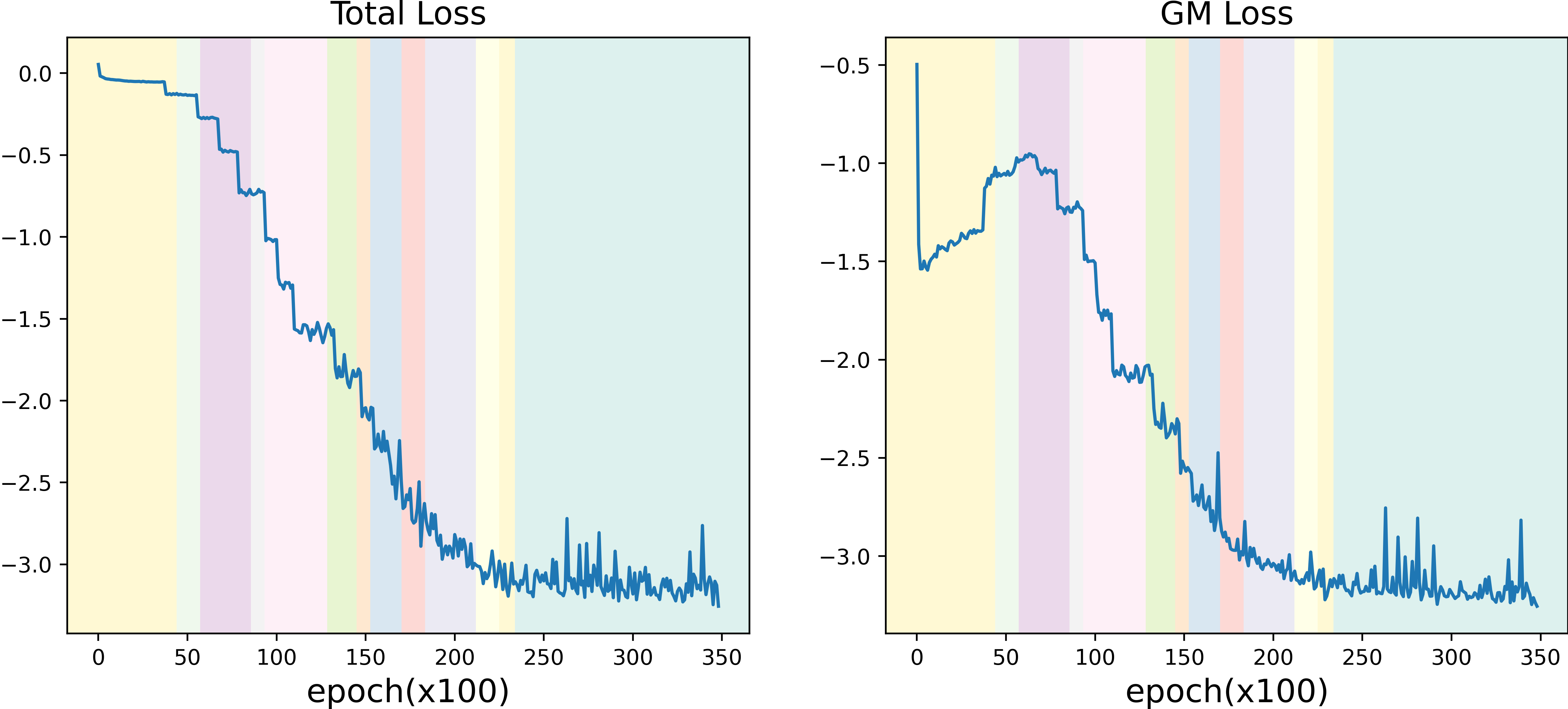}}
	\caption{Total loss and GM loss during training with the annealing factor $1/2$ for the examples without an exact solution (for the factor $1/10$, see Figure~\ref{fig:loss_curves}).}
	\label{fig:half_loss}
\end{figure}

The main text shows the training curves with the annealing factor $1/10$ of $\epsilon$ for all examples (Figures~\ref{fig:rand_eps}, \ref{fig:curves_exact}, \ref{fig:loss_curves}, \ref{fig:lambda_aniso1} and~\ref{fig:lambda_aniso2}), and with the factor $1/2$ for the first set of random source points of Section~\ref{ex:random_pnts} and for Examples~\ref{ex:aniso_ex2} and~\ref{ex:aniso_ex3} (Figures~\ref{fig:rand_eps}(a), \ref{fig:lambda_aniso1}(e) and~\ref{fig:lambda_aniso2}(e)). We give the remaining curves: the mean squared errors $E_2(u_\theta)$ and $E_2(\nabla u_\theta)$ of Section~\ref{sec:experiments} (labelled ``$L^2$ error''), the total loss and the GM loss for the second set of random source points (Figure~\ref{fig:rand}(b)) with both factors (Figure~\ref{fig:rand_eps_2}) and for Examples~\ref{ex:sine}, \ref{ex:shape_shading} and~\ref{ex:aniso_ex1} with the factor $1/2$ (Figure~\ref{fig:half_exact}), and the total loss and the GM loss for Examples~\ref{ex:disconti1}, \ref{ex:sinusoidal}, \ref{ex:single_source} and~\ref{ex:aniso_ex3_rectangle}, without an exact solution, with the factor $1/2$ (Figure~\ref{fig:half_loss}). In all graphs the $y$-axis is on a logarithmic scale and the background colors mark the successive values of $\epsilon$. For the examples with an exact or reference solution, the errors decrease as $\epsilon$ decreases, alike for the two schedules.

\section*{Acknowledgments}
The authors used the AI assistant Claude (Anthropic) to suggest edits to the wording and presentation of the manuscript, to check its internal consistency, and to help locate and check related literature. All suggestions were reviewed by the authors, and every cited work was checked by the authors. The authors assume responsibility for all content.

\bibliographystyle{siamplain}
\bibliography{references}








\end{document}